\documentclass[final,3p]{elsarticle}
\usepackage[latin1]{inputenc}
 \usepackage{graphics}
 \usepackage{graphicx}
 \usepackage{epsfig}
\usepackage{amssymb}
 \usepackage{amsthm}
 \usepackage{lineno}
 \usepackage{amsmath}
   \numberwithin{equation}{section}
\usepackage{mathrsfs}

\newtheorem{thm}{Theorem}[section]
\newtheorem{cor}[thm]{Corollary}
\newtheorem{lem}[thm]{Lemma}

\newtheorem{defn}[thm]{Definition}

\biboptions{sort&compress,square}
\begin{document}
\begin{frontmatter}
\author{Tong Wu}
\ead{wut977@nenu.edu.cn}
\author{Yong Wang\corref{cor2}}
\ead{wangy581@nenu.edu.cn}
\cortext[cor2]{Corresponding author.}

\address{School of Mathematics and Statistics, Northeast Normal University,
Changchun, 130024, China}

\title{Affine Ricci Solitons associated to the Bott connection on three-dimensional Lorentzian Lie groups}
\begin{abstract}
In this paper, we compute the Bott connections and their curvature on three-dimensional Lorentzian Lie groups with three different distributions, and we classify affine Ricci solitons associated to the Bott connection on three-dimensional Lorentzian Lie groups with three different distributions.
\end{abstract}
\end{frontmatter}
\section{Introduction}
\indent In geometry, Einstein metrics have been studied in three-dimensional Lorentzian manifolds. As a natural generalization of Einstein metrics, the definition of the Ricci soliton was introduced by Hamilton in \cite{RS}. Therefore, study of Ricci solitons over different geometric spaces has become more significant. Einstein manifolds associated to affine connections have been studied by many geometers. In \cite {YH} and \cite{YY1}, Wang studied Einstein manifolds associated to semi-symmetric non-metric connections and  semi-symmetric metric connections respectively.\\
 \indent Naturally, mathematicians start to study Ricci solitons associated to different affine connections. In \cite{M},  Crasmareanu gave several generalizations of Ricci solitons with linear connections. In \cite{N}, it is proved that the equivalent conditions of the Ricci Soliton in Kenmotsu manifold associated to the Schouten-Van Kampen connection to be steady. Hui-Prasad-Chakraborty studied Ricci solitons on Kenmotsu manifolds with respect to quarter symmetric non-metric $\phi$-connection in \cite{S1}. Perktas-Yildiz studied several soliton types on a quasi-Sasakian 3-manifold with respect to the Schouten-Van Kampen connection in \cite{S2}. In \cite{YY2}, Wang classify affine Ricci solitons associated to canonical connections and Kobayashi-Nomizu connections and perturbed canonical connections and perturbed Kobayashi-Nomizu connections on three-dimensional Lorentzian Lie groups with some product structures. In \cite{F,J,RK}, the definition of the Bott connection was introduced. In this paper, we compute the Bott connections and their curvature on three-dimensional Lorentzian Lie groups with three different distributions, and we classify affine Ricci solitons associated to the Bott connection on three-dimensional Lorentzian Lie groups with three different distributions.\\
\indent In section 2 and section 3, we recall the definition of the Bott connection, and give affine Ricci solitons associated to the Bott connection on three-dimensional Lorentzian unimodular and non-unimodular Lie groups with the first distribution. In section 4, we give affine Ricci solitons associated to the perturbed Bott connection on three-dimensional Lorentzian Lie groups with the first distribution. In section 5, we give affine Ricci solitons associated to the Bott connection on three-dimensional Lorentzian Lie groups with the second distribution.
In section 6, we give affine Ricci solitons associated to the perturbed Bott connection on three-dimensional Lorentzian Lie groups with the second distribution. In section 7, we give affine Ricci solitons associated to the Bott connection on three-dimensional Lorentzian Lie groups with the third distribution. In section 8, we give affine Ricci solitons associated to the perturbed Bott connection on three-dimensional Lorentzian Lie groups with the third distribution.

\section{Affine Ricci Solitons associated to the Bott connection on three-dimensional Lorentzian Unimodular Lie groups with the first distribution}
\indent Firstly, throughout this paper, we shall by $\{G_i\}_{i=1,\cdot\cdot\cdot,7}$, denote the connected, simply connected three-dimensional Lie group equipped with a left-invariant Lorentzian metric $g$ and having Lie algebra $\{\mathfrak{g}_i\}_{i=1,\cdot\cdot\cdot,7}$. Let $\nabla^L$ be the Levi-Civita connection of $G_i$ and $R^L$ its curvature tensor, then \\
\begin{equation}
R^L(X,Y)Z=\nabla^L_X\nabla^L_YZ-\nabla^L_Y\nabla^L_XZ-\nabla_{[X,Y]}Z.
\end{equation}
The Ricci tensor of $(G_i,g)$ is defined by
\begin{equation}
\rho^L(X,Y)=-g(R^L(X,\widehat{e}_1)Y,\widehat{e}_1)-g(R^L(X,\widehat{e}_2)Y,\widehat{e}_2)+g(R^L(X,\widehat{e}_3)Y,\widehat{e}_3),
\end{equation}
where ${\widehat{e}_1,\widehat{e}_2,\widehat{e}_3}$ is a pseudo-orthonormal basis, with $\widehat{e}_3$ timelike.\\
\indent Nextly, we recall the definition of the Bott connection $\nabla^B$. Let $M$ be a smooth manifold, and let $TM=span\{\widehat{e}_1,\widehat{e}_2,\widehat{e}_3\}$, then took the distribution $D=span\{\widehat{e}_1,\widehat{e}_2\}$ and $D^\bot=span\{\widehat{e}_3\}$.\\
The definition of the Bott connection $\nabla^B$ is given as follows: (see \cite{G}, \cite{RK}, \cite{W})
\begin{eqnarray}
\nabla^B_XY=
       \begin{cases}
        \pi_D(\nabla^L_XY),~~~&X,Y\in\Gamma^\infty(D) \\[2pt]
       \pi_D([X,Y]),~~~&X\in\Gamma^\infty(D^\bot),Y\in\Gamma^\infty(D)\\[2pt]
       \pi_{D^\bot}([X,Y]),~~~&X\in\Gamma^\infty(D),Y\in\Gamma^\infty(D^\bot)\\[2pt]
       \pi_{D^\bot}(\nabla^L_XY),~~~&X,Y\in\Gamma^\infty(D^\bot)\\[2pt]
       \end{cases}
\end{eqnarray}
where $\pi_D$(resp. $\pi_D^\bot$) the projection on $D$ (resp. $D^\bot$).\\
We define
\begin{equation}
R^B(X,Y)Z=\nabla^B_X\nabla^B_YZ-\nabla^B_Y\nabla^B_XZ-\nabla^B_{[X,Y]}Z.
\end{equation}
The Ricci tensor of $(G_i,g)$ associated to the Bott connection $\nabla^B$ is defined by
\begin{equation}
\rho^B(X,Y)=-g(R^B(X,\widehat{e}_1)Y,\widehat{e}_1)-g(R^B(X,\widehat{e}_2)Y,\widehat{e}_2)+g(R^B(X,\widehat{e}_3)Y,\widehat{e}_3).
\end{equation}
Let
\begin{equation}
\widetilde{\rho}^B(X,Y)=\frac{\rho^B(X,Y)+\rho^B(Y,X)}{2}.
\end{equation}
We define:
\begin{equation}
(L^B_Vg)(X,Y):=g(\nabla^B_XV,Y)+g(X,\nabla^B_YV),
\end{equation}
for vector $X,Y,V$.
\begin{defn}$(G_i,g)$ is called the affine Ricci soliton associated to the connection $\nabla^B$ if it satisfies
\begin{equation}
(L^B_Vg)(X,Y)+2\widetilde{\rho}^B(X,Y)+2\mu g(X,Y)=0,
\end{equation}
where $\mu$ is a real number and $V=\mu_1\widehat{e}_1+\mu_2\widehat{e}_2+\mu_3\widehat{e}_3$ and $\mu_1,\mu_2,\mu_3$ are real numbers.
\end{defn}

\vskip 0.5 true cm
\noindent{\bf 2.1 Affine Ricci solitons of $G_1$}\\
\vskip 0.5 true cm
 By \cite{W}, we have the following Lie algebra of $G_1$ satisfies
\begin{equation}
[\widehat{e}_1, \widehat{e}_2]=\alpha\widehat{e}_1-\beta\widehat{e}_3,~~~[\widehat{e}_1, \widehat{e}_3]=-\alpha\widehat{e}_1-\beta\widehat{e}_2,~~~[\widehat{e}_2, \widehat{e}_3]=\beta\widehat{e}_1+\alpha\widehat{e}_2+\alpha\widehat{e}_3,~\alpha\neq0.
\end{equation}
where ${\widehat{e}_1,\widehat{e}_2,\widehat{e}_3}$ is a pseudo-orthonormal basis, with $\widehat{e}_3$ timelike.
\begin{lem}(\cite{G},\cite{W}) The Levi-Civita connection $\nabla^L$ of $G_1$ is given by
\begin{align}
&\nabla^L_{\widehat{e}_1}\widehat{e}_1=-\alpha\widehat{e}_2-\alpha\widehat{e}_3,~~~\nabla^L_{\widehat{e}_1}\widehat{e}_2=\alpha\widehat{e}_1-\frac{\beta}{2}\widehat{e}_3,~~~\nabla^L_{\widehat{e}_1}\widehat{e}_3=-\alpha\widehat{e}_1-\frac{\beta}{2}\widehat{e}_2,\nonumber\\
&\nabla^L_{\widehat{e}_2}\widehat{e}_1=\frac{\beta}{2}\widehat{e}_3,~~~\nabla^L_{\widehat{e}_2}\widehat{e}_2=\alpha\widehat{e}_3,~~~\nabla^L_{\widehat{e}_2}\widehat{e}_3=\frac{\beta}{2}\widehat{e}_1+\alpha\widehat{e}_2,\nonumber\\
&\nabla^L_{\widehat{e}_3}\widehat{e}_1=\frac{\beta}{2}\widehat{e}_2,~~~\nabla^L_{\widehat{e}_3}\widehat{e}_2=-\frac{\beta}{2}\widehat{e}_1-\alpha\widehat{e}_3,~~~\nabla^L_{\widehat{e}_3}\widehat{e}_3=-\alpha\widehat{e}_2.
\end{align}
\end{lem}
\begin{lem} The Bott connection $\nabla^B$ of $G_1$ is given by
\begin{align}
&\nabla^B_{\widehat{e}_1}\widehat{e}_1=-\alpha\widehat{e}_2,~~~\nabla^B_{\widehat{e}_1}\widehat{e}_2=\alpha\widehat{e}_1,~~~\nabla^B_{\widehat{e}_1}\widehat{e}_3=0,\nonumber\\
&\nabla^B_{\widehat{e}_2}\widehat{e}_1=0,~~~\nabla^L_{\widehat{e}_2}\widehat{e}_2=0,~~~\nabla^B_{\widehat{e}_2}\widehat{e}_3=\alpha\widehat{e}_3,\nonumber\\
&\nabla^B_{\widehat{e}_3}\widehat{e}_1=\alpha\widehat{e}_1+\beta\widehat{e}_2,~~~\nabla^B_{\widehat{e}_3}\widehat{e}_2=-\beta\widehat{e}_1-\alpha\widehat{e}_2,~~~\nabla^B_{\widehat{e}_3}\widehat{e}_3=0.
\end{align}
\end{lem}
\begin{lem} The curvature $R^B$ of the Bott connection $\nabla^B$ of $(G_1,g)$ is given by
\begin{align}
&R^B(\widehat{e}_1,\widehat{e}_2)\widehat{e}_1=\alpha\beta\widehat{e}_1+(\alpha^2+\beta^2)\widehat{e}_2,~~~R^B(\widehat{e}_1,\widehat{e}_2)\widehat{e}_2=-(\alpha^2+\beta^2)\widehat{e}_1-\alpha\beta\widehat{e}_2,~~~R^B(\widehat{e}_1,\widehat{e}_2)\widehat{e}_3=0,\nonumber\\
&R^B(\widehat{e}_1,\widehat{e}_3)\widehat{e}_1=-3\alpha^2\widehat{e}_2,~~~R^B(\widehat{e}_1,\widehat{e}_3)\widehat{e}_2=-\alpha^2\widehat{e}_1,~~~R^B(\widehat{e}_1,\widehat{e}_3)\widehat{e}_3=\alpha\beta\widehat{e}_3,\nonumber\\
&R^B(\widehat{e}_2,\widehat{e}_3)\widehat{e}_1=-\alpha^2\widehat{e}_1,~~~R^B(\widehat{e}_2,\widehat{e}_3)\widehat{e}_2=\alpha^2\widehat{e}_2,~~~R^B(\widehat{e}_2,\widehat{e}_3)\widehat{e}_3=-\alpha^2\widehat{e}_3.
\end{align}
\end{lem}
By (2.5), we have
\begin{align}
&\rho^B(\widehat{e}_1,\widehat{e}_1)=-(\alpha^2+\beta^2),~~~\rho^B(\widehat{e}_1,\widehat{e}_2)=\alpha\beta,~~~\rho^B(\widehat{e}_1,\widehat{e}_3)=-\alpha\beta,\nonumber\\
&\rho^B(\widehat{e}_2,\widehat{e}_1)=\alpha\beta,~~~\rho^B(\widehat{e}_2,\widehat{e}_2)=-(\alpha^2+\beta^2),~~~\rho^B(\widehat{e}_2,\widehat{e}_3)=\alpha^2,\nonumber\\
&\rho^B(\widehat{e}_3,\widehat{e}_1)=\rho^B(\widehat{e}_3,\widehat{e}_2)=\rho^B(\widehat{e}_3,\widehat{e}_3)=0.
\end{align}
Then,
\begin{align}
&\widetilde{\rho}^B(\widehat{e}_1,\widehat{e}_1)=-(\alpha^2+\beta^2),~~~\widetilde{\rho}^B(\widehat{e}_1,\widehat{e}_2)=\alpha\beta,~~~\widetilde{\rho}^B(\widehat{e}_1,\widehat{e}_3)=-\frac{\alpha\beta}{2},\nonumber\\
&\widetilde{\rho}^B(\widehat{e}_2,\widehat{e}_2)=-(\alpha^2+\beta^2),~~~\widetilde{\rho}^B(\widehat{e}_2,\widehat{e}_3)=\frac{\alpha^2}{2},~~~\widetilde{\rho}^B(\widehat{e}_3,\widehat{e}_3)=0.
\end{align}
By (2.7), we have
\begin{align}
&(L^B_Vg)(\widehat{e}_1,\widehat{e}_1)=2\mu_2\alpha,~~~(L^B_Vg)(\widehat{e}_1,\widehat{e}_2)=-\mu_1\alpha,~~~(L^B_Vg)(\widehat{e}_1,\widehat{e}_3)=\mu_1\alpha-\mu_2\beta,\nonumber\\
&(L^B_Vg)(\widehat{e}_2,\widehat{e}_2)=0,~~~(L^B_Vg)(\widehat{e}_2,\widehat{e}_3)=\mu_1\beta-(\mu_2+\mu_3)\alpha,~~~(L^B_Vg)(\widehat{e}_3,\widehat{e}_3)=0.
\end{align}
Then, if $(G_1,g,V)$  is an affine Ricci soliton associated to the Bott connection $\nabla^B$, by (2.8), we have the following six equations:
\begin{eqnarray}
       \begin{cases}
        \mu_2\alpha-\alpha^2-\beta^2+\mu=0 \\[2pt]
       2\alpha\beta-\mu_1\alpha=0\\[2pt]
       \mu_1\alpha-\mu_2\beta-\alpha\beta=0\\[2pt]
       \alpha^2+\beta^2-\mu=0\\[2pt]
       (\mu_2+\mu_3)\alpha-\mu_1\beta-\alpha^2=0\\[2pt]
      \mu=0\\[2pt]
       \end{cases}
\end{eqnarray}
By solving (2.16) , we get $\alpha=0$, there is a contradiction. So\\
\begin{thm}
$(G_1, g, V)$ is not an affine Ricci soliton associated to the Bott connection $\nabla^B$.
\end{thm}

\vskip 0.5 true cm
\noindent{\bf 2.2 Affine Ricci solitons of $G_2$}\\
\vskip 0.5 true cm
\indent By \cite{W}, we have the following Lie algebra of $G_2$ satisfies
\begin{equation}
[\widehat{e}_1, \widehat{e}_2]=\gamma\widehat{e}_2-\beta\widehat{e}_3,~~~[\widehat{e}_1, \widehat{e}_3]=-\beta\widehat{e}_2-\gamma\widehat{e}_3,~~~[\widehat{e}_2, \widehat{e}_3]=\alpha\widehat{e}_1,~\gamma\neq0.
\end{equation}
where ${\widehat{e}_1,\widehat{e}_2,\widehat{e}_3}$ is a pseudo-orthonormal basis, with $\widehat{e}_3$ timelike.
\begin{lem}(\cite{G},\cite{W}) The Levi-Civita connection $\nabla^L$ of $G_2$ is given by
\begin{align}
&\nabla^L_{\widehat{e}_1}\widehat{e}_1=0,~~~\nabla^L_{\widehat{e}_1}\widehat{e}_2=(\frac{\alpha}{2}-\beta)\widehat{e}_3,~~~\nabla^L_{\widehat{e}_1}\widehat{e}_3=(\frac{\alpha}{2}-\beta)\widehat{e}_2,\nonumber\\
&\nabla^L_{\widehat{e}_2}\widehat{e}_1=-\gamma\widehat{e}_2+\frac{\alpha}{2}\widehat{e}_3,~~~\nabla^L_{\widehat{e}_2}\widehat{e}_2=\gamma\widehat{e}_1,~~~\nabla^L_{\widehat{e}_2}\widehat{e}_3=\frac{\alpha}{2}\widehat{e}_1,\nonumber\\
&\nabla^L_{\widehat{e}_3}\widehat{e}_1=\frac{\alpha}{2}\widehat{e}_2+\gamma\widehat{e}_3,~~~\nabla^L_{\widehat{e}_3}\widehat{e}_2=-\frac{\alpha}{2}\widehat{e}_1,~~~\nabla^L_{\widehat{e}_3}\widehat{e}_3=\gamma\widehat{e}_1.
\end{align}
\end{lem}
\begin{lem} The Bott connection $\nabla^B$ of $G_2$ is given by
\begin{align}
&\nabla^B_{\widehat{e}_1}\widehat{e}_1=0,~~~\nabla^B_{\widehat{e}_1}\widehat{e}_2=0,~~~\nabla^B_{\widehat{e}_1}\widehat{e}_3=-\gamma\widehat{e}_3,\nonumber\\
&\nabla^B_{\widehat{e}_2}\widehat{e}_1=-\gamma\widehat{e}_2,~~~\nabla^B_{\widehat{e}_2}\widehat{e}_2=\gamma\widehat{e}_1,~~~\nabla^B_{\widehat{e}_2}\widehat{e}_3=0,\nonumber\\
&\nabla^B_{\widehat{e}_3}\widehat{e}_1=\beta\widehat{e}_2,~~~\nabla^B_{\widehat{e}_3}\widehat{e}_2=-\alpha\widehat{e}_1,~~~\nabla^B_{\widehat{e}_3}\widehat{e}_3=0.
\end{align}
\end{lem}
\begin{lem} The curvature $R^B$ of the Bott connection $\nabla^B$ of $(G_2,g)$ is given by
\begin{align}
&R^B(\widehat{e}_1,\widehat{e}_2)\widehat{e}_1=(\beta^2+\gamma^2)\widehat{e}_2,~~~R^B(\widehat{e}_1,\widehat{e}_2)\widehat{e}_2=-(\gamma^2+\alpha\beta)\widehat{e}_1,~~~R^B(\widehat{e}_1,\widehat{e}_2)\widehat{e}_3=0,\nonumber\\
&R^B(\widehat{e}_1,\widehat{e}_3)\widehat{e}_1=0,~~~R^B(\widehat{e}_1,\widehat{e}_3)\widehat{e}_2=\gamma(\alpha-\beta)\widehat{e}_1,~~~R^B(\widehat{e}_1,\widehat{e}_3)\widehat{e}_3=0,\nonumber\\
&R^B(\widehat{e}_2,\widehat{e}_3)\widehat{e}_1=\gamma(\beta-\alpha)\widehat{e}_1,~~~R^B(\widehat{e}_2,\widehat{e}_3)\widehat{e}_2=\gamma(\alpha-\beta)\widehat{e}_2,~~~R^B(\widehat{e}_2,\widehat{e}_3)\widehat{e}_3=\alpha\gamma\widehat{e}_3.
\end{align}
\end{lem}
By (2.5), we have
\begin{align}
&\rho^B(\widehat{e}_1,\widehat{e}_1)=-(\beta^2+\gamma^2),~~~\rho^B(\widehat{e}_1,\widehat{e}_2)=0,~~~\rho^B(\widehat{e}_1,\widehat{e}_3)=0,\nonumber\\
&\rho^B(\widehat{e}_2,\widehat{e}_1)=0,~~~\rho^B(\widehat{e}_2,\widehat{e}_2)=-(\gamma^2+\alpha\beta),~~~\rho^B(\widehat{e}_2,\widehat{e}_3)=-\alpha\gamma,\nonumber\\
&\rho^B(\widehat{e}_3,\widehat{e}_1)=\rho^B(\widehat{e}_3,\widehat{e}_2)=\rho^B(\widehat{e}_3,\widehat{e}_3)=0.
\end{align}
Then,
\begin{align}
&\widetilde{\rho}^B(\widehat{e}_1,\widehat{e}_1)=-(\beta^2+\gamma^2),~~~\widetilde{\rho}^B(\widehat{e}_1,\widehat{e}_2)=0,~~~\widetilde{\rho}^B(\widehat{e}_1,\widehat{e}_3)=0,\nonumber\\
&\widetilde{\rho}^B(\widehat{e}_2,\widehat{e}_2)=-(\gamma^2+\alpha\beta),~~~\widetilde{\rho}^B(\widehat{e}_2,\widehat{e}_3)=-\frac{\alpha\gamma}{2},~~~\widetilde{\rho}^B(\widehat{e}_3,\widehat{e}_3)=0.
\end{align}
By (2.7), we have
\begin{align}
&(L^B_Vg)(\widehat{e}_1,\widehat{e}_1)=0,~~~(L^B_Vg)(\widehat{e}_1,\widehat{e}_2)=\mu_2\gamma,~~~(L^B_Vg)(\widehat{e}_1,\widehat{e}_3)=\mu_3\gamma-\mu_2\alpha,\nonumber\\
&(L^B_Vg)(\widehat{e}_2,\widehat{e}_2)=-2\mu_1\gamma,~~~(L^B_Vg)(\widehat{e}_2,\widehat{e}_3)=\mu_1\beta,~~~(L^B_Vg)(\widehat{e}_3,\widehat{e}_3)=0.
\end{align}
Then, if $(G_2,g,V)$  is an affine Ricci soliton associated to the Bott connection $\nabla^B$, by (2.8), we have the following six equations:
\begin{eqnarray}
       \begin{cases}
        \beta^2+\gamma^2-\mu=0 \\[2pt]
       \mu_2\gamma=0\\[2pt]
       \mu_3\gamma-\mu_2\alpha=0\\[2pt]
       \mu_1\gamma+\gamma^2+\alpha\beta-\mu=0\\[2pt]
       \mu_1\beta-\alpha\gamma=0\\[2pt]
       \mu=0\\[2pt]
       \end{cases}
\end{eqnarray}
By solving (2.24), we get $\beta=\gamma=0$, there is a contradiction. So\\
\begin{thm}
$(G_2, g, V)$ is not an affine Ricci soliton associated to the Bott connection $\nabla^B$.
\end{thm}

\vskip 0.5 true cm
\noindent{\bf 2.3 Affine Ricci solitons of $G_3$}\\
\vskip 0.5 true cm
\indent By \cite{W}, we have the following Lie algebra of $G_3$ satisfies
\begin{equation}
[\widehat{e}_1, \widehat{e}_2]=-\gamma\widehat{e}_3,~~~[\widehat{e}_1, \widehat{e}_3]=-\beta\widehat{e}_2,~~~[\widehat{e}_2, \widehat{e}_3]=\alpha\widehat{e}_1.
\end{equation}
where ${\widehat{e}_1,\widehat{e}_2,\widehat{e}_3}$ is a pseudo-orthonormal basis, with $\widehat{e}_3$ timelike.
\begin{lem}(\cite{G},\cite{W}) The Levi-Civita connection $\nabla^L$ of $G_3$ is given by
\begin{align}
&\nabla^L_{\widehat{e}_1}\widehat{e}_1=0,~~~\nabla^L_{\widehat{e}_1}\widehat{e}_2=m_1\widehat{e}_3,~~~\nabla^L_{\widehat{e}_1}\widehat{e}_3=m_1\widehat{e}_2,\nonumber\\
&\nabla^L_{\widehat{e}_2}\widehat{e}_1=m_2\widehat{e}_3,~~~\nabla^L_{\widehat{e}_2}\widehat{e}_2=0,~~~\nabla^L_{\widehat{e}_2}\widehat{e}_3=m_2\widehat{e}_1,\nonumber\\
&\nabla^L_{\widehat{e}_3}\widehat{e}_1=m_3\widehat{e}_2,~~~\nabla^L_{\widehat{e}_3}\widehat{e}_2=-m_3\widehat{e}_1,~~~\nabla^L_{\widehat{e}_3}\widehat{e}_3=0.
\end{align}
where\\
 \begin{align}
 m_1=\frac{1}{2}(\alpha-\beta-\gamma),~~~m_2=\frac{1}{2}(\alpha-\beta+\gamma),~~~m_3=\frac{1}{2}(\alpha+\beta-\gamma).
 \end{align}
\end{lem}
\begin{lem} The Bott connection $\nabla^B$ of $G_3$ is given by
\begin{align}
&\nabla^B_{\widehat{e}_1}\widehat{e}_1=0,~~~\nabla^B_{\widehat{e}_1}\widehat{e}_2=0,~~~\nabla^B_{\widehat{e}_1}\widehat{e}_3=-\gamma\widehat{e}_3,\nonumber\\
&\nabla^B_{\widehat{e}_2}\widehat{e}_1=0,~~~\nabla^B_{\widehat{e}_2}\widehat{e}_2=0,~~~\nabla^B_{\widehat{e}_2}\widehat{e}_3=0,\nonumber\\
&\nabla^B_{\widehat{e}_3}\widehat{e}_1=\beta\widehat{e}_2,~~~\nabla^B_{\widehat{e}_3}\widehat{e}_2=-\alpha\widehat{e}_1,~~~\nabla^B_{\widehat{e}_3}\widehat{e}_3=0.
\end{align}
\end{lem}
\begin{lem} The curvature $R^B$ of the Bott connection $\nabla^B$ of $(G_3,g)$ is given by
\begin{align}
&R^B(\widehat{e}_1,\widehat{e}_2)\widehat{e}_1=\beta\gamma\widehat{e}_2,~~~R^B(\widehat{e}_1,\widehat{e}_2)\widehat{e}_2=-\alpha\gamma\widehat{e}_1,~~~R^B(\widehat{e}_1,\widehat{e}_2)\widehat{e}_3=0,\nonumber\\
&R^B(\widehat{e}_1,\widehat{e}_3)\widehat{e}_1=0,~~~R^B(\widehat{e}_1,\widehat{e}_3)\widehat{e}_2=0,~~~R^B(\widehat{e}_1,\widehat{e}_3)\widehat{e}_3=0,\nonumber\\
&R^B(\widehat{e}_2,\widehat{e}_3)\widehat{e}_1=0,~~~R^B(\widehat{e}_2,\widehat{e}_3)\widehat{e}_2=0,~~~R^B(\widehat{e}_2,\widehat{e}_3)\widehat{e}_3=0.
\end{align}
\end{lem}
By (2.5), we have
\begin{align}
&\rho^B(\widehat{e}_1,\widehat{e}_1)=-\beta\gamma,~~~\rho^B(\widehat{e}_1,\widehat{e}_2)=0,~~~\rho^B(\widehat{e}_1,\widehat{e}_3)=0,\nonumber\\
&\rho^B(\widehat{e}_2,\widehat{e}_1)=0,~~~\rho^B(\widehat{e}_2,\widehat{e}_2)=-\alpha\gamma,~~~\rho^B(\widehat{e}_2,\widehat{e}_3)=0,\nonumber\\
&\rho^B(\widehat{e}_3,\widehat{e}_1)=\rho^B(\widehat{e}_3,\widehat{e}_2)=\rho^B(\widehat{e}_3,\widehat{e}_3)=0.
\end{align}
Then,
\begin{align}
&\widetilde{\rho}^B(\widehat{e}_1,\widehat{e}_1)=-\beta\gamma,~~~\widetilde{\rho}^B(\widehat{e}_1,\widehat{e}_2)=0,~~~\widetilde{\rho}^B(\widehat{e}_1,\widehat{e}_3)=0,\nonumber\\
&\widetilde{\rho}^B(\widehat{e}_2,\widehat{e}_2)=-\alpha\gamma,~~~\widetilde{\rho}^B(\widehat{e}_2,\widehat{e}_3)=0,~~~\widetilde{\rho}^B(\widehat{e}_3,\widehat{e}_3)=0.
\end{align}
By (2.7), we have
\begin{align}
&(L^B_Vg)(\widehat{e}_1,\widehat{e}_1)=0,~~~(L^B_Vg)(\widehat{e}_1,\widehat{e}_2)=0,~~~(L^B_Vg)(\widehat{e}_1,\widehat{e}_3)=-\mu_2\alpha,\nonumber\\
&(L^B_Vg)(\widehat{e}_2,\widehat{e}_2)=0,~~~(L^B_Vg)(\widehat{e}_2,\widehat{e}_3)=\mu_1\beta,~~~(L^B_Vg)(\widehat{e}_3,\widehat{e}_3)=0.
\end{align}
Then, if $(G_3,g,V)$  is an affine Ricci soliton associated to the Bott connection $\nabla^B$, by (2.8), we have the following five equations:
\begin{eqnarray}
       \begin{cases}
        \mu-\beta\gamma=0 \\[2pt]
       \mu_2\alpha=0\\[2pt]
       \mu-\alpha\gamma=0\\[2pt]
       \mu_1\beta=0\\[2pt]
       \mu=0\\[2pt]
       \end{cases}
\end{eqnarray}
By solving (2.33), we get \\
\begin{thm}
$(G_3, g, V)$ is an affine Ricci soliton associated to the Bott connection $\nabla^B$ if and only if\\
\begin{eqnarray}
&&(1)\mu=\gamma=\alpha\mu_2=\mu_1\beta=0;\nonumber\\
&&(2)\mu=\alpha=\beta=0,~~~\gamma\neq0\nonumber.
\end{eqnarray}
\end{thm}
\vskip 0.5 true cm
\noindent{\bf 2.4 Affine Ricci solitons of $G_4$}\\
\vskip 0.5 true cm
\indent By \cite{W}, we have the following Lie algebra of $G_4$ satisfies
\begin{equation}
[\widehat{e}_1, \widehat{e}_2]=-\widehat{e}_2+(2\eta-\beta)\widehat{e}_3,~~\eta=\pm 1,~~~[\widehat{e}_1, \widehat{e}_3]=-\beta\widehat{e}_2+\widehat{e}_3,~~~[\widehat{e}_2, \widehat{e}_3]=\alpha\widehat{e}_1.
\end{equation}
where ${\widehat{e}_1,\widehat{e}_2,\widehat{e}_3}$ is a pseudo-orthonormal basis, with $\widehat{e}_3$ timelike.
\begin{lem}(\cite{G},\cite{W}) The Levi-Civita connection $\nabla^L$ of $G_4$ is given by
\begin{align}
&\nabla^L_{\widehat{e}_1}\widehat{e}_1=0,~~~\nabla^L_{\widehat{e}_1}\widehat{e}_2=n_1\widehat{e}_3,~~~\nabla^L_{\widehat{e}_1}\widehat{e}_3=n_1\widehat{e}_2,\nonumber\\
&\nabla^L_{\widehat{e}_2}\widehat{e}_1=\widehat{e}_2+n_2\widehat{e}_3,~~~\nabla^L_{\widehat{e}_2}\widehat{e}_2=-\widehat{e}_1,~~~\nabla^L_{\widehat{e}_2}\widehat{e}_3=n_2\widehat{e}_1,\nonumber\\
&\nabla^L_{\widehat{e}_3}\widehat{e}_1=n_3\widehat{e}_2-\widehat{e}_3,~~~\nabla^L_{\widehat{e}_3}\widehat{e}_2=-n_3\widehat{e}_1,~~~\nabla^L_{\widehat{e}_3}\widehat{e}_3=-\widehat{e}_1.
\end{align}
where\\
 \begin{align}
 n_1=\frac{\alpha}{2}+\eta-\beta,~~~n_2=\frac{\alpha}{2}-\gamma,~~~n_3=\frac{\alpha}{2}+\gamma.
 \end{align}
\end{lem}
\begin{lem} The Bott connection $\nabla^B$ of $G_4$ is given by
\begin{align}
&\nabla^B_{\widehat{e}_1}\widehat{e}_1=0,~~~\nabla^B_{\widehat{e}_1}\widehat{e}_2=0,~~~\nabla^B_{\widehat{e}_1}\widehat{e}_3=\widehat{e}_3,\nonumber\\
&\nabla^B_{\widehat{e}_2}\widehat{e}_1=\widehat{e}_2,~~~\nabla^B_{\widehat{e}_2}\widehat{e}_2=-\widehat{e}_1,~~~\nabla^B_{\widehat{e}_2}\widehat{e}_3=0,\nonumber\\
&\nabla^B_{\widehat{e}_3}\widehat{e}_1=\beta\widehat{e}_2,~~~\nabla^B_{\widehat{e}_3}\widehat{e}_2=-\alpha\widehat{e}_1,~~~\nabla^B_{\widehat{e}_3}\widehat{e}_3=0.
\end{align}
\end{lem}
\begin{lem} The curvature $R^B$ of the Bott connection $\nabla^B$ of $(G_4,g)$ is given by
\begin{align}
&R^B(\widehat{e}_1,\widehat{e}_2)\widehat{e}_1=(\beta-\eta)^2\widehat{e}_2,~~~R^B(\widehat{e}_1,\widehat{e}_2)\widehat{e}_2=(2\alpha\eta-\alpha\beta-1)\widehat{e}_1,~~~R^B(\widehat{e}_1,\widehat{e}_2)\widehat{e}_3=0,\nonumber\\
&R^B(\widehat{e}_1,\widehat{e}_3)\widehat{e}_1=0,~~~R^B(\widehat{e}_1,\widehat{e}_3)\widehat{e}_2=(\alpha-\beta)\widehat{e}_1,~~~R^B(\widehat{e}_1,\widehat{e}_3)\widehat{e}_3=0,\nonumber\\
&R^B(\widehat{e}_2,\widehat{e}_3)\widehat{e}_1=(\alpha-\beta)\widehat{e}_1,~~~R^B(\widehat{e}_2,\widehat{e}_3)\widehat{e}_2=(\beta-\alpha)\widehat{e}_2,~~~R^B(\widehat{e}_2,\widehat{e}_3)\widehat{e}_3=-\alpha\widehat{e}_3.
\end{align}
\end{lem}
By (2.5), we have
\begin{align}
&\rho^B(\widehat{e}_1,\widehat{e}_1)=-(\beta-\eta)^2,~~~\rho^B(\widehat{e}_1,\widehat{e}_2)=0,~~~\rho^B(\widehat{e}_1,\widehat{e}_3)=0,\nonumber\\
&\rho^B(\widehat{e}_2,\widehat{e}_1)=(2\alpha\eta-\alpha\beta-1),~~~\rho^B(\widehat{e}_2,\widehat{e}_2)=\alpha,~~~\rho^B(\widehat{e}_2,\widehat{e}_3)=0,\nonumber\\
&\rho^B(\widehat{e}_3,\widehat{e}_1)=\rho^B(\widehat{e}_3,\widehat{e}_2)=\rho^B(\widehat{e}_3,\widehat{e}_3)=0.
\end{align}
Then,
\begin{align}
&\widetilde{\rho}^B(\widehat{e}_1,\widehat{e}_1)=-(\beta-\eta)^2,~~~\widetilde{\rho}^B(\widehat{e}_1,\widehat{e}_2)=0,~~~\widetilde{\rho}^B(\widehat{e}_1,\widehat{e}_3)=0,\nonumber\\
&\widetilde{\rho}^B(\widehat{e}_2,\widehat{e}_2)=(2\alpha\eta-\alpha\beta-1),~~~\widetilde{\rho}^B(\widehat{e}_2,\widehat{e}_3)=\frac{\alpha}{2},~~~\widetilde{\rho}^B(\widehat{e}_3,\widehat{e}_3)=0.
\end{align}
By (2.7), we have
\begin{align}
&(L^B_Vg)(\widehat{e}_1,\widehat{e}_1)=0,~~~(L^B_Vg)(\widehat{e}_1,\widehat{e}_2)=-\mu_2,~~~(L^B_Vg)(\widehat{e}_1,\widehat{e}_3)=-\mu_3-\mu_2\alpha,\nonumber\\
&(L^B_Vg)(\widehat{e}_2,\widehat{e}_2)=2\mu_1,~~~(L^B_Vg)(\widehat{e}_2,\widehat{e}_3)=\mu_1\beta,~~~(L^B_Vg)(\widehat{e}_3,\widehat{e}_3)=0.
\end{align}
Then, if $(G_4,g,V)$  is an affine Ricci soliton associated to the Bott connection $\nabla^B$, by (2.8), we have the following six equations:
\begin{eqnarray}
       \begin{cases}
        (\beta-\eta)^2-\mu=0 \\[2pt]
       \mu_2=0\\[2pt]
       \mu_2\alpha+\mu_3=0\\[2pt]
       \mu_1+2\alpha\eta-\alpha\beta-1+\mu=0\\[2pt]
              \mu_1\beta+\alpha=0\\[2pt]
       \mu=0\\[2pt]
       \end{cases}
\end{eqnarray}
By solving (2.42), we get
\begin{thm}
$(G_4, g, V)$ is not an affine Ricci soliton associated to the Bott connection $\nabla^B$.
\end{thm}

\section{Affine Ricci Solitons associated to the Bott connection on three-dimensional Lorentzian Non-unimodular Lie groups with the first distribution}

\vskip 0.5 true cm
\noindent{\bf 3.1 Affine Ricci solitons of $G_5$}\\
\vskip 0.5 true cm
\indent By \cite{W}, we have the following Lie algebra of $G_5$ satisfies
\begin{equation}
[\widehat{e}_1, \widehat{e}_2]=0,~~~[\widehat{e}_1, \widehat{e}_3]=\alpha\widehat{e}_1+\beta\widehat{e}_2,~~~[\widehat{e}_2, \widehat{e}_3]=\gamma\widehat{e}_1+\delta\widehat{e}_2,~~~\alpha+\delta\neq 0,~~~\alpha\gamma+\beta\delta=0.
\end{equation}
where ${\widehat{e}_1,\widehat{e}_2,\widehat{e}_3}$ is a pseudo-orthonormal basis, with $\widehat{e}_3$ timelike.
\begin{lem}(\cite{G},\cite{W}) The Levi-Civita connection $\nabla^L$ of $G_5$ is given by
\begin{align}
&\nabla^L_{\widehat{e}_1}\widehat{e}_1=\alpha\widehat{e}_3,~~~\nabla^L_{\widehat{e}_1}\widehat{e}_2=\frac{\beta+\gamma}{2}\widehat{e}_3,~~~\nabla^L_{\widehat{e}_1}\widehat{e}_3=\alpha\widehat{e}_1+\frac{\beta+\gamma}{2}\widehat{e}_2,\nonumber\\
&\nabla^L_{\widehat{e}_2}\widehat{e}_1=\frac{\beta+\gamma}{2}\widehat{e}_3,~~~\nabla^L_{\widehat{e}_2}\widehat{e}_2=\delta\widehat{e}_3,~~~\nabla^L_{\widehat{e}_2}\widehat{e}_3=\frac{\beta+\gamma}{2}\widehat{e}_1+\delta\widehat{e}_2,\nonumber\\
&\nabla^L_{\widehat{e}_3}\widehat{e}_1=\frac{\gamma-\beta}{2}\widehat{e}_2,~~~\nabla^L_{\widehat{e}_3}\widehat{e}_2=\frac{\beta-\gamma}{2}\widehat{e}_1,~~~\nabla^L_{\widehat{e}_3}\widehat{e}_3=0.
\end{align}
\end{lem}
\begin{lem} The Bott connection $\nabla^B$ of $G_5$ is given by
\begin{align}
&\nabla^B_{\widehat{e}_1}\widehat{e}_1=0,~~~\nabla^B_{\widehat{e}_1}\widehat{e}_2=0,~~~\nabla^B_{\widehat{e}_1}\widehat{e}_3=0,\nonumber\\
&\nabla^B_{\widehat{e}_2}\widehat{e}_1=0,~~~\nabla^B_{\widehat{e}_2}\widehat{e}_2=0,~~~\nabla^B_{\widehat{e}_2}\widehat{e}_3=0,\nonumber\\
&\nabla^B_{\widehat{e}_3}\widehat{e}_1=-\alpha\widehat{e}_1-\beta\widehat{e}_2,~~~\nabla^B_{\widehat{e}_3}\widehat{e}_2=-\gamma\widehat{e}_1-\delta\widehat{e}_2,~~~\nabla^B_{\widehat{e}_3}\widehat{e}_3=0.
\end{align}
\end{lem}
\begin{lem} The curvature $R^B$ of the Bott connection $\nabla^B$ of $(G_5,g)$ is given by
\begin{align}
R^B(\widehat{e}_s,\widehat{e}_t)\widehat{e}_p=0,
\end{align}
for any $(s,t,p)$.
\end{lem}
By (2.5), we have
\begin{align}
\rho^B(\widehat{e}_s,\widehat{e}_t)=0,
\end{align}
Then,
\begin{align}
&\widetilde{\rho}^B(\widehat{e}_s,\widehat{e}_t)=0,
\end{align}
for any pairs $(s,t)$.\\
By (2.7), we have
\begin{align}
&(L^B_Vg)(\widehat{e}_1,\widehat{e}_1)=0,~~~(L^B_Vg)(\widehat{e}_1,\widehat{e}_2)=-\mu_2,~~~(L^B_Vg)(\widehat{e}_1,\widehat{e}_3)=-\mu_1\alpha-\mu_2\gamma,\nonumber\\
&(L^B_Vg)(\widehat{e}_2,\widehat{e}_2)=2\mu_1,~~~(L^B_Vg)(\widehat{e}_2,\widehat{e}_3)=\mu_1\beta-\mu_2\delta,~~~(L^B_Vg)(\widehat{e}_3,\widehat{e}_3)=0.
\end{align}
Then, if $(G_5,g,V)$  is an affine Ricci soliton associated to the Bott connection $\nabla^B$, by (2.8), we have the following three equations:
\begin{eqnarray}
       \begin{cases}
       \mu=0\\[2pt]
       \mu_1\alpha+\mu_2\gamma=0\\[2pt]
       \mu_1\beta+\mu_2\delta=0\\[2pt]
       \end{cases}
\end{eqnarray}
By solving (3.8), we get \\
\begin{thm}
$(G_5, g, V)$ is an affine Ricci soliton associated to the Bott connection $\nabla^B$ if and only if\\
\begin{eqnarray}
&&(1)\mu=\mu_1=\mu_2=0,~~~ \alpha+\delta\neq 0, ~~~\alpha\gamma+\beta\delta=0;\nonumber\\
&&(2)\mu=\mu_2=\alpha=\beta=0,~~~\mu_1\neq0,~~~\delta\neq0;\nonumber\\
&&(3)\mu=\mu_1=\gamma=\delta=0,~~~\mu_2\neq0,~~~\alpha\neq0\nonumber.
\end{eqnarray}
\end{thm}
\vskip 0.5 true cm
\noindent{\bf 3.2 Affine Ricci solitons of $G_6$}\\
\vskip 0.5 true cm
\indent By \cite{W}, we have the following Lie algebra of $G_6$ satisfies
\begin{equation}
[\widehat{e}_1, \widehat{e}_2]=\alpha\widehat{e}_2+\beta\widehat{e}_3,~~~[\widehat{e}_1, \widehat{e}_3]=\gamma\widehat{e}_2+\delta\widehat{e}_3,~~~[\widehat{e}_2, \widehat{e}_3]=0,~~~\alpha+\delta\neq 0,~~~\alpha\gamma-\beta\delta=0.
\end{equation}
where ${\widehat{e}_1,\widehat{e}_2,\widehat{e}_3}$ is a pseudo-orthonormal basis, with $\widehat{e}_3$ timelike.
\begin{lem}(\cite{G},\cite{W}) The Levi-Civita connection $\nabla^L$ of $G_6$ is given by
\begin{align}
&\nabla^L_{\widehat{e}_1}\widehat{e}_1=0,~~~\nabla^L_{\widehat{e}_1}\widehat{e}_2=\frac{\beta+\gamma}{2}\widehat{e}_3,~~~\nabla^L_{\widehat{e}_1}\widehat{e}_3=\frac{\beta+\gamma}{2}\widehat{e}_2,\nonumber\\
&\nabla^L_{\widehat{e}_2}\widehat{e}_1=-\alpha\widehat{e}_2-\frac{\beta-\gamma}{2}\widehat{e}_3,~~~\nabla^L_{\widehat{e}_2}\widehat{e}_2=\alpha\widehat{e}_1,~~~\nabla^L_{\widehat{e}_2}\widehat{e}_3=\frac{\gamma-\beta}{2}\widehat{e}_1,\nonumber\\
&\nabla^L_{\widehat{e}_3}\widehat{e}_1=\frac{\beta-\gamma}{2}\widehat{e}_2-\delta\widehat{e}_3,~~~\nabla^L_{\widehat{e}_3}\widehat{e}_2=-\frac{\gamma-\beta}{2}\widehat{e}_1,~~~\nabla^L_{\widehat{e}_3}\widehat{e}_3=-\delta\widehat{e}_1.
\end{align}
\end{lem}
\begin{lem} The Bott connection $\nabla^B$ of $G_6$ is given by
\begin{align}
&\nabla^B_{\widehat{e}_1}\widehat{e}_1=0,~~~\nabla^B_{\widehat{e}_1}\widehat{e}_2=0,~~~\nabla^B_{\widehat{e}_1}\widehat{e}_3=\delta\widehat{e}_3,\nonumber\\
&\nabla^B_{\widehat{e}_2}\widehat{e}_1=-\alpha\widehat{e}_2,~~~\nabla^B_{\widehat{e}_2}\widehat{e}_2=\alpha\widehat{e}_1,~~~\nabla^B_{\widehat{e}_2}\widehat{e}_3=0,\nonumber\\
&\nabla^B_{\widehat{e}_3}\widehat{e}_1=-\gamma\widehat{e}_2,~~~\nabla^B_{\widehat{e}_3}\widehat{e}_2=0,~~~\nabla^B_{\widehat{e}_3}\widehat{e}_3=0.
\end{align}
\end{lem}
\begin{lem} The curvature $R^B$ of the Bott connection $\nabla^B$ of $(G_6,g)$ is given by
\begin{align}
&R^B(\widehat{e}_1,\widehat{e}_2)\widehat{e}_1=(\alpha^2+\beta\gamma)\widehat{e}_2,~~~R^B(\widehat{e}_1,\widehat{e}_2)\widehat{e}_2=-\alpha^2\widehat{e}_1,~~~R^B(\widehat{e}_1,\widehat{e}_2)\widehat{e}_3=0,\nonumber\\
&R^B(\widehat{e}_1,\widehat{e}_3)\widehat{e}_1=\gamma(\alpha+\delta)\widehat{e}_2,~~~R^B(\widehat{e}_1,\widehat{e}_3)\widehat{e}_2=-\alpha\gamma\widehat{e}_1,~~~R^B(\widehat{e}_1,\widehat{e}_3)\widehat{e}_3=0,\nonumber\\
&R^B(\widehat{e}_2,\widehat{e}_3)\widehat{e}_1=-\alpha\gamma\widehat{e}_1,~~~R^B(\widehat{e}_2,\widehat{e}_3)\widehat{e}_2=\alpha\gamma\widehat{e}_2,~~~R^B(\widehat{e}_2,\widehat{e}_3)\widehat{e}_3=0.
\end{align}
\end{lem}
By (2.5), we have
\begin{align}
&\rho^B(\widehat{e}_1,\widehat{e}_1)=-(\alpha^2+\beta\gamma),~~~\rho^B(\widehat{e}_1,\widehat{e}_2)=\rho^B(\widehat{e}_1,\widehat{e}_3)=0,\nonumber\\
&\rho^B(\widehat{e}_2,\widehat{e}_1)=0,~~~\rho^B(\widehat{e}_2,\widehat{e}_2)=-\alpha^2,~~~\rho^B(\widehat{e}_2,\widehat{e}_3)=0,\nonumber\\
&\rho^B(\widehat{e}_3,\widehat{e}_1)=\rho^B(\widehat{e}_3,\widehat{e}_2)=\rho^B(\widehat{e}_3,\widehat{e}_3)=0.
\end{align}
Then,
\begin{align}
&\widetilde{\rho}^B(\widehat{e}_1,\widehat{e}_1)=-(\alpha^2+\beta\gamma),~~~\widetilde{\rho}^B(\widehat{e}_1,\widehat{e}_2)=\widetilde{\rho}^B(\widehat{e}_1,\widehat{e}_3)=0,\nonumber\\
&\widetilde{\rho}^B(\widehat{e}_2,\widehat{e}_2)=-\alpha^2,~~~\widetilde{\rho}^B(\widehat{e}_2,\widehat{e}_3)=0,~~~\widetilde{\rho}^B(\widehat{e}_3,\widehat{e}_3)=0.
\end{align}
By (2.7), we have
\begin{align}
&(L^B_Vg)(\widehat{e}_1,\widehat{e}_1)=0,~~~(L^B_Vg)(\widehat{e}_1,\widehat{e}_2)=\mu_2\alpha,~~~(L^B_Vg)(\widehat{e}_1,\widehat{e}_3)=-\mu_3\delta,\nonumber\\
&(L^B_Vg)(\widehat{e}_2,\widehat{e}_2)=-2\mu_1\alpha,~~~(L^B_Vg)(\widehat{e}_2,\widehat{e}_3)=-\mu_1\gamma,~~~(L^B_Vg)(\widehat{e}_3,\widehat{e}_3)=0.
\end{align}
Then, if $(G_6,g,V)$  is an affine Ricci soliton associated to the Bott connection $\nabla^B$, by (2.8), we have the following six equations:
\begin{eqnarray}
       \begin{cases}
       \alpha^2+\beta\gamma-\mu=0\\[2pt]
       \mu_2\alpha=0\\[2pt]
       \mu_3\delta=0\\[2pt]
       \mu_1\alpha+\alpha^2-\mu=0\\[2pt]
       \mu_1\gamma=0\\[2pt]
       \mu=0\\[2pt]
       \end{cases}
\end{eqnarray}
By solving (3.16), we get \\
\begin{thm}
$(G_6, g, V)$ is an affine Ricci soliton associated to the Bott connection $\nabla^B$ if and only if\\
\begin{eqnarray}
&&(1)\mu=\mu_1=\mu_3=\alpha=\beta=0,~~~ \delta\neq0;\nonumber\\
&&(2)\mu=\mu_3=\alpha=\beta=\gamma==0,~~~\mu_1\neq0,~~~\delta\neq0\nonumber.
\end{eqnarray}
\end{thm}
\vskip 0.5 true cm
\noindent{\bf 3.3 Affine Ricci solitons of $G_7$}\\
\vskip 0.5 true cm\indent By \cite{W}, we have the following Lie algebra of $G_7$ satisfies
\begin{equation}
[\widehat{e}_1, \widehat{e}_2]=-\alpha\widehat{e}_1-\beta\widehat{e}_2-\beta\widehat{e}_3,~~~[\widehat{e}_1, \widehat{e}_3]=\alpha\widehat{e}_1+\beta\widehat{e}_2+\beta\widehat{e}_3,~~~[\widehat{e}_2, \widehat{e}_3]=\gamma\widehat{e}_1+\delta\widehat{e}_2+\delta\widehat{e}_3,~~~\alpha+\delta\neq 0,~~~\alpha\gamma=0.
\end{equation}
where ${\widehat{e}_1,\widehat{e}_2,\widehat{e}_3}$ is a pseudo-orthonormal basis, with $\widehat{e}_3$ timelike.
\begin{lem}(\cite{G},\cite{W}) The Levi-Civita connection $\nabla^L$ of $G_7$ is given by
\begin{align}
&\nabla^L_{\widehat{e}_1}\widehat{e}_1=\alpha\widehat{e}_2+\alpha\widehat{e}_3,~~~\nabla^L_{\widehat{e}_1}\widehat{e}_2=-\alpha\widehat{e}_1+\frac{\gamma}{2}\widehat{e}_3,~~~\nabla^L_{\widehat{e}_1}\widehat{e}_3=\alpha\widehat{e}_1+\frac{\gamma}{2}\widehat{e}_2,\nonumber\\
&\nabla^L_{\widehat{e}_2}\widehat{e}_1=\beta\widehat{e}_2+(\beta+\frac{\gamma}{2})\widehat{e}_3,~~~\nabla^L_{\widehat{e}_2}\widehat{e}_2=-\beta\widehat{e}_1+\delta\widehat{e}_3,~~~\nabla^L_{\widehat{e}_2}\widehat{e}_3=(\beta+\frac{\gamma}{2})\widehat{e}_1+\delta\widehat{e}_2,\nonumber\\
&\nabla^L_{\widehat{e}_3}\widehat{e}_1=(\frac{\gamma}{2}-\beta)\widehat{e}_2-\beta\widehat{e}_3,~~~\nabla^L_{\widehat{e}_3}\widehat{e}_2=(\beta-\frac{\gamma}{2})\widehat{e}_1-\delta\widehat{e}_3,~~~\nabla^L_{\widehat{e}_3}\widehat{e}_3=-\beta\widehat{e}_1-\delta\widehat{e}_2.
\end{align}
\end{lem}
\begin{lem} The Bott connection $\nabla^B$ of $G_7$ is given by
\begin{align}
&\nabla^B_{\widehat{e}_1}\widehat{e}_1=\alpha\widehat{e}_2,~~~\nabla^B_{\widehat{e}_1}\widehat{e}_2=-\alpha\widehat{e}_1,~~~\nabla^B_{\widehat{e}_1}\widehat{e}_3=\beta\widehat{e}_3,\nonumber\\
&\nabla^B_{\widehat{e}_2}\widehat{e}_1=\beta\widehat{e}_2,~~~\nabla^B_{\widehat{e}_2}\widehat{e}_2=-\beta\widehat{e}_1,~~~\nabla^B_{\widehat{e}_2}\widehat{e}_3=\delta\widehat{e}_3,\nonumber\\
&\nabla^B_{\widehat{e}_3}\widehat{e}_1=-\alpha\widehat{e}_1-\beta\widehat{e}_2,~~~\nabla^B_{\widehat{e}_3}\widehat{e}_2=-\gamma\widehat{e}_1-\delta\widehat{e}_2,~~~\nabla^B_{\widehat{e}_3}\widehat{e}_3=0.
\end{align}
\end{lem}
\begin{lem} The curvature $R^B$ of the Bott connection $\nabla^B$ of $(G_7,g)$ is given by
\begin{align}
&R^B(\widehat{e}_1,\widehat{e}_2)\widehat{e}_1=-\alpha\beta\widehat{e}_1+\alpha^2\widehat{e}_2,~~~R^B(\widehat{e}_1,\widehat{e}_2)\widehat{e}_2=-(\alpha^2+\beta^2+\beta\gamma)\widehat{e}_1-\beta\delta\widehat{e}_2,~~~R^B(\widehat{e}_1,\widehat{e}_2)\widehat{e}_3=\beta(\alpha-\delta)\widehat{e}_3,\nonumber\\
&R^B(\widehat{e}_1,\widehat{e}_3)\widehat{e}_1=\alpha(2\beta+\gamma)\widehat{e}_1+(\alpha\delta-2\alpha^2)\widehat{e}_2,~~~R^B(\widehat{e}_1,\widehat{e}_3)\widehat{e}_2=(\alpha\delta+\beta^2+\beta\gamma)\widehat{e}_1+(\beta\delta-\alpha\beta-\alpha\gamma)\widehat{e}_2,\nonumber\\
&R^B(\widehat{e}_1,\widehat{e}_3)\widehat{e}_3=-\beta(\alpha+\delta)\widehat{e}_3,~~~R^B(\widehat{e}_2,\widehat{e}_3)\widehat{e}_1=(\beta^2+\beta\gamma+\alpha\delta)\widehat{e}_1+(\beta\delta-\alpha\beta-\alpha\gamma)\widehat{e}_2,\nonumber\\
&R^B(\widehat{e}_2,\widehat{e}_3)\widehat{e}_2=(2\beta\delta+\delta\gamma+\alpha\gamma-\alpha\beta)\widehat{e}_1+(\delta^2-\beta^2-\beta\gamma)\widehat{e}_2,~~~R^B(\widehat{e}_2,\widehat{e}_3)\widehat{e}_3=-(\beta\gamma+\delta^2)\widehat{e}_3.
\end{align}
\end{lem}
By (2.5), we have
\begin{align}
&\rho^B(\widehat{e}_1,\widehat{e}_1)=-\alpha^2,~~~\rho^B(\widehat{e}_1,\widehat{e}_2)=\beta\delta,~~~\rho^B(\widehat{e}_1,\widehat{e}_3)=\beta(\alpha+\delta),\nonumber\\
&\rho^B(\widehat{e}_2,\widehat{e}_1)=-\alpha\beta,~~~\rho^B(\widehat{e}_2,\widehat{e}_2)=-(\alpha^2+\beta^2+\beta\gamma),~~~\rho^B(\widehat{e}_2,\widehat{e}_3)=(\beta\gamma+\delta^2),\nonumber\\
&\rho^B(\widehat{e}_3,\widehat{e}_1)=\beta(\alpha+\delta),~~~\rho^B(\widehat{e}_3,\widehat{e}_2)=\delta(\alpha+\delta),~~~\rho^B(\widehat{e}_3,\widehat{e}_3)=0.
\end{align}
Then,
\begin{align}
&\widetilde{\rho}^B(\widehat{e}_1,\widehat{e}_1)=-\alpha^2,~~~\widetilde{\rho}^B(\widehat{e}_1,\widehat{e}_2)=\frac{\beta(\delta-\alpha)}{2},~~~\widetilde{\rho}^B(\widehat{e}_1,\widehat{e}_3)=\delta(\alpha+\delta),\nonumber\\
&\widetilde{\rho}^B(\widehat{e}_2,\widehat{e}_2)=-(\alpha^2+\beta^2+\beta\gamma),~~~\widetilde{\rho}^B(\widehat{e}_2,\widehat{e}_3)=\delta^2+\frac{\beta\gamma+\alpha\delta}{2},~~~\widetilde{\rho}^B(\widehat{e}_3,\widehat{e}_3)=0.
\end{align}
By (2.7), we have
\begin{align}
&(L^B_Vg)(\widehat{e}_1,\widehat{e}_1)=-2\mu_2\alpha,~~~(L^B_Vg)(\widehat{e}_1,\widehat{e}_2)=\mu_1\alpha-\mu_2\beta,~~~(L^B_Vg)(\widehat{e}_1,\widehat{e}_3)=-\mu_1\alpha-\mu_2\gamma-\mu_3\beta,\nonumber\\
&(L^B_Vg)(\widehat{e}_2,\widehat{e}_2)=2\mu_1\beta,~~~(L^B_Vg)(\widehat{e}_2,\widehat{e}_3)=-\mu_1\beta-\mu_2\delta-\mu_3\delta,~~~(L^B_Vg)(\widehat{e}_3,\widehat{e}_3)=0.
\end{align}
Then, if $(G_7,g,V)$  is an affine Ricci soliton associated to the Bott connection $\nabla^B$, by (2.8), we have the following six equations:
\begin{eqnarray}
       \begin{cases}
       \mu_2\alpha+\alpha^2-\mu=0\\[2pt]
       \mu_1\alpha-\mu_2\beta+\beta\delta-\alpha\beta=0\\[2pt]
       \mu_3\beta+\mu_1\alpha+\mu_2\gamma-2\alpha\beta-2\delta\beta=0\\[2pt]
       \mu_1\beta-\alpha^2-\beta^2-\beta\gamma+\mu=0\\[2pt]
       \mu_3\delta+\mu_1\beta+\mu_2\delta-2\delta^2-\beta\gamma-\alpha\delta=0\\[2pt]
       \mu=0\\[2pt]
       \end{cases}
\end{eqnarray}
By solving (3.24), we get \\
\begin{thm}
$(G_7, g, V)$ is an affine Ricci soliton associated to the Bott connection $\nabla^B$ if and only if\\
\begin{eqnarray*}
&&(1)\mu=\alpha=\beta=\gamma=0,~~~ \delta\neq 0, ~~~\mu_2+\mu_3-2\delta=0;\nonumber\\
&&(2)\mu=\mu_2=\alpha=\beta=0,~~~\gamma\neq0,~~~\delta\neq0,~~~\mu_3-2\delta=0;\nonumber\\
&&(3)\mu=\alpha=0,~~~\beta\neq0,~~~\delta\neq0,~~~\mu_1=\gamma+\beta,~~~\mu_2=\delta,~~~\mu_3=\frac{\delta(2\beta-\gamma)}{\beta},~~~\gamma=\frac{\beta(\beta^2+\delta^2)}{\delta^2}.\nonumber\\
\end{eqnarray*}
\end{thm}

\indent Specially, let $V=0$, we get the following corollary:
\begin{cor}
(I)~~$(G_1, g, V)$ is not an affine Einstein associated to the Bott connection $\nabla^B$;\nonumber\\
(II)~~$(G_2, g, V)$ is not an affine Einstein associated to the Bott connection $\nabla^B$;\nonumber\\
(III)~~$(G_3, g, V)$ is an affine Einstein associated to the Bott connection $\nabla^B$ if and only if
$(1)\mu=\alpha=\beta=0,~~~\gamma\neq0;$
$(2)\mu=\gamma=0$;\nonumber\\
(IV)~~$(G_4 ,g, V)$ is not an affine Einstein associated to the Bott connection $\nabla^B$;\nonumber\\
(V)~~$(G_5, g, V)$ is an affine Einstein associated to the Bott connection $\nabla^B$ if and only if $\mu=0,~~~\alpha+\delta\neq0,~~~\alpha\gamma+\beta\delta=0$;\nonumber\\
(VI)~~$(G_6, g, V)$ is an affine Einstein associated to the Bott connection $\nabla^B$ if and only if $\mu=\alpha=\beta=0,~~~\delta\neq0$;\nonumber\\
(VII)~~$(G_7, g, V)$ is not an affine Einstein associated to the Bott connection $\nabla^B$.\nonumber\\
\end{cor}
\section{Affine Ricci Solitons associated to the perturbed Bott connection on three-dimensional Lorentzian Lie groups with the first distribution}
By the above calculations, we always obtain $\mu=0$. In order to get the affine Ricci soliton with non zero $\mu$, we introduce the perturbed Bott connection $\widetilde{\nabla}^B$ in the following. Let $\widehat{e}_3^*$ be the dual base of $e_3$. We define on $G_{i=1,\cdot\cdot\cdot,7}$
\begin{align}
\widetilde{\nabla}^B_XY=\nabla^B_XY+a_0\widehat{e}_3^*(X)\widehat{e}_3^*(Y)e_3,
\end{align}
where $a_0$ is a non zero real number. Then
\begin{align}
\widetilde{\nabla}^B_{\widehat{e}_3}\widehat{e}_3=a_0\widehat{e}_3,~~~\widetilde{\nabla}^B_{\widehat{e}_s}\widehat{e}_t=\nabla^B_{\widehat{e}_s}\widehat{e}_t,
\end{align}
where $s$ and $t$ does not equal 3. We define
\begin{align}
(\widetilde{L}_V^Bg)(X,Y):=g(\widetilde{\nabla}^B_XV,Y)+g(X,\widetilde{\nabla}^B_YV),
\end{align}
for vector fields $X,Y,V$. Then we have for $G_{i=1,\cdot\cdot\cdot,7}$
\begin{align}
(\widetilde{L}_V^Bg)(\widehat{e}_3,\widehat{e}_3)=-2a_0\mu_3,~~~(\widetilde{L}_V^Bg)(\widehat{e}_s,\widehat{e}_t)=(L_V^Bg)(\widehat{e}_s,\widehat{e}_t),
\end{align}
where $s$ and $t$ does not equal 3.
\begin{defn}$(G_i,V,g)$ is called the affine Ricci soliton associated to the connection $\widetilde{\nabla}^B$ if it satisfies
\begin{align}
(\widetilde{L}_V^Bg)(X,Y)+2\widetilde{\overline{\rho}}^B(X,Y)+2\mu g(X,Y)=0.
\end{align}
\end{defn}
For $(G_1,\widetilde{\nabla}^B)$, we have
\begin{align}
\widetilde{R}^B(\widehat{e}_1,\widehat{e}_2)\widehat{e}_3=a_0\beta \widehat{e}_3,~~~\widetilde{R}^B(\widehat{e}_2,\widehat{e}_3)\widehat{e}_3=-\alpha(a_0+\alpha)\widehat{e}_3,~~~\widetilde{R}^B(\widehat{e}_s,\widehat{e}_t)\widehat{e}_p=R^B(\widehat{e}_s,\widehat{e}_t)\widehat{e}_p,
\end{align}
for $(s,t,p)\neq(1,2,3),(2,3,3)$.\\
\begin{align}
\widetilde{\overline{\rho}}^B(\widehat{e}_2,\widehat{e}_3)=\frac{\alpha(a_0+\alpha)}{2},~~~\widetilde{\overline{\rho}}^B(\widehat{e}_s,\widehat{e}_t)=\widetilde{\rho}^B(\widehat{e}_s,\widehat{e}_t),
\end{align}
for the pair $(s,t)\neq(2,3)$. If $(G_1,g,V)$ is an affine Ricci soliton associated to the connection $\widetilde{\nabla}^B$, then by (4.5), we have
\begin{eqnarray}
       \begin{cases}
        \mu_2\alpha-\alpha^2-\beta^2+\mu=0 \\[2pt]
       2\alpha\beta-\mu_1\alpha=0\\[2pt]
       \mu_1\alpha-\mu_2\beta-\alpha\beta=0\\[2pt]
       \alpha^2+\beta^2-\mu=0\\[2pt]
       (\mu_2+\mu_3)\alpha-\mu_1\beta-a_0\alpha-\alpha^2=0\\[2pt]
      a_0\mu_3+\mu=0\\[2pt]
       \end{cases}
\end{eqnarray}
Solve (4.8), we get $a_0=\alpha=0$, there is a contradiction. So\\
\begin{thm}
$(G_1,V,g)$ is not an affine Ricci soliton associated to the connection $\widetilde{\nabla}^B$.
\end{thm}
For $(G_2,\widetilde{\nabla}^B)$, we have
\begin{align}
\widetilde{R}^B(\widehat{e}_1,\widehat{e}_2)\widehat{e}_3=a_0\beta \widehat{e}_3,~~~\widetilde{R}^B(\widehat{e}_1,\widehat{e}_3)\widehat{e}_3=a_0\gamma\widehat{e}_3,~~~\widetilde{R}^B(\widehat{e}_s,\widehat{e}_t)\widehat{e}_p=R^B(\widehat{e}_s,\widehat{e}_t)\widehat{e}_p,
\end{align}
for $(s,t,p)\neq(1,2,3),(1,3,3)$.\\
\begin{align}
\widetilde{\overline{\rho}}^B(\widehat{e}_1,\widehat{e}_3)=\frac{a_0\gamma}{2},~~~\widetilde{\overline{\rho}}^B(\widehat{e}_s,\widehat{e}_t)=\widetilde{\rho}^B(\widehat{e}_s,\widehat{e}_t),
\end{align}
for the pair $(s,t)\neq(1,3)$. If $(G_2,g,V)$ is an affine Ricci soliton associated to the connection $\widetilde{\nabla}^B$, then by (4.5), we have
\begin{eqnarray}
       \begin{cases}
       \gamma^2+\beta^2-\mu=0 \\[2pt]
      \mu_2\gamma=0\\[2pt]
       \mu_2\alpha-\mu_3\gamma+a_0\gamma=0\\[2pt]
       \gamma^2+\mu_1\gamma+\alpha\beta-\mu=0\\[2pt]
       \mu_1\beta-\alpha\gamma=0\\[2pt]
      a_0\mu_3+\mu=0\\[2pt]
       \end{cases}
\end{eqnarray}
Solve (4.11), we get $a_0=\beta=\gamma=0$, there is a contradiction. So\\
\begin{thm}
$(G_2,V,g)$ is not an affine Ricci soliton associated to the connection $\widetilde{\nabla}^B$.
\end{thm}
For $(G_3,\widetilde{\nabla}^B)$, we have
\begin{align}
\widetilde{R}^B(\widehat{e}_1,\widehat{e}_2)\widehat{e}_3=-a_0\beta \widehat{e}_3,~~~\widetilde{R}^B(\widehat{e}_s,\widehat{e}_t)\widehat{e}_p=R^B(\widehat{e}_s,\widehat{e}_t)\widehat{e}_p,
\end{align}
for $(s,t,p)\neq(1,2,3)$.\\
\begin{align}
\widetilde{\overline{\rho}}^B(\widehat{e}_s,\widehat{e}_t)=\widetilde{\rho}^B(\widehat{e}_s,\widehat{e}_t),
\end{align}
for any pairs $(s,t)$. If $(G_3,g,V)$ is an affine Ricci soliton associated to the connection $\widetilde{\nabla}^B$, then by (4.5), we have
\begin{eqnarray}
       \begin{cases}
       \beta\gamma-\mu=0 \\[2pt]
      \mu_2\alpha=0\\[2pt]
       \mu-\alpha\gamma=0\\[2pt]
       \mu_1\beta=0\\[2pt]
      a_0\mu_3+\mu=0\\[2pt]
       \end{cases}
\end{eqnarray}
Solve (4.14), we get
\begin{thm}
$(G_3,V,g)$ is an affine Ricci soliton associated to the connection $\widetilde{\nabla}^B$ if and only if
\begin{eqnarray*}
&&(1)\mu=\gamma=\mu_3=\mu_2\alpha=\mu_1\beta=0;\nonumber\\
&&(2)\gamma\neq 0,~~~\alpha=\beta=\mu=\mu_3=0;\nonumber\\
&&(3)\gamma\neq0,~~~\alpha=\beta\neq0,~~~\mu_1=\mu_2=0,~~~\mu=\alpha\gamma,~~~\mu_3=-\frac{\alpha\gamma}{a_0}.
\end{eqnarray*}
\end{thm}
For $(G_4,\widetilde{\nabla}^B)$, we have
\begin{align}
\widetilde{R}^B(\widehat{e}_1,\widehat{e}_2)\widehat{e}_3=a_0(\beta-2\eta) \widehat{e}_3,~~~\widetilde{R}^B(\widehat{e}_1,\widehat{e}_3)\widehat{e}_3=-a_0\widehat{e}_3,~~~\widetilde{R}^B(\widehat{e}_s,\widehat{e}_t)\widehat{e}_p=R^B(\widehat{e}_s,\widehat{e}_t)\widehat{e}_p,
\end{align}
for $(s,t,p)\neq(1,2,3),(1,3,3)$.\\
\begin{align}
\widetilde{\overline{\rho}}^B(\widehat{e}_1,\widehat{e}_3)=\frac{a_0}{2},~~~\widetilde{\overline{\rho}}^B(\widehat{e}_s,\widehat{e}_t)=\widetilde{\rho}^B(\widehat{e}_s,\widehat{e}_t),
\end{align}
for the pair $(s,t)\neq(1,3)$. If $(G_4,g,V)$ is an affine Ricci soliton associated to the connection $\widetilde{\nabla}^B$, then by (4.5), we have
\begin{eqnarray}
       \begin{cases}
       \beta^2-2\beta\eta+1-\mu=0 \\[2pt]
      \mu_2=0\\[2pt]
       \mu_2\alpha+\mu_3-a_0=0\\[2pt]
       \mu_1+2\alpha\gamma-\alpha\beta-1+\mu=0\\[2pt]
       \mu_1\beta+\alpha=0\\[2pt]
      a_0\mu_3+\mu=0\\[2pt]
       \end{cases}
\end{eqnarray}
Solve (4.17), we get $a_0=0$, there is a contradiction. So\\
\begin{thm}
$(G_4,V,g)$ is not an affine Ricci soliton associated to the connection $\widetilde{\nabla}^B$.
\end{thm}

For $(G_5,\widetilde{\nabla}^B)$, we have
\begin{align}
\widetilde{R}^B(\widehat{e}_s,\widehat{e}_t)\widehat{e}_p=R^B(\widehat{e}_s,\widehat{e}_t)\widehat{e}_p,
\end{align}
for any $(s,t,p)$.\\
\begin{align}
\widetilde{\overline{\rho}}^B(\widehat{e}_s,\widehat{e}_t)=\widetilde{\rho}^B(\widehat{e}_s,\widehat{e}_t),
\end{align}
for any pairs $(s,t)$. If $(G_5,g,V)$ is an affine Ricci soliton associated to the connection $\widetilde{\nabla}^B$, then by (4.5), we have
\begin{eqnarray}
       \begin{cases}
       \mu=0 \\[2pt]
       \mu_1\alpha+\mu_2\gamma=0\\[2pt]
       \mu_1\beta+\mu_2\delta=0\\[2pt]
      a_0\mu_3+\mu=0\\[2pt]
       \end{cases}
\end{eqnarray}
Solve (4.20), we get
\begin{thm}
$(G_5,V,g)$ is an affine Ricci soliton associated to the connection $\widetilde{\nabla}^B$ if and only if
\begin{eqnarray*}
&&(1)\mu=\mu_1=\mu_2=\mu_3=0,~~~\alpha+\delta\neq0,~~~\alpha\gamma+\beta\delta=0;\nonumber\\
&&(2)\mu=\mu_2=\mu_3=\alpha=\beta=0,~~~\mu_1\neq0,~~~\delta\neq0;\nonumber\\
&&(3)\mu=\mu_1=\mu_3=\delta=\gamma=0=0,~~~\mu_2\neq0,~~~\alpha\neq0.
\end{eqnarray*}
\end{thm}
For $(G_6,\widetilde{\nabla}^B)$, we have
\begin{align}
\widetilde{R}^B(\widehat{e}_1,\widehat{e}_2)\widehat{e}_3=a_0\gamma \widehat{e}_3,~~~\widetilde{R}^B(\widehat{e}_1,\widehat{e}_3)\widehat{e}_3=-a_0\delta\widehat{e}_3,~~~\widetilde{R}^B(\widehat{e}_s,\widehat{e}_t)\widehat{e}_p=R^B(\widehat{e}_s,\widehat{e}_t)\widehat{e}_p,
\end{align}
for $(s,t,p)\neq(1,2,3),(1,3,3)$.\\
\begin{align}
\widetilde{\overline{\rho}}^B(\widehat{e}_1,\widehat{e}_3)=\frac{a_0\delta}{2},~~~\widetilde{\overline{\rho}}^B(\widehat{e}_s,\widehat{e}_t)=\widetilde{\rho}^B(\widehat{e}_s,\widehat{e}_t),
\end{align}
for the pair $(s,t)\neq(1,3)$. If $(G_6,g,V)$ is an affine Ricci soliton associated to the connection $\widetilde{\nabla}^B$, then by (4.5), we have
\begin{eqnarray}
       \begin{cases}
       \alpha^2+\beta\gamma-\mu=0 \\[2pt]
      \mu_2\alpha=0\\[2pt]
       \mu_3\delta-a_0\delta=0\\[2pt]
       \mu-\mu_1\alpha-\alpha^2=0\\[2pt]
       \mu_1\gamma=0\\[2pt]
      a_0\mu_3+\mu=0\\[2pt]
       \end{cases}
\end{eqnarray}
Solve (4.23), we get
\begin{thm}
$(G_6,V,g)$ is an affine Ricci soliton associated to the connection $\widetilde{\nabla}^B$ if and only if
$\alpha\neq0,~~~\mu_1=\mu_2=\gamma=\delta=0,~~~\mu=\alpha^2,~~~\mu_3=-\frac{\alpha^2}{a_0}$.
\end{thm}
For $(G_7,\widetilde{\nabla}^B)$, we have
\begin{align}
&\widetilde{R}^B(\widehat{e}_1,\widehat{e}_2)\widehat{e}_3=\beta(\alpha-\delta+a_0) \widehat{e}_3,~~~\widetilde{R}^B(\widehat{e}_1,\widehat{e}_3)\widehat{e}_3=-\beta(a_0+\alpha+\delta)\widehat{e}_3,\nonumber\\
&\widetilde{R}^B(\widehat{e}_2,\widehat{e}_3)\widehat{e}_3=-(a_0\delta+\delta^2+\beta\gamma)\widehat{e}_3,~~~\widetilde{R}^B(\widehat{e}_s,\widehat{e}_t)\widehat{e}_p=R^B(\widehat{e}_s,\widehat{e}_t)\widehat{e}_p,
\end{align}
for $(s,t,p)\neq(1,2,3),(1,3,3),(2,3,3)$.\\
\begin{align}
\widetilde{\overline{\rho}}^B(\widehat{e}_1,\widehat{e}_3)=\frac{\beta(a_0+2\alpha+2\delta)}{2},~~~\widetilde{\overline{\rho}}^B(\widehat{e}_2,\widehat{e}_3)=\frac{a_0\delta+2\delta^2+\alpha\delta+\beta\gamma}{2},~~~\widetilde{\overline{\rho}}^B(\widehat{e}_s,\widehat{e}_t)=\widetilde{\rho}^B(\widehat{e}_s,\widehat{e}_t),
\end{align}
for the pair $(s,t)\neq(1,3),(2,3)$. If $(G_7,g,V)$ is an affine Ricci soliton associated to the connection $\widetilde{\nabla}^B$, then by (4.5), we have
\begin{eqnarray}
       \begin{cases}
        \mu_2\alpha+\alpha^2-\mu=0 \\[2pt]
       \mu_2\beta-\mu_1\alpha-\beta\delta+\alpha\beta=0\\[2pt]
       \mu_1\alpha+\mu_2\gamma+\mu_3\beta-(2\alpha\beta+2\beta\delta+a_0\beta)=0\\[2pt]
       \mu_1\beta-\alpha^2-\beta^2-\beta\gamma-\mu=0\\[2pt]
       (\mu_2+\mu_3)\delta+\mu_1\beta-(\beta\gamma+\alpha\delta+2\delta^2+a_0\delta)=0\\[2pt]
      a_0\mu_3+\mu=0\\[2pt]
       \end{cases}
\end{eqnarray}
Solve (4.26), we get
\begin{thm}
$(G_7,V,g)$ is an affine Ricci soliton associated to the connection $\widetilde{\nabla}^B$ if and only if
\begin{eqnarray*}
&&(1)\mu=\alpha=\beta=\gamma=\mu_3=0,~~~\delta\neq0,~~~\mu_2=2\delta+a_0;\nonumber\\
&&(2)\alpha=\beta=\mu=\mu_2=\mu_3=0,~~~\gamma\neq0,~~~\delta\neq0,~~~a_0=-2\delta;\nonumber\\
&&(3)\alpha=\mu=\mu_3=0,~~~\beta\neq0,~~~\delta\neq0,~~~\mu_1=\beta+\gamma,~~~\mu_2=\delta,~~~a_0=\frac{\delta(\gamma-2\beta)}{\beta},~~~\gamma=\frac{\beta(\beta^2+\delta^2)}{\delta^2};\nonumber\\
&&(4)\alpha\neq0,~~~\mu_1=\mu_2=\beta=\gamma=\delta=0,~~~\mu=\alpha^2,~~~\mu_3=-\frac{\alpha^2}{a_0};\nonumber\\
&&(5)\alpha\neq0,~~~\mu_1=\mu_2=\beta=\gamma=0,~~~\mu=\alpha^2,~~~\mu_3=\alpha+2\delta+a_0,~~~a_0^2+a_0\alpha+2a_0\delta+\alpha^2=0.
\end{eqnarray*}
\end{thm}
\section{Affine Ricci Solitons associated to the Bott connection on three-dimensional Lorentzian Lie groups with the second distribution}
Let $TM=span\{\widehat{e}_1,\widehat{e}_2,\widehat{e}_3\}$, then took the distribution $D_1=span\{\widehat{e}_1,\widehat{e}_3\}$ and $D_1^\bot=span\{\widehat{e}_2\}$.\\
Similar to (2.3), we have\\
\begin{eqnarray}
\nabla^{B_1}_XY=
       \begin{cases}
        \pi_{D_1}(\nabla^L_XY),~~~&X,Y\in\Gamma^\infty({D_1}) \\[2pt]
       \pi_{D_1}([X,Y]),~~~&X\in\Gamma^\infty({D_1}^\bot),Y\in\Gamma^\infty({D_1})\\[2pt]
       \pi_{{D_1}^\bot}([X,Y]),~~~&X\in\Gamma^\infty({D_1}),Y\in\Gamma^\infty({D_1}^\bot)\\[2pt]
       \pi_{{D_1}^\bot}(\nabla^L_XY),~~~&X,Y\in\Gamma^\infty({D_1}^\bot)\\[2pt]
       \end{cases}
\end{eqnarray}
where $\pi_{D_1}$(resp. $\pi_{D_1}^\bot$) the projection on ${D_1}$ (resp. ${D_1}^\bot$).\\
\vskip 0.5 true cm
\noindent{\bf 5.1 Affine Ricci solitons of $G_1$}\\
\vskip 0.5 true cm
\begin{lem} The Bott connection $\nabla^{B_1}$ of $G_1$ is given by
\begin{align}
&\nabla^{B_1}_{\widehat{e}_1}\widehat{e}_1=-\alpha\widehat{e}_3,~~~\nabla^{B_1}_{\widehat{e}_1}\widehat{e}_2=0,~~~\nabla^{B_1}_{\widehat{e}_1}\widehat{e}_3=-\alpha\widehat{e}_1,\nonumber\\
&\nabla^{B_1}_{\widehat{e}_2}\widehat{e}_1=-\alpha\widehat{e}_1+\beta\widehat{e}_3,~~~\nabla^{B_1}_{\widehat{e}_2}\widehat{e}_2=0,~~~\nabla^{B_1}_{\widehat{e}_2}\widehat{e}_3=\beta\widehat{e}_1+\alpha\widehat{e}_3,\nonumber\\
&\nabla^{B_1}_{\widehat{e}_3}\widehat{e}_1=0,~~~\nabla^{B_1}_{\widehat{e}_3}\widehat{e}_2=-\alpha\widehat{e}_2,~~~\nabla^{B_1}_{\widehat{e}_3}\widehat{e}_3=0.
\end{align}
\end{lem}
\begin{lem} The curvature $R^{B_1}$ of the Bott connection $\nabla^{B_1}$ of $(G_1,g)$ is given by
\begin{align}
&R^{B_1}(\widehat{e}_1,\widehat{e}_2)\widehat{e}_1=3\alpha^2\widehat{e}_3,~~~R^{B_1}(\widehat{e}_1,\widehat{e}_2)\widehat{e}_2=-\alpha\beta\widehat{e}_2,~~~R^{B_1}(\widehat{e}_1,\widehat{e}_2)\widehat{e}_3=-\alpha^2\widehat{e}_1,\nonumber\\
&R^{B_1}(\widehat{e}_1,\widehat{e}_3)\widehat{e}_1=-\alpha\beta\widehat{e}_1+(\beta^2-\alpha^2)\widehat{e}_3,~~~R^{B_1}(\widehat{e}_1,\widehat{e}_3)\widehat{e}_2=0,~~~R^{B_1}(\widehat{e}_1,\widehat{e}_3)\widehat{e}_3=(\beta^2-\alpha^2)\widehat{e}_1+\alpha\beta\widehat{e}_3,\nonumber\\
&R^{B_1}(\widehat{e}_2,\widehat{e}_3)\widehat{e}_1=\alpha^2\widehat{e}_1,~~~R^{B_1}(\widehat{e}_2,\widehat{e}_3)\widehat{e}_2=\alpha^2\widehat{e}_2,~~~R^{B_1}(\widehat{e}_2,\widehat{e}_3)\widehat{e}_3=-\alpha^2\widehat{e}_3.
\end{align}
\end{lem}
By (2.5), we have
\begin{align}
&\rho^{B_1}(\widehat{e}_1,\widehat{e}_1)=\alpha^2-\beta^2,~~~\rho^{B_1}(\widehat{e}_1,\widehat{e}_2)=\alpha\beta,~~~\rho^{B_1}(\widehat{e}_1,\widehat{e}_3)=-\alpha\beta,\nonumber\\
&\rho^{B_1}(\widehat{e}_2,\widehat{e}_1)=\rho^{B_1}(\widehat{e}_2,\widehat{e}_2)=\rho^{B_1}(\widehat{e}_2,\widehat{e}_3)=0,\nonumber\\
&\rho^{B_1}(\widehat{e}_3,\widehat{e}_1)=-\alpha\beta,~~~\rho^{B_1}(\widehat{e}_3,\widehat{e}_2)=\alpha^2,~~~\rho^{B_1}(\widehat{e}_3,\widehat{e}_3)=(\beta^2-\alpha^2).
\end{align}
Then,
\begin{align}
&\widetilde{\rho}^{B_1}(\widehat{e}_1,\widehat{e}_1)=\alpha^2-\beta^2,~~~\widetilde{\rho}^{B_1}(\widehat{e}_1,\widehat{e}_2)=\frac{\alpha\beta}{2},~~~\widetilde{\rho}^{B_1}(\widehat{e}_1,\widehat{e}_3)=-\alpha\beta,\nonumber\\
&\widetilde{\rho}^{B_1}(\widehat{e}_2,\widehat{e}_2)=0,~~~\widetilde{\rho}^{B_1}(\widehat{e}_2,\widehat{e}_3)=\frac{\alpha^2}{2},~~~\widetilde{\rho}^{B_1}(\widehat{e}_3,\widehat{e}_3)=(\beta^2-\alpha^2).
\end{align}
By (2.7), we have
\begin{align}
&(L^{B_1}_Vg)(\widehat{e}_1,\widehat{e}_1)=-2\mu_3\alpha,~~~(L^{B_1}_Vg)(\widehat{e}_1,\widehat{e}_2)=-\mu_1\alpha+\mu_3\beta,~~~(L^{B_1}_Vg)(\widehat{e}_1,\widehat{e}_3)=\mu_1\alpha,\nonumber\\
&(L^{B_1}_Vg)(\widehat{e}_2,\widehat{e}_2)=0,~~~(L^{B_1}_Vg)(\widehat{e}_2,\widehat{e}_2)=-(\mu_2+\mu_3)\alpha-\mu_1\beta,~~~(L^{B_1}_Vg)(\widehat{e}_3,\widehat{e}_3)=0.
\end{align}
Then, if $(G_1,g,V)$  is an affine Ricci soliton associated to the Bott connection $\nabla^{B_1}$, by (2.8), we have the following six equations:
\begin{eqnarray}
       \begin{cases}
        \mu_3\alpha-\alpha^2+\beta^2-\mu=0 \\[2pt]
      \mu_3\beta-\mu_1\alpha+\alpha\beta=0\\[2pt]
       \mu_1\alpha-2\alpha\beta=0\\[2pt]
       \mu=0\\[2pt]
       \mu_1\beta+\mu_3\alpha+\mu_2\alpha-\alpha^2=0\\[2pt]
       \beta^2-\alpha^2-\mu=0\\[2pt]
       \end{cases}
\end{eqnarray}
By solving (5.7) , we get $\alpha=0$, there is a contradiction. So\\
\begin{thm}
$(G_1, g, V)$ is not an affine Ricci soliton associated to the Bott connection $\nabla^{B_1}$.
\end{thm}

\vskip 0.5 true cm
\noindent{\bf 5.2 Affine Ricci solitons of $G_2$}\\
\vskip 0.5 true cm
\begin{lem} The Bott connection $\nabla^{B_1}$ of $G_2$ is given by
\begin{align}
&\nabla^{B_1}_{\widehat{e}_1}\widehat{e}_1=0,~~~\nabla^{B_1}_{\widehat{e}_1}\widehat{e}_2=\gamma\widehat{e}_2,~~~\nabla^{B_1}_{\widehat{e}_1}\widehat{e}_3=0,\nonumber\\
&\nabla^{B_1}_{\widehat{e}_2}\widehat{e}_1=\beta\widehat{e}_3,~~~\nabla^{B_1}_{\widehat{e}_2}\widehat{e}_2=0,~~~\nabla^{B_1}_{\widehat{e}_2}\widehat{e}_3=\alpha\widehat{e}_1,\nonumber\\
&\nabla^{B_1}_{\widehat{e}_3}\widehat{e}_1=\gamma\widehat{e}_3,~~~\nabla^{B_1}_{\widehat{e}_3}\widehat{e}_2=0,~~~\nabla^{B_1}_{\widehat{e}_3}\widehat{e}_3=\gamma\widehat{e}_1.
\end{align}
\end{lem}
\begin{lem} The curvature $R^{B_1}$ of the Bott connection $\nabla^{B_1}$ of $(G_2,g)$ is given by
\begin{align}
&R^{B_1}(\widehat{e}_1,\widehat{e}_2)\widehat{e}_1=0,~~~R^{B_1}(\widehat{e}_1,\widehat{e}_2)\widehat{e}_2=0,~~~R^{B_1}(\widehat{e}_1,\widehat{e}_2)\widehat{e}_3=\gamma(\beta-\alpha)\widehat{e}_1,\nonumber\\
&R^{B_1}(\widehat{e}_1,\widehat{e}_3)\widehat{e}_1=(\beta^2+\gamma^2)\widehat{e}_3,~~~R^{B_1}(\widehat{e}_1,\widehat{e}_3)\widehat{e}_2=0,~~~R^{B_1}(\widehat{e}_1,\widehat{e}_3)\widehat{e}_3=(\alpha\beta+\gamma^2)\widehat{e}_1,\nonumber\\
&R^{B_1}(\widehat{e}_2,\widehat{e}_3)\widehat{e}_1=\gamma(\alpha-\beta)\widehat{e}_1,~~~R^{B_1}(\widehat{e}_2,\widehat{e}_3)\widehat{e}_2=-\alpha\gamma\widehat{e}_2,~~~R^{B_1}(\widehat{e}_2,\widehat{e}_3)\widehat{e}_3=\gamma(\beta-\alpha)\widehat{e}_3.
\end{align}
\end{lem}
By (2.5), we have
\begin{align}
&\rho^{B_1}(\widehat{e}_1,\widehat{e}_1)=-(\beta^2+\gamma^2),~~~\rho^{B_1}(\widehat{e}_1,\widehat{e}_2)=0,~~~\rho^{B_1}(\widehat{e}_1,\widehat{e}_3)=0,\nonumber\\
&\rho^{B_1}(\widehat{e}_2,\widehat{e}_1)=\rho^{B_1}(\widehat{e}_2,\widehat{e}_2)=\rho^{B_1}(\widehat{e}_2,\widehat{e}_3)=0,\nonumber\\
&\rho^{B_1}(\widehat{e}_3,\widehat{e}_1)=0,~~~\rho^{B_1}(\widehat{e}_3,\widehat{e}_2)=-\alpha\gamma,~~~\rho^{B_1}(\widehat{e}_3,\widehat{e}_3)=\alpha\beta+\gamma^2.
\end{align}
Then,
\begin{align}
&\widetilde{\rho}^{B_1}(\widehat{e}_1,\widehat{e}_1)=-(\beta^2+\gamma^2),~~~\widetilde{\rho}^{B_1}(\widehat{e}_1,\widehat{e}_2)=0,~~~\widetilde{\rho}^{B_1}(\widehat{e}_1,\widehat{e}_3)=0,\nonumber\\
&\widetilde{\rho}^{B_1}(\widehat{e}_2,\widehat{e}_2)=0,~~~\widetilde{\rho}^{B_1}(\widehat{e}_2,\widehat{e}_3)=-\frac{\alpha\gamma}{2},~~~\widetilde{\rho}^{B_1}(\widehat{e}_3,\widehat{e}_3)=\alpha\beta+\gamma^2.
\end{align}
By (2.7), we have
\begin{align}
&(L^{B_1}_Vg)(\widehat{e}_1,\widehat{e}_1)=0,~~~(L^{B_1}_Vg)(\widehat{e}_1,\widehat{e}_2)=\mu_2\gamma+\mu_3\alpha,~~~(L^{B_1}_Vg)(\widehat{e}_1,\widehat{e}_3)=\mu_3\gamma,\nonumber\\
&(L^{B_1}_Vg)(\widehat{e}_2,\widehat{e}_2)=0,~~~(L^{B_1}_Vg)(\widehat{e}_2,\widehat{e}_3)=-\mu_1\beta,~~~(L^{B_1}_Vg)(\widehat{e}_3,\widehat{e}_3)=-2\mu_1\gamma.
\end{align}
Then, if $(G_2,g,V)$  is an affine Ricci soliton associated to the Bott connection $\nabla^{B_1}$, by (2.8), we have the following six equations:
\begin{eqnarray}
       \begin{cases}
        \beta^2+\gamma^2-\mu=0 \\[2pt]
       \mu_2\gamma+\mu_3\alpha=0\\[2pt]
       \mu_3\gamma=0\\[2pt]
       \mu=0\\[2pt]
       \mu_1\beta+\alpha\gamma=0\\[2pt]
       \mu_1\gamma-\alpha\beta-\gamma^2+\mu=0\\[2pt]
       \end{cases}
\end{eqnarray}
By solving (5.13), we get $\beta=\gamma=0$, there is a contradiction. So\\
\begin{thm}
$(G_2, g, V)$ is not an affine Ricci soliton associated to the Bott connection $\nabla^{B_1}$.
\end{thm}

\vskip 0.5 true cm
\noindent{\bf 5.3 Affine Ricci solitons of $G_3$}\\
\vskip 0.5 true cm

\begin{lem} The Bott connection $\nabla^{B_1}$ of $G_3$ is given by
\begin{align}
&\nabla^{B_1}_{\widehat{e}_1}\widehat{e}_1=\nabla^{B_1}_{\widehat{e}_1}\widehat{e}_2=\nabla^{B_1}_{\widehat{e}_1}\widehat{e}_3=-\gamma\widehat{e}_3,\nonumber\\
&\nabla^{B_1}_{\widehat{e}_2}\widehat{e}_1=\gamma\widehat{e}_3,~~~\nabla^{B_1}_{\widehat{e}_2}\widehat{e}_2=0,~~~\nabla^{B_1}_{\widehat{e}_2}\widehat{e}_3=\alpha\widehat{e}_1,\nonumber\\
&\nabla^{B_1}_{\widehat{e}_3}\widehat{e}_1=\nabla^{B_1}_{\widehat{e}_3}\widehat{e}_2=\nabla^{B_1}_{\widehat{e}_3}\widehat{e}_3=0.
\end{align}
\end{lem}
\begin{lem} The curvature $R^{B_1}$ of the Bott connection $\nabla^{B_1}$ of $(G_3,g)$ is given by
\begin{align}
&R^{B_1}(\widehat{e}_1,\widehat{e}_2)\widehat{e}_1=R^{B_1}(\widehat{e}_1,\widehat{e}_2)\widehat{e}_2=R^{B_1}(\widehat{e}_1,\widehat{e}_2)\widehat{e}_3=0,\nonumber\\
&R^{B_1}(\widehat{e}_1,\widehat{e}_3)\widehat{e}_1=\beta\gamma\widehat{e}_3,~~~R^{B_1}(\widehat{e}_1,\widehat{e}_3)\widehat{e}_2=0,~~~R^{B_1}(\widehat{e}_1,\widehat{e}_3)\widehat{e}_3=\alpha\beta\widehat{e}_1,\nonumber\\
&R^{B_1}(\widehat{e}_2,\widehat{e}_3)\widehat{e}_1=R^{B_1}(\widehat{e}_2,\widehat{e}_3)\widehat{e}_2=R^{B_1}(\widehat{e}_2,\widehat{e}_3)\widehat{e}_3=0.
\end{align}
\end{lem}
By (2.5), we have
\begin{align}
&\rho^{B_1}(\widehat{e}_1,\widehat{e}_1)=-\beta\gamma,~~~\rho^{B_1}(\widehat{e}_1,\widehat{e}_2)=0,~~~\rho^{B_1}(\widehat{e}_1,\widehat{e}_3)=0,\nonumber\\
&\rho^{B_1}(\widehat{e}_2,\widehat{e}_1)=\rho^{B_1}(\widehat{e}_2,\widehat{e}_2)=\rho^{B_1}(\widehat{e}_2,\widehat{e}_3)=0,\nonumber\\
&\rho^{B_1}(\widehat{e}_3,\widehat{e}_1)=\rho^{B_1}(\widehat{e}_3,\widehat{e}_2)=0,~~~\rho^{B_1}(\widehat{e}_3,\widehat{e}_3)=\alpha\beta.
\end{align}
Then,
\begin{align}
&\widetilde{\rho}^{B_1}(\widehat{e}_1,\widehat{e}_1)=-\beta\gamma,~~~\widetilde{\rho}^{B_1}(\widehat{e}_1,\widehat{e}_2)=\widetilde{\rho}^{B_1}(\widehat{e}_1,\widehat{e}_3)=0,\nonumber\\
&\widetilde{\rho}^{B_1}(\widehat{e}_2,\widehat{e}_2)=\widetilde{\rho}^{B_1}(\widehat{e}_2,\widehat{e}_3)=0,~~~\widetilde{\rho}^{B_1}(\widehat{e}_3,\widehat{e}_3)=\alpha\beta.
\end{align}
By (2.7), we have
\begin{align}
&(L^{B_1}_Vg)(\widehat{e}_1,\widehat{e}_1)=0,~~~(L^{B_1}_Vg)(\widehat{e}_1,\widehat{e}_2)=\mu_3\alpha,~~~(L^{B_1}_Vg)(\widehat{e}_1,\widehat{e}_3)=0,\nonumber\\
&(L^{B_1}_Vg)(\widehat{e}_2,\widehat{e}_2)=0,~~~(L^{B_1}_Vg)(\widehat{e}_2,\widehat{e}_3)=-\mu_1\gamma,~~~(L^{B_1}_Vg)(\widehat{e}_3,\widehat{e}_3)=0.
\end{align}
Then, if $(G_3,g,V)$  is an affine Ricci soliton associated to the Bott connection $\nabla^{B_1}$, by (2.8), we have the following five equations:
\begin{eqnarray}
       \begin{cases}
        \mu-\beta\gamma=0 \\[2pt]
       \mu_3\alpha=0\\[2pt]
       \mu=0\\[2pt]
       \mu_1\gamma=0\\[2pt]
       \alpha\beta-\mu=0\\[2pt]
       \end{cases}
\end{eqnarray}
By solving (5.19), we get \\
\begin{thm}
$(G_3, g, V)$ is an affine Ricci soliton associated to the Bott connection $\nabla^{B_1}$ if and only if\\
\begin{eqnarray*}
&&(1)\mu=\beta=\mu_3\alpha=\mu_1\gamma=0;\nonumber\\
&&(2)\mu=\alpha=\gamma=0,~~~\beta\neq0.
\end{eqnarray*}
\end{thm}

\vskip 0.5 true cm
\noindent{\bf 5.4 Affine Ricci solitons of $G_4$}\\
\vskip 0.5 true cm

\begin{lem} The Bott connection $\nabla^{B_1}$ of $G_4$ is given by
\begin{align}
&\nabla^{B_1}_{\widehat{e}_1}\widehat{e}_1=0,~~~\nabla^{B_1}_{\widehat{e}_1}\widehat{e}_2=-\widehat{e}_2,~~~\nabla^{B_1}_{\widehat{e}_1}\widehat{e}_3=0,\nonumber\\
&\nabla^{B_1}_{\widehat{e}_2}\widehat{e}_1=(\beta-2\eta)\widehat{e}_3,~~~\nabla^{B_1}_{\widehat{e}_2}\widehat{e}_2=0,~~~\nabla^{B_1}_{\widehat{e}_2}\widehat{e}_3=\alpha\widehat{e}_1,\nonumber\\
&\nabla^{B_1}_{\widehat{e}_3}\widehat{e}_1=-\widehat{e}_3,~~~\nabla^{B_1}_{\widehat{e}_3}\widehat{e}_2=0,~~~\nabla^{B_1}_{\widehat{e}_3}\widehat{e}_3=-\widehat{e}_1.
\end{align}
\end{lem}
\begin{lem} The curvature $R^{B_1}$ of the Bott connection $\nabla^{B_1}$ of $(G_4,g)$ is given by
\begin{align}
&R^{B_1}(\widehat{e}_1,\widehat{e}_2)\widehat{e}_1=0,~~~R^{B_1}(\widehat{e}_1,\widehat{e}_2)\widehat{e}_2=0,~~~R^{B_1}(\widehat{e}_1,\widehat{e}_2)\widehat{e}_3=(2\eta+\alpha-\beta)\widehat{e}_1,\nonumber\\
&R^{B_1}(\widehat{e}_1,\widehat{e}_3)\widehat{e}_1=(\beta-\eta)^2\widehat{e}_3,~~~R^{B_1}(\widehat{e}_1,\widehat{e}_3)\widehat{e}_2=0,~~~R^{B_1}(\widehat{e}_1,\widehat{e}_3)\widehat{e}_3=(1+\alpha\beta)\widehat{e}_1,\nonumber\\
&R^{B_1}(\widehat{e}_2,\widehat{e}_3)\widehat{e}_1=(\beta-2\eta-\alpha)\widehat{e}_1,~~~R^{B_1}(\widehat{e}_2,\widehat{e}_3)\widehat{e}_2=\alpha\widehat{e}_2,~~~R^{B_1}(\widehat{e}_2,\widehat{e}_3)\widehat{e}_3=(2\eta+\alpha-\beta)\widehat{e}_3.
\end{align}
\end{lem}
By (2.5), we have
\begin{align}
&\rho^{B_1}(\widehat{e}_1,\widehat{e}_1)=-(\beta-\eta)^2,~~~\rho^{B_1}(\widehat{e}_1,\widehat{e}_2)=\rho^{B_1}(\widehat{e}_1,\widehat{e}_3)=0,\nonumber\\
&\rho^{B_1}(\widehat{e}_2,\widehat{e}_1)=\rho^{B_1}(\widehat{e}_2,\widehat{e}_2)=\rho^{B_1}(\widehat{e}_2,\widehat{e}_3)=0,\nonumber\\
&\rho^{B_1}(\widehat{e}_3,\widehat{e}_1)=0,~~~\rho^{B_1}(\widehat{e}_3,\widehat{e}_2)=\alpha,~~~\rho^{B_1}(\widehat{e}_3,\widehat{e}_3)=\alpha\beta+1.
\end{align}
Then,
\begin{align}
&\widetilde{\rho}^{B_1}(\widehat{e}_1,\widehat{e}_1)=-(\beta-\eta)^2,~~~\widetilde{\rho}^{B_1}(\widehat{e}_1,\widehat{e}_2)=0,~~~\widetilde{\rho}^{B_1}(\widehat{e}_1,\widehat{e}_3)=0,\nonumber\\
&\widetilde{\rho}^{B_1}(\widehat{e}_2,\widehat{e}_2)=0,~~~\widetilde{\rho}^{B_1}(\widehat{e}_2,\widehat{e}_3)=\frac{\alpha}{2},~~~\widetilde{\rho}^{B_1}(\widehat{e}_3,\widehat{e}_3)=\alpha\beta+1.
\end{align}
By (2.7), we have
\begin{align}
&(L^{B_1}_Vg)(\widehat{e}_1,\widehat{e}_1)=0,~~~(L^{B_1}_Vg)(\widehat{e}_1,\widehat{e}_2)=-\mu_2+\mu_3\alpha,~~~(L^{B_1}_Vg)(\widehat{e}_1,\widehat{e}_3)=-\mu_3,\nonumber\\
&(L^{B_1}_Vg)(\widehat{e}_2,\widehat{e}_2)=0,~~~(L^{B_1}_Vg)(\widehat{e}_2,\widehat{e}_3)=\mu_1(2\eta-\beta),~~~(L^{B_1}_Vg)(\widehat{e}_3,\widehat{e}_3)=2\mu_1.
\end{align}
Then, if $(G_4,g,V)$  is an affine Ricci soliton associated to the Bott connection $\nabla^{B_1}$, by (2.8), we have the following six equations:
\begin{eqnarray}
       \begin{cases}
        (\beta-\eta)^2-\mu=0 \\[2pt]
       \mu_2-\mu_3\alpha=0\\[2pt]
       \mu_3=0\\[2pt]
        \mu=0\\[2pt]
       \mu_1(2\eta-\beta)+\alpha=0\\[2pt]
              \mu_1+\alpha\beta+1-\mu=0\\[2pt]
       \end{cases}
\end{eqnarray}
By solving (5.25), we get \\
\begin{thm}
$(G_4, g, V)$ is not an affine Ricci soliton associated to the Bott connection $\nabla^{B_1}$.
\end{thm}

\vskip 0.5 true cm
\noindent{\bf 5.5 Affine Ricci solitons of $G_5$}\\
\vskip 0.5 true cm

\begin{lem} The Bott connection $\nabla^{B_1}$ of $G_5$ is given by
\begin{align}
&\nabla^{B_1}_{\widehat{e}_1}\widehat{e}_1=\alpha\widehat{e}_3,~~~\nabla^{B_1}_{\widehat{e}_1}\widehat{e}_2=0,~~~\nabla^{B_1}_{\widehat{e}_1}\widehat{e}_3=\alpha\widehat{e}_1,\nonumber\\
&\nabla^{B_1}_{\widehat{e}_2}\widehat{e}_1=0,~~~\nabla^{B_1}_{\widehat{e}_2}\widehat{e}_2=0,~~~\nabla^{B_1}_{\widehat{e}_2}\widehat{e}_3=\gamma\widehat{e}_1,\nonumber\\
&\nabla^{B_1}_{\widehat{e}_3}\widehat{e}_1=0,~~~\nabla^{B_1}_{\widehat{e}_3}\widehat{e}_2=-\delta\widehat{e}_2,~~~\nabla^{B_1}_{\widehat{e}_3}\widehat{e}_3=0.
\end{align}
\end{lem}
\begin{lem} The curvature $R^{B_1}$ of the Bott connection $\nabla^{B_1}$ of $(G_5,g)$ is given by
\begin{align}
&R^{B_1}(\widehat{e}_1,\widehat{e}_2)\widehat{e}_1=-\alpha\gamma\widehat{e}_1,~~~R^{B_1}(\widehat{e}_1,\widehat{e}_2)\widehat{e}_2=0,~~~R^{B_1}(\widehat{e}_1,\widehat{e}_2)\widehat{e}_3=\alpha\gamma\widehat{e}_3,\nonumber\\
&R^{B_1}(\widehat{e}_1,\widehat{e}_3)\widehat{e}_1=-\alpha^2\widehat{e}_3,~~~R^{B_1}(\widehat{e}_1,\widehat{e}_3)\widehat{e}_2=0,~~~R^{B_1}(\widehat{e}_1,\widehat{e}_3)\widehat{e}_3=-(\beta\gamma+\alpha^2)\widehat{e}_2,\nonumber\\
&R^{B_1}(\widehat{e}_2,\widehat{e}_3)\widehat{e}_1=-\alpha\gamma\widehat{e}_3,~~~R^{B_1}(\widehat{e}_2,\widehat{e}_3)\widehat{e}_2=0,~~~R^{B_1}(\widehat{e}_2,\widehat{e}_3)\widehat{e}_3=-\gamma(\alpha+\delta)\widehat{e}_1.
\end{align}
\end{lem}
By (2.5), we have
\begin{align}
&\rho^{B_1}(\widehat{e}_1,\widehat{e}_1)=\alpha^2,~~~\rho^{B_1}(\widehat{e}_1,\widehat{e}_2)=\rho^{B_1}(\widehat{e}_1,\widehat{e}_3)=0,\nonumber\\
&\rho^{B_1}(\widehat{e}_2,\widehat{e}_1)=\rho^{B_1}(\widehat{e}_2,\widehat{e}_2)=\rho^{B_1}(\widehat{e}_2,\widehat{e}_3)=0,\nonumber\\
&\rho^{B_1}(\widehat{e}_3,\widehat{e}_1)=\rho^{B_1}(\widehat{e}_3,\widehat{e}_2)=0,~~~\rho^{B_1}(\widehat{e}_3,\widehat{e}_3)=-(\beta\gamma+\alpha^2).
\end{align}
Then,
\begin{align}
&\widetilde{\rho}^{B_1}(\widehat{e}_1,\widehat{e}_1)=\alpha^2,~~~\widetilde{\rho}^{B_1}(\widehat{e}_1,\widehat{e}_2)=\widetilde{\rho}^{B_1}(\widehat{e}_1,\widehat{e}_3)=0,\nonumber\\
&\widetilde{\rho}^{B_1}(\widehat{e}_2,\widehat{e}_2)=\widetilde{\rho}^{B_1}(\widehat{e}_2,\widehat{e}_3)=0,~~~\widetilde{\rho}^{B_1}(\widehat{e}_3,\widehat{e}_3)=-(\beta\gamma+\alpha^2).
\end{align}
By (2.7), we have
\begin{align}
&(L^{B_1}_Vg)(\widehat{e}_1,\widehat{e}_1)=2\mu_3\alpha,~~~(L^{B_1}_Vg)(\widehat{e}_1,\widehat{e}_2)=\mu_3\gamma,~~~(L^{B_1}_Vg)(\widehat{e}_1,\widehat{e}_3)=-\mu_1\alpha,\nonumber\\
&(L^{B_1}_Vg)(\widehat{e}_2,\widehat{e}_2)=0,~~~(L^{B_1}_Vg)(\widehat{e}_2,\widehat{e}_3)=-\mu_2\delta,~~~(L^{B_1}_Vg)(\widehat{e}_3,\widehat{e}_3)=0.
\end{align}
Then, if $(G_5,g,V)$  is an affine Ricci soliton associated to the Bott connection $\nabla^{B_1}$, by (2.8), we have the following six equations:
\begin{eqnarray}
       \begin{cases}
       \mu_3\alpha+\alpha^2+\mu=0\\[2pt]
       \mu_3\gamma=0\\[2pt]
       \mu_1\alpha=0\\[2pt]
       \mu=0\\[2pt]
       \mu_2\delta=0\\[2pt]
       \alpha^2+\beta\gamma+\mu=0\\[2pt]
       \end{cases}
\end{eqnarray}
By solving (5.31), we get \\
\begin{thm}
$(G_5, g, V)$ is an affine Ricci soliton associated to the Bott connection $\nabla^{B_1}$ if and only if
$\mu=\alpha=\beta=\mu_2=\mu_3\gamma=0,~~~ \delta\neq 0.$
\end{thm}

\vskip 0.5 true cm
\noindent{\bf 5.6 Affine Ricci solitons of $G_6$}\\
\vskip 0.5 true cm
\begin{lem} The Bott connection $\nabla^{B_1}$ of $G_6$ is given by
\begin{align}
&\nabla^{B_1}_{\widehat{e}_1}\widehat{e}_1=0,~~~\nabla^{B_1}_{\widehat{e}_1}\widehat{e}_2=\alpha\widehat{e}_2,~~~\nabla^{B_1}_{\widehat{e}_1}\widehat{e}_3=0,\nonumber\\
&\nabla^{B_1}_{\widehat{e}_2}\widehat{e}_1=-\beta\widehat{e}_3,~~~\nabla^{B_1}_{\widehat{e}_2}\widehat{e}_2=\nabla^{B_1}_{\widehat{e}_2}\widehat{e}_3=0,\nonumber\\
&\nabla^{B_1}_{\widehat{e}_3}\widehat{e}_1=-\delta\widehat{e}_3,~~~\nabla^{B_1}_{\widehat{e}_3}\widehat{e}_2=0,~~~\nabla^{B_1}_{\widehat{e}_3}\widehat{e}_3=-\delta\widehat{e}_1.
\end{align}
\end{lem}
\begin{lem} The curvature $R^{B_1}$ of the Bott connection $\nabla^{B_1}$ of $(G_6,g)$ is given by
\begin{align}
&R^{B_1}(\widehat{e}_1,\widehat{e}_2)\widehat{e}_1=\beta(\alpha+\delta)\widehat{e}_3,~~~R^{B_1}(\widehat{e}_1,\widehat{e}_2)\widehat{e}_2=0,~~~R^{B_1}(\widehat{e}_1,\widehat{e}_2)\widehat{e}_3=\beta\delta\widehat{e}_1,\nonumber\\
&R^{B_1}(\widehat{e}_1,\widehat{e}_3)\widehat{e}_1=(\beta\gamma+\delta^2)\widehat{e}_3,~~~R^{B_1}(\widehat{e}_1,\widehat{e}_3)\widehat{e}_2=0,~~~R^{B_1}(\widehat{e}_1,\widehat{e}_3)\widehat{e}_3=\delta^2\widehat{e}_1,\nonumber\\
&R^{B_1}(\widehat{e}_2,\widehat{e}_3)\widehat{e}_1=-\beta\delta\widehat{e}_1,~~~R^{B_1}(\widehat{e}_2,\widehat{e}_3)\widehat{e}_2=0,~~~R^{B_1}(\widehat{e}_2,\widehat{e}_3)\widehat{e}_3=\beta\delta\widehat{e}_3.
\end{align}
\end{lem}
By (2.5), we have
\begin{align}
&\rho^{B_1}(\widehat{e}_1,\widehat{e}_1)=-(\beta\gamma+\delta^2),~~~\rho^{B_1}(\widehat{e}_1,\widehat{e}_2)=\rho^{B_1}(\widehat{e}_1,\widehat{e}_3)=0,\nonumber\\
&\rho^{B_1}(\widehat{e}_2,\widehat{e}_1)=\rho^{B_1}(\widehat{e}_2,\widehat{e}_2)=\rho^{B_1}(\widehat{e}_2,\widehat{e}_3)=0,\nonumber\\
&\rho^{B_1}(\widehat{e}_3,\widehat{e}_1)=\rho^{B_1}(\widehat{e}_3,\widehat{e}_2)=0,~~~\rho^{B_1}(\widehat{e}_3,\widehat{e}_3)=\delta^2.
\end{align}
Then,
\begin{align}
&\widetilde{\rho}^{B_1}(\widehat{e}_1,\widehat{e}_1)=-(\delta^2+\beta\gamma),~~~\widetilde{\rho}^{B_1}(\widehat{e}_1,\widehat{e}_2)=\widetilde{\rho}^{B_1}(\widehat{e}_1,\widehat{e}_3)=0,\nonumber\\
&\widetilde{\rho}^{B_1}(\widehat{e}_2,\widehat{e}_2)=\widetilde{\rho}^{B_1}(\widehat{e}_2,\widehat{e}_3)=0,~~~\widetilde{\rho}^{B_1}(\widehat{e}_3,\widehat{e}_3)=\delta^2.
\end{align}
By (2.7), we have
\begin{align}
&(L^{B_1}_Vg)(\widehat{e}_1,\widehat{e}_1)=0,~~~(L^{B_1}_Vg)(\widehat{e}_1,\widehat{e}_2)=\mu_2\alpha,~~~(L^{B_1}_Vg)(\widehat{e}_1,\widehat{e}_3)=-\mu_3\delta,\nonumber\\
&(L^{B_1}_Vg)(\widehat{e}_2,\widehat{e}_2)=0,~~~(L^{B_1}_Vg)(\widehat{e}_2,\widehat{e}_3)=\mu_1\beta,~~~(L^{B_1}_Vg)(\widehat{e}_3,\widehat{e}_3)=2\mu_1\delta.
\end{align}
Then, if $(G_6,g,V)$  is an affine Ricci soliton associated to the Bott connection $\nabla^{B_1}$, by (2.8), we have the following six equations:
\begin{eqnarray}
       \begin{cases}
       \delta^2+\beta\gamma-\mu=0\\[2pt]
       \mu_2\alpha=0\\[2pt]
       \mu_3\delta=0\\[2pt]
       \mu=0\\[2pt]
       \mu_1\beta=0\\[2pt]
       \mu_1\delta+\delta^2-\mu=0\\[2pt]
       \end{cases}
\end{eqnarray}
By solving (5.37), we get \\
\begin{thm}
$(G_6, g, V)$ is an affine Ricci soliton associated to the Bott connection $\nabla^{B_1}$ if and only if
$\mu=\mu_2=\delta=\gamma=\mu_1\beta=0,~~~ \alpha\neq0.$
\end{thm}
\vskip 0.5 true cm
\noindent{\bf 5.7 Affine Ricci solitons of $G_7$}\\
\vskip 0.5 true cm
\begin{lem} The Bott connection $\nabla^{B_1}$ of $G_7$ is given by
\begin{align}
&\nabla^{B_1}_{\widehat{e}_1}\widehat{e}_1=\alpha\widehat{e}_3,~~~\nabla^{B_1}_{\widehat{e}_1}\widehat{e}_2=-\beta\widehat{e}_2,~~~\nabla^{B_1}_{\widehat{e}_1}\widehat{e}_3=\alpha\widehat{e}_1,\nonumber\\
&\nabla^{B_1}_{\widehat{e}_2}\widehat{e}_1=\alpha\widehat{e}_1+\beta\widehat{e}_3,~~~\nabla^{B_1}_{\widehat{e}_2}\widehat{e}_2=0,~~~\nabla^{B_1}_{\widehat{e}_2}\widehat{e}_3=\gamma\widehat{e}_1+\delta\widehat{e}_3,\nonumber\\
&\nabla^{B_1}_{\widehat{e}_3}\widehat{e}_1=-\beta\widehat{e}_3,~~~\nabla^{B_1}_{\widehat{e}_3}\widehat{e}_2=-\delta\widehat{e}_2,~~~\nabla^{B_1}_{\widehat{e}_3}\widehat{e}_3=-\beta\widehat{e}_1.
\end{align}
\end{lem}
\begin{lem} The curvature $R^{B_1}$ of the Bott connection $\nabla^{B_1}$ of $(G_7,g)$ is given by
\begin{align}
&R^{B_1}(\widehat{e}_1,\widehat{e}_2)\widehat{e}_1=\alpha(2\beta-\gamma)\widehat{e}_1+\alpha(2\alpha-\delta)\widehat{e}_3,~~~R^{B_1}(\widehat{e}_1,\widehat{e}_2)\widehat{e}_2=-\beta(\alpha+\delta)\widehat{e}_2,\nonumber\\
&R^{B_1}(\widehat{e}_1,\widehat{e}_2)\widehat{e}_3=(\alpha\delta+\beta\gamma-\beta^2)\widehat{e}_1+(\alpha\gamma+\beta\delta-\alpha\beta)\widehat{e}_3,~~~R^{B_1}(\widehat{e}_1,\widehat{e}_3)\widehat{e}_1=-\alpha\beta\widehat{e}_1-\alpha^2\widehat{e}_3,\nonumber\\
&R^{B_1}(\widehat{e}_1,\widehat{e}_3)\widehat{e}_2=\beta(\alpha+\delta)\widehat{e}_2,~~~R^{B_1}(\widehat{e}_1,\widehat{e}_3)\widehat{e}_3=(\beta^2-\alpha^2-\beta\gamma)\widehat{e}_1-\beta\delta\widehat{e}_3,\nonumber\\
&R^{B_1}(\widehat{e}_2,\widehat{e}_3)\widehat{e}_1=(\beta^2-\beta\gamma-\alpha\delta)\widehat{e}_1+(\alpha\beta-\beta\delta-\alpha\gamma)\widehat{e}_3,~~~R^{B_1}(\widehat{e}_2,\widehat{e}_3)\widehat{e}_2=(\beta\gamma+\delta^2)\widehat{e}_2,\nonumber\\
&R^{B_1}(\widehat{e}_2,\widehat{e}_3)\widehat{e}_3=(2\beta\delta-\delta\gamma-\alpha\gamma-\alpha\beta)\widehat{e}_1+(\beta\gamma-\beta^2-\delta^2)\widehat{e}_3.
\end{align}
\end{lem}
By (2.5), we have
\begin{align}
&\rho^{B_1}(\widehat{e}_1,\widehat{e}_1)=\alpha^2,~~~\rho^{B_1}(\widehat{e}_1,\widehat{e}_2)=\beta(\alpha+\delta),~~~\rho^{B_1}(\widehat{e}_1,\widehat{e}_3)=\beta\delta,\nonumber\\
&\rho^{B_1}(\widehat{e}_2,\widehat{e}_1)=\beta(\alpha+\delta),~~~\rho^{B_1}(\widehat{e}_2,\widehat{e}_2)=0,~~~\rho^{B_1}(\widehat{e}_2,\widehat{e}_3)=\delta(\alpha+\delta),\nonumber\\
&\rho^{B_1}(\widehat{e}_3,\widehat{e}_1)=-\alpha\beta,~~~\rho^{B_1}(\widehat{e}_3,\widehat{e}_2)=\beta\gamma+\delta^2,~~~\rho^{B_1}(\widehat{e}_3,\widehat{e}_3)=\beta^2-\alpha^2-\beta\gamma.
\end{align}
Then,
\begin{align}
&\widetilde{\rho}^{B_1}(\widehat{e}_1,\widehat{e}_1)=\alpha^2,~~~\widetilde{\rho}^{B_1}(\widehat{e}_1,\widehat{e}_2)=\beta(\alpha+\delta),~~~\widetilde{\rho}^{B_1}(\widehat{e}_1,\widehat{e}_3)=\frac{\beta(\delta-\alpha)}{2},\nonumber\\
&\widetilde{\rho}^{B_1}(\widehat{e}_2,\widehat{e}_2)=0,~~~\widetilde{\rho}^{B_1}(\widehat{e}_2,\widehat{e}_3)=\delta^2+\frac{\beta\gamma+\alpha\delta}{2},~~~\widetilde{\rho}^{B_1}(\widehat{e}_3,\widehat{e}_3)=\beta^2-\alpha^2-\beta\gamma.
\end{align}
By (2.7), we have
\begin{align}
&(L^{B_1}_Vg)(\widehat{e}_1,\widehat{e}_1)=2\mu_3\alpha,~~~(L^{B_1}_Vg)(\widehat{e}_1,\widehat{e}_2)=\mu_1\alpha-\mu_2\beta+\mu_3\gamma,~~~(L^{B_1}_Vg)(\widehat{e}_1,\widehat{e}_3)=-\mu_1\alpha-\mu_3\beta,\nonumber\\
&(L^{B_1}_Vg)(\widehat{e}_2,\widehat{e}_2)=0,~~~(L^{B_1}_Vg)(\widehat{e}_2,\widehat{e}_3)=-\mu_1\beta-\mu_2\delta-\mu_3\delta,~~~(L^{B_1}_Vg)(\widehat{e}_3,\widehat{e}_3)=2\mu_1\beta.
\end{align}
Then, if $(G_7,g,V)$  is an affine Ricci soliton associated to the Bott connection $\nabla^{B_1}$, by (2.8), we have the following six equations:
\begin{eqnarray}
       \begin{cases}
       \mu_3\alpha+\alpha^2+\mu=0\\[2pt]
       \mu_1\alpha-\mu_2\beta+\mu_3\gamma+2\beta(\delta+\alpha)=0\\[2pt]
       \mu_3\beta+\mu_1\alpha+\beta(\alpha-\delta)=0\\[2pt]
       \mu=0\\[2pt]
       \mu_3\delta+\mu_1\beta+\mu_2\delta-2\delta^2-\beta\gamma-\alpha\delta=0\\[2pt]
       \mu_1\beta+\beta^2-\alpha^2-\beta\gamma-\mu=0\\[2pt]
       \end{cases}
\end{eqnarray}
By solving (5.43), we get \\
\begin{thm}
$(G_7, g, V)$ is an affine Ricci soliton associated to the Bott connection $\nabla^{B_1}$ if and only if\\
\begin{eqnarray*}
&&(1)\mu=\alpha=\beta=\gamma=0,~~~ \delta\neq 0, ~~~\mu_2+\mu_3-2\delta=0;\nonumber\\
&&(2)\mu=\mu_3=\alpha=\beta=0,~~~\gamma\neq0,~~~\delta\neq0,~~~\mu_2-2\delta=0;\nonumber\\
&&(3)\mu=\alpha=0,~~~\beta\neq0,~~~\delta\neq0,~~~\mu_1=\gamma-\beta,~~~\mu_3=\delta,~~~\mu_2=\frac{\gamma\delta+2\beta\delta}{\beta},~~~\gamma=\frac{\beta(\beta^2-\beta\delta)}{\delta};\nonumber\\
\end{eqnarray*}
\end{thm}

\indent Specially, let $V=0$, we get the following corollary:
\begin{cor}
(I)~~$(G_1, g, V)$ is not an affine Einstein associated to the Bott connection $\nabla^{B_1}$;\nonumber\\
(II)~~$(G_2, g, V)$ is not an affine Einstein associated to the Bott connection $\nabla^{B_1}$;\nonumber\\
(III)~~$(G_3, g, V)$ is an affine Einstein associated to the Bott connection $\nabla^{B_1}$ if and only if
$(1)\mu=\alpha=\gamma=0,~~~\beta\neq0;$
$(2)\mu=\beta=0;$\nonumber\\
(IV)~~$(G_4, g, V)$ is not an affine Einstein associated to the Bott connection $\nabla^{B_1}$;\nonumber\\
(V)~~$(G_5, g, V)$ is an affine Einstein associated to the Bott connection $\nabla^{B_1}$ if and only if $\mu=\alpha=\beta=0,~~~\delta\neq0$;\nonumber\\
(VI)~~$(G_6, g, V)$ is an affine Einstein associated to the Bott connection $\nabla^{B_1}$ if and only if $\mu=\delta=\gamma=0,~~~\alpha\neq0;$\nonumber\\
(VII)~~$(G_7, g, V)$ is not an affine Einstein associated to the Bott connection $\nabla^{B_1}$.\nonumber\\
\end{cor}
\section{Affine Ricci Solitons associated to the perturbed Bott connection on three-dimensional Lorentzian Lie groups with the second distribution}
Similarly, by the above calculations, we always obtain $\mu=0$. In order to get the affine Ricci soliton with non zero $\mu$, we introduce the perturbed Bott connection $\widetilde{\nabla}^{B_1}$ in the following. Let $\widehat{e}_2^*$ be the dual base of $e_2$. We define on $G_{i=1,\cdot\cdot\cdot,7}$
\begin{align}
\widetilde{\nabla}^{B_1}_XY=\nabla^{B_1}_XY+a_0\widehat{e}_2^*(X)\widehat{e}_2^*(Y)e_2,
\end{align}
where $a_0$ is a non zero real number. Then
\begin{align}
\widetilde{\nabla}^{B_1}_{\widehat{e}_2}\widehat{e}_2=a_0\widehat{e}_2,~~~\widetilde{\nabla}^{B_1}_{\widehat{e}_s}\widehat{e}_t=\nabla^{B_1}_{\widehat{e}_s}\widehat{e}_t,
\end{align}
where $s$ and $t$ does not equal 2. We define
\begin{align}
(\widetilde{L}_V^{B_1}g)(X,Y):=g(\widetilde{\nabla}^{B_1}_XV,Y)+g(X,\widetilde{\nabla}^{B_1}_YV),
\end{align}
for vector fields $X,Y,V$. Then we have for $G_{i=1,\cdot\cdot\cdot,7}$
\begin{align}
(\widetilde{L}_V^{B_1}g)(\widehat{e}_2,\widehat{e}_2)=2a_0\mu_2,~~~(\widetilde{L}_V^{B_1}g)(\widehat{e}_s,\widehat{e}_t)=(L_V^{B_1}g)(\widehat{e}_s,\widehat{e}_t),
\end{align}
where $s$ and $t$ does not equal 2.
\begin{defn}$(G_i,V,g)$ is called the affine Ricci soliton associated to the connection $\widetilde{\nabla}^{B_1}$ if it satisfies
\begin{align}
(\widetilde{L}_V^{B_1}g)(X,Y)+2\widetilde{\overline{\rho}}^{B_1}(X,Y)+2\mu g(X,Y)=0.
\end{align}
\end{defn}
For $(G_1,\widetilde{\nabla}^{B_1})$, we have
\begin{align}
\widetilde{R}^{B_1}(\widehat{e}_1,\widehat{e}_3)\widehat{e}_2=a_0\beta \widehat{e}_2,~~~\widetilde{R}^{B_1}(\widehat{e}_2,\widehat{e}_3)\widehat{e}_2=(\alpha^2-a_0\alpha)\widehat{e}_2,~~~\widetilde{R}^{B_1}(\widehat{e}_s,\widehat{e}_t)\widehat{e}_p=R^{B_1}(\widehat{e}_s,\widehat{e}_t)\widehat{e}_p,
\end{align}
for $(s,t,p)\neq(1,3,2),(2,3,2)$.\\
\begin{align}
\widetilde{\overline{\rho}}^{B_1}(\widehat{e}_2,\widehat{e}_3)=\frac{\alpha(\alpha-a_0)}{2},~~~\widetilde{\overline{\rho}}^{B_1}(\widehat{e}_s,\widehat{e}_t)=\widetilde{\rho}^{B_1}(\widehat{e}_s,\widehat{e}_t),
\end{align}
for the pair $(s,t)\neq(2,3)$. If $(G_1,g,V)$ is an affine Ricci soliton associated to the connection $\widetilde{\nabla}^{B_1}$, then by (6.5), we have
\begin{eqnarray}
       \begin{cases}
        \mu_3\alpha-\alpha^2+\beta^2-\mu=0 \\[2pt]
       \alpha\beta-\mu_1\alpha+\mu_3\beta=0\\[2pt]
       \mu_1\alpha-2\alpha\beta=0\\[2pt]
       a_0\mu_2+\mu=0\\[2pt]
       (\mu_2+\mu_3)\alpha+\mu_1\beta-\alpha^2+a_0\alpha=0\\[2pt]
      \beta^2-\alpha^2-\mu=0\\[2pt]
       \end{cases}
\end{eqnarray}
Solve (6.8), we get $a_0=\alpha=0$, there is a contradiction. So\\
\begin{thm}
$(G_1,V,g)$ is not an affine Ricci soliton associated to the connection $\widetilde{\nabla}^{B_1}$.
\end{thm}
For $(G_2,\widetilde{\nabla}^{B_1})$, we have
\begin{align}
\widetilde{R}^{B_1}(\widehat{e}_1,\widehat{e}_2)\widehat{e}_2=-a_0\gamma \widehat{e}_2,~~~\widetilde{R}^{B_1}(\widehat{e}_1,\widehat{e}_3)\widehat{e}_2=a_0\beta\widehat{e}_2,~~~\widetilde{R}^{B_1}(\widehat{e}_s,\widehat{e}_t)\widehat{e}_p=R^{B_1}(\widehat{e}_s,\widehat{e}_t)\widehat{e}_p,
\end{align}
for $(s,t,p)\neq(1,2,2),(1,3,2)$.\\
\begin{align}
\widetilde{\overline{\rho}}^{B_1}(\widehat{e}_1,\widehat{e}_2)=\frac{a_0\gamma}{2},~~~\widetilde{\overline{\rho}}^{B_1}(\widehat{e}_s,\widehat{e}_t)=\widetilde{\rho}^{B_1}(\widehat{e}_s,\widehat{e}_t),
\end{align}
for the pair $(s,t)\neq(1,2)$. If $(G_2,g,V)$ is an affine Ricci soliton associated to the connection $\widetilde{\nabla}^{B_1}$, then by (6.5), we have
\begin{eqnarray}
       \begin{cases}
       \gamma^2+\beta^2-\mu=0 \\[2pt]
       \mu_2\gamma+\mu_3\alpha+a_0\gamma=0\\[2pt]
       \mu_3\gamma=0\\[2pt]
        a_0\mu_2+\mu=0\\[2pt]
       \mu_1\beta+\alpha\gamma=0\\[2pt]
       \mu_1\gamma-\alpha\beta-\gamma^2+\mu=0\\[2pt]
       \end{cases}
\end{eqnarray}
Solve (6.11), we get
\begin{thm}
$(G_2,V,g)$ is an affine Ricci soliton associated to the connection $\widetilde{\nabla}^{B_1}$ if and only if
\begin{eqnarray*}
&&(1)\beta=\mu_3=0,~~~\gamma\neq0,~~~a_0=\pm\gamma,~~~\mu_2=\mp\gamma;\nonumber\\
&&(2)\mu_3=0,~~~\gamma\neq0,~~~\beta\neq0,~~~\mu_1=-\frac{\alpha\gamma}{\beta},~~~\beta^3-\alpha\gamma^2-\alpha\beta^2=0.
\end{eqnarray*}
\end{thm}
For $(G_3,\widetilde{\nabla}^{B_1})$, we have
\begin{align}
\widetilde{R}^{B_1}(\widehat{e}_1,\widehat{e}_3)\widehat{e}_2=a_0\beta \widehat{e}_2,~~~\widetilde{R}^{B_1}(\widehat{e}_s,\widehat{e}_t)\widehat{e}_p=R^{B_1}(\widehat{e}_s,\widehat{e}_t)\widehat{e}_p,
\end{align}
for $(s,t,p)\neq(1,3,2)$.\\
\begin{align}
\widetilde{\overline{\rho}}^{B_1}(\widehat{e}_s,\widehat{e}_t)=\widetilde{\rho}^{B_1}(\widehat{e}_s,\widehat{e}_t),
\end{align}
for any pairs $(s,t)$. If $(G_3,g,V)$ is an affine Ricci soliton associated to the connection $\widetilde{\nabla}^{B_1}$, then by (6.5), we have
\begin{eqnarray}
       \begin{cases}
       \beta\gamma-\mu=0 \\[2pt]
      \mu_3\alpha=0\\[2pt]
      a_0\mu_2+\mu=0\\[2pt]
       \mu_1\gamma=0\\[2pt]
       \alpha\beta-\mu=0\\[2pt]
       \end{cases}
\end{eqnarray}
Solve (6.14), we get
\begin{thm}
$(G_3,V,g)$ is an affine Ricci soliton associated to the connection $\widetilde{\nabla}^{B_1}$ if and only if
\begin{eqnarray*}
&&(1)\beta=\mu=\mu_2=\mu_3\alpha=\mu_1\gamma=0;\nonumber\\
&&(2)\beta\neq 0,~~~\gamma=\alpha=\mu=\mu_2=0;\nonumber\\
&&(3)\beta\neq0,~~~\gamma=\alpha\neq0,~~~\mu_1=\mu_3=0,~~~\mu_2=-\frac{\alpha\beta}{a_0},~~~\mu=\alpha\beta.\nonumber\\
\end{eqnarray*}
\end{thm}
For $(G_4,\widetilde{\nabla}^{B_1})$, we have
\begin{align}
\widetilde{R}^{B_1}(\widehat{e}_1,\widehat{e}_2)\widehat{e}_2=a_0 \widehat{e}_2,~~~\widetilde{R}^{B_1}(\widehat{e}_1,\widehat{e}_3)\widehat{e}_2=a_0\beta\widehat{e}_2,~~~\widetilde{R}^{B_1}(\widehat{e}_s,\widehat{e}_t)\widehat{e}_p=R^B(\widehat{e}_s,\widehat{e}_t)\widehat{e}_p,
\end{align}
for $(s,t,p)\neq(1,2,2),(1,3,2)$.\\
\begin{align}
\widetilde{\overline{\rho}}^{B_1}(\widehat{e}_1,\widehat{e}_2)=-\frac{a_0}{2},~~~\widetilde{\overline{\rho}}^{B_1}(\widehat{e}_s,\widehat{e}_t)=\widetilde{\rho}^{B_1}(\widehat{e}_s,\widehat{e}_t),
\end{align}
for the pair $(s,t)\neq(1,2)$. If $(G_4,g,V)$ is an affine Ricci soliton associated to the connection $\widetilde{\nabla}^{B_1}$, then by (6.5), we have
\begin{eqnarray}
       \begin{cases}
       (\beta-\eta)^2-\mu=0 \\[2pt]
       \mu_3\alpha-\mu_2-a_0=0\\[2pt]
       \mu_3=0\\[2pt]
        a_0\mu_2+\mu=0\\[2pt]
       \mu_1(2\eta-\beta)+\alpha=0\\[2pt]
       \mu_1+\alpha\beta+1-\mu=0\\[2pt]
       \end{cases}
\end{eqnarray}
Solve (6.17), we get
\begin{thm}
$(G_4,V,g)$ is an affine Ricci soliton associated to the connection $\widetilde{\nabla}^{B_1}$ if and only if
$\beta\neq\eta,~~~a_0=\pm(\beta-\eta),~~~\mu_3=0,~~~\mu_2=\mp(\beta-\eta),~~~\mu_1=1-\frac{1}{(\beta-\eta)^2},~~~\mu=(\beta-\eta)^2$.
\end{thm}

For $(G_5,\widetilde{\nabla}^{B_1})$, we have
\begin{align}
\widetilde{R}^{B_1}(\widehat{e}_1,\widehat{e}_3)\widehat{e}_2=-a_0\beta \widehat{e}_2,~~~\widetilde{R}^{B_1}(\widehat{e}_2,\widehat{e}_3)\widehat{e}_2=-a_0\delta\widehat{e}_2,~~~
\widetilde{R}^{B_1}(\widehat{e}_s,\widehat{e}_t)\widehat{e}_p=R^{B_1}(\widehat{e}_s,\widehat{e}_t)\widehat{e}_p,
\end{align}
for $(s,t,p)\neq(1,3,2),(2,3,2)$.\\
\begin{align}
\widetilde{\overline{\rho}}^{B_1}(\widehat{e}_2,\widehat{e}_3)=-\frac{a_0\delta}{2},~~~
\widetilde{\overline{\rho}}^{B_1}(\widehat{e}_s,\widehat{e}_t)=\widetilde{\rho}^{B_1}(\widehat{e}_s,\widehat{e}_t),
\end{align}
for the pair $(s,t)\neq(2,3)$. If $(G_5,g,V)$ is an affine Ricci soliton associated to the connection $\widetilde{\nabla}^{B_1}$, then by (6.5), we have
\begin{eqnarray}
       \begin{cases}
       \mu_3\alpha-\alpha^2-\mu=0 \\[2pt]
       \mu_1\alpha-\mu_3\beta=0 \\[2pt]
      \mu_1\alpha=0\\[2pt]
       a_0\mu_2+\mu=0\\[2pt]
       \mu_1\beta+\alpha(\mu_2+\mu_3)+a_0\delta=0\\[2pt]
      \alpha^2+\beta\gamma+\mu=0\\[2pt]
       \end{cases}
\end{eqnarray}
Solve (6.20), we get
\begin{thm}
$(G_5,V,g)$ is an affine Ricci soliton associated to the connection $\widetilde{\nabla}^{B_1}$ if and only if
$\mu_3=\beta=\gamma=0,~~~\alpha\neq0,~~~\mu=-\alpha^2,~~~\mu_2=-\frac{a_0\delta}{\alpha},~~~\alpha^3+a_0^2\delta=0$.
\end{thm}
For $(G_6,\widetilde{\nabla}^{B_1})$, we have
\begin{align}
\widetilde{R}^{B_1}(\widehat{e}_1,\widehat{e}_2)\widehat{e}_2=-a_0\alpha \widehat{e}_2,~~~\widetilde{R}^{B_1}(\widehat{e}_1,\widehat{e}_3)\widehat{e}_2=-a_0\gamma\widehat{e}_2,~~~\widetilde{R}^{B_1}(\widehat{e}_s,\widehat{e}_t)\widehat{e}_p=R^{B_1}(\widehat{e}_s,\widehat{e}_t)\widehat{e}_p,
\end{align}
for $(s,t,p)\neq(1,2,2),(1,3,2)$.\\
\begin{align}
\widetilde{\overline{\rho}}^{B_1}(\widehat{e}_1,\widehat{e}_2)=\frac{a_0\alpha}{2},~~~\widetilde{\overline{\rho}}^{B_1}(\widehat{e}_s,\widehat{e}_t)=\widetilde{\rho}^{B_1}(\widehat{e}_s,\widehat{e}_t),
\end{align}
for the pair $(s,t)\neq(1,2)$. If $(G_6,g,V)$ is an affine Ricci soliton associated to the connection $\widetilde{\nabla}^{B_1}$, then by (6.5), we have
\begin{eqnarray}
       \begin{cases}
       \delta^2+\beta\gamma-\mu=0 \\[2pt]
      \mu_2\alpha+a_0\alpha=0\\[2pt]
       \mu_3\delta=0\\[2pt]
        a_0\mu_2+\mu=0\\[2pt]
       \mu_1\beta=0\\[2pt]
       \mu_1\delta+\delta^2-\mu=0\\[2pt]
       \end{cases}
\end{eqnarray}
Solve (6.23), we get
\begin{thm}
$(G_6,V,g)$ is an affine Ricci soliton associated to the connection $\widetilde{\nabla}^{B_1}$ if and only if
\begin{eqnarray*}
&&(1)\beta=\alpha=\mu_3=\mu_1=0,~~~\delta\neq0,~~~\mu_2=-\frac{\delta^2}{a_0};\nonumber\\
&&(2)\alpha\neq0,~~~\delta\neq0,~~~\mu_1=0,~~~a_0=\pm\delta,~~~\mu_2=\mp\delta.
\end{eqnarray*}
\end{thm}
For $(G_7,\widetilde{\nabla}^{B_1}$, we have
\begin{align}
&\widetilde{R}^{B_1}(\widehat{e}_1,\widehat{e}_2)\widehat{e}_2=(a_0-\alpha-\delta)\beta \widehat{e}_2,~~~\widetilde{R}^{B_1}(\widehat{e}_1,\widehat{e}_3)\widehat{e}_2=\beta(\alpha+\delta-a_0)\widehat{e}_2,\nonumber\\
&\widetilde{R}^{B_1}(\widehat{e}_2,\widehat{e}_3)\widehat{e}_2=(\delta^2+\beta\gamma-a_0\delta)\widehat{e}_2,~~~\widetilde{R}^{B_1}(\widehat{e}_s,\widehat{e}_t)\widehat{e}_p=R^{B_1}(\widehat{e}_s,\widehat{e}_t)\widehat{e}_p,
\end{align}
for $(s,t,p)\neq(1,2,2),(1,3,2),(2,3,2)$.\\
\begin{align}
\widetilde{\overline{\rho}}^{B_1}(\widehat{e}_1,\widehat{e}_2)=\frac{2\beta(\alpha\delta)-a_0\beta}{2},~~~\widetilde{\overline{\rho}}^{B_1}(\widehat{e}_2,\widehat{e}_3)=\frac{2\delta^2+\alpha\delta+\beta\gamma-a_0\delta}{2},~~~\widetilde{\overline{\rho}}^{B_1}(\widehat{e}_s,\widehat{e}_t)=\widetilde{\rho}^{B_1}(\widehat{e}_s,\widehat{e}_t),
\end{align}
for the pair $(s,t)\neq(1,2),(2,3)$. If $(G_7,g,V)$ is an affine Ricci soliton associated to the connection $\widetilde{\nabla}^{B_1}$, then by (6.5), we have
\begin{eqnarray}
       \begin{cases}
        \mu_3\alpha+\alpha^2+\mu=0 \\[2pt]
       \mu_2\beta-\mu_1\alpha-\mu_3\gamma-2\beta\delta-2\alpha\beta+a_0\beta=0\\[2pt]
       \mu_1\alpha+\mu_3\beta+\alpha\beta-\beta\delta=0\\[2pt]
        a_0\mu_2+\mu=0\\[2pt]
       \mu_1\beta+\mu_2\delta+\mu_3\delta-2\delta^2-\alpha\delta-\beta\gamma+a_0\delta=0\\[2pt]
      \mu_1\beta+\beta^2-\alpha^2-\beta\gamma-\mu=0\\[2pt]
       \end{cases}
\end{eqnarray}
Solve (6.26), we get
\begin{thm}
$(G_7,V,g)$ is an affine Ricci soliton associated to the connection $\widetilde{\nabla}^{B_1}$ if and only if
\begin{eqnarray*}
&&(1)\mu=\alpha=\beta=\gamma=\mu_2=0,~~~\delta\neq0,~~~\mu_3=2\delta-a_0;\nonumber\\
&&(2)\alpha=\beta=\mu=\mu_2=\mu_3=0,~~~\gamma\neq0,~~~\delta\neq0,~~~a_0=2\delta;\nonumber\\
&&(3)\alpha=\mu=\mu_2=0,~~~\beta\neq0,~~~\delta\neq0,~~~\mu_1=\gamma-\beta,~~~\mu_3=\delta,~~~a_0=\frac{2\beta\delta+\gamma\delta}{\beta},~~~\gamma=\frac{\beta^3-\beta\delta^2}{\delta^2};\nonumber\\
&&(4)\alpha\neq0,~~~\beta=\gamma=\delta=\mu_1=\mu_3=0,~~~\mu=-\alpha^2,~~~\mu_2=\frac{\alpha^2}{a_0};\nonumber\\
&&(5)\alpha\neq0,~~~\beta=\gamma=\mu_1=\mu_3=0,~~~\mu=-\alpha^2,~~~\mu_2=\alpha+2\delta-a_0,~~~a_0^2-a_0\alpha-2a_0\delta+\alpha^2=0;\nonumber\\
&&(6)\alpha\neq0,~~~\beta\neq0,~~~\gamma=0,~~~\alpha^2\neq\beta^2,~~~\mu=\frac{\alpha\beta^2\delta}{\beta^2-\alpha^2}-\alpha^2,~~~\mu=\frac{\alpha\beta\delta}{\beta^2-\alpha^2}-\beta,\nonumber\\
&&\mu_2=\alpha+2\delta+\frac{\alpha^2\delta}{\beta^2-\alpha^2}-a_0,~~~\alpha\beta^2\delta+\alpha^2\delta^2-\beta^2\delta^2+\alpha^2\beta^2-\beta^4=0.
\end{eqnarray*}
\end{thm}
\section{Affine Ricci Solitons associated to the Bott connection on three-dimensional Lorentzian Lie groups with the third distribution}
Let $TM=span\{\widehat{e}_1,\widehat{e}_2,\widehat{e}_3\}$, then took the distribution $D_2=span\{\widehat{e}_2,\widehat{e}_3\}$ and $D_2^\bot=span\{\widehat{e}_1\}$.\\
Similar to (2.3), we have\\
\begin{eqnarray}
\nabla^{B_2}_XY=
       \begin{cases}
        \pi_{D_2}(\nabla^L_XY),~~~&X,Y\in\Gamma^\infty({D_2}) \\[2pt]
       \pi_{D_2}([X,Y]),~~~&X\in\Gamma^\infty({D_2}^\bot),Y\in\Gamma^\infty({D_2})\\[2pt]
       \pi_{{D_2}^\bot}([X,Y]),~~~&X\in\Gamma^\infty({D_2}),Y\in\Gamma^\infty({D_2}^\bot)\\[2pt]
       \pi_{{D_2}^\bot}(\nabla^L_XY),~~~&X,Y\in\Gamma^\infty({D_2}^\bot)\\[2pt]
       \end{cases}
\end{eqnarray}
where $\pi_{D_2}$(resp. $\pi_{D_2}^\bot$) the projection on ${D_2}$ (resp. ${D_2}^\bot$).\\
\vskip 0.5 true cm
\noindent{\bf 7.1 Affine Ricci solitons of $G_1$}\\
\vskip 0.5 true cm
\begin{lem} The Bott connection $\nabla^{B_2}$ of $G_1$ is given by
\begin{align}
&\nabla^{B_2}_{\widehat{e}_1}\widehat{e}_1=0,~~~\nabla^{B_2}_{\widehat{e}_1}\widehat{e}_2=-\beta\widehat{e}_3,~~~\nabla^{B_2}_{\widehat{e}_1}\widehat{e}_3=-\beta\widehat{e}_2,\nonumber\\
&\nabla^{B_2}_{\widehat{e}_2}\widehat{e}_1=-\alpha\widehat{e}_1,~~~\nabla^{B_2}_{\widehat{e}_2}\widehat{e}_2=\alpha\widehat{e}_3,~~~\nabla^{B_2}_{\widehat{e}_2}\widehat{e}_3=\alpha\widehat{e}_2,\nonumber\\
&\nabla^{B_2}_{\widehat{e}_3}\widehat{e}_1=\alpha\widehat{e}_1,~~~\nabla^{B_2}_{\widehat{e}_3}\widehat{e}_2=-\alpha\widehat{e}_3,~~~\nabla^{B_2}_{\widehat{e}_3}\widehat{e}_3=-\alpha\widehat{e}_2.
\end{align}
\end{lem}
\begin{lem} The curvature $R^{B_2}$ of the Bott connection $\nabla^{B_2}$ of $(G_1,g)$ is given by
\begin{align}
&R^{B_2}(\widehat{e}_1,\widehat{e}_2)\widehat{e}_1=\alpha\beta\widehat{e}_1,~~~R^{B_2}(\widehat{e}_1,\widehat{e}_2)\widehat{e}_2=0,~~~R^{B_2}(\widehat{e}_1,\widehat{e}_2)\widehat{e}_3=0,\nonumber\\
&R^{B_2}(\widehat{e}_1,\widehat{e}_3)\widehat{e}_1=-\alpha\beta\widehat{e}_1,~~~R^{B_2}(\widehat{e}_1,\widehat{e}_3)\widehat{e}_2=0,~~~R^{B_2}(\widehat{e}_1,\widehat{e}_3)\widehat{e}_3=0,\nonumber\\
&R^{B_2}(\widehat{e}_2,\widehat{e}_3)\widehat{e}_1=0,~~~R^{B_2}(\widehat{e}_2,\widehat{e}_3)\widehat{e}_2=\beta^2\widehat{e}_3,~~~R^{B_2}(\widehat{e}_2,\widehat{e}_3)\widehat{e}_3=\beta^2\widehat{e}_2.
\end{align}
\end{lem}
By (2.5), we have
\begin{align}
&\rho^{B_2}(\widehat{e}_1,\widehat{e}_1)=\rho^{B_2}(\widehat{e}_1,\widehat{e}_2)=\rho^{B_2}(\widehat{e}_1,\widehat{e}_3)=0,\nonumber\\
&\rho^{B_2}(\widehat{e}_2,\widehat{e}_1)=\alpha\beta,~~~\rho^{B_2}(\widehat{e}_2,\widehat{e}_2)=-\beta^2,~~~\rho^{B_2}(\widehat{e}_2,\widehat{e}_3)=0,\nonumber\\
&\rho^{B_2}(\widehat{e}_3,\widehat{e}_1)=-\alpha\beta,~~~\rho^{B_2}(\widehat{e}_3,\widehat{e}_2)=0,~~~\rho^{B_2}(\widehat{e}_3,\widehat{e}_3)=\beta^2.
\end{align}
Then,
\begin{align}
&\widetilde{\rho}^{B_2}(\widehat{e}_1,\widehat{e}_1)=0,~~~\widetilde{\rho}^{B_2}(\widehat{e}_1,\widehat{e}_2)=\frac{\alpha\beta}{2},~~~\widetilde{\rho}^{B_2}(\widehat{e}_1,\widehat{e}_3)=-\frac{\alpha\beta}{2},\nonumber\\
&\widetilde{\rho}^{B_2}(\widehat{e}_2,\widehat{e}_2)=-\beta^2,~~~\widetilde{\rho}^{B_2}(\widehat{e}_2,\widehat{e}_3)=0,~~~\widetilde{\rho}^{B_2}(\widehat{e}_3,\widehat{e}_3)=\beta^2.
\end{align}
By (2.7), we have
\begin{align}
&(L^{B_2}_Vg)(\widehat{e}_1,\widehat{e}_1)=0,~~~(L^{B_2}_Vg)(\widehat{e}_1,\widehat{e}_2)=-\mu_1\alpha-\mu_3\beta,~~~(L^{B_2}_Vg)(\widehat{e}_1,\widehat{e}_3)=\mu_1\alpha+\mu_2\beta,\nonumber\\
&(L^{B_2}_Vg)(\widehat{e}_2,\widehat{e}_2)=2\mu_3\alpha,~~~(L^{B_2}_Vg)(\widehat{e}_2,\widehat{e}_3)=-(\mu_2+\mu_3)\alpha,~~~(L^{B_2}_Vg)(\widehat{e}_3,\widehat{e}_3)=2\mu_2\alpha.
\end{align}
Then, if $(G_1,g,V)$  is an affine Ricci soliton associated to the Bott connection $\nabla^{B_2}$, by (2.8), we have the following six equations:
\begin{eqnarray}
       \begin{cases}
        \mu=0 \\[2pt]
      \mu_3\beta+\mu_1\alpha-\alpha\beta=0\\[2pt]
       \mu_1\alpha+\mu_2\beta-\alpha\beta=0\\[2pt]
       (\mu_2+\mu_3)\alpha=0\\[2pt]
       \mu_3\alpha-\beta^2+\mu=0\\[2pt]
       \beta^2+\mu_2\alpha-\mu=0\\[2pt]
       \end{cases}
\end{eqnarray}
By solving (7.7) , we get
\begin{thm}
$(G_1, g, V)$ is an affine Ricci soliton associated to the Bott connection $\nabla^{B_2}$ if and only if \\
$\mu=\beta=\mu_2=\mu_3=0,~~~\alpha\neq0$.
\end{thm}

\vskip 0.5 true cm
\noindent{\bf 7.2 Affine Ricci solitons of $G_2$}\\
\vskip 0.5 true cm
\begin{lem} The Bott connection $\nabla^{B_2}$ of $G_2$ is given by
\begin{align}
&\nabla^{B_2}_{\widehat{e}_1}\widehat{e}_1=0,~~~\nabla^{B_2}_{\widehat{e}_1}\widehat{e}_2=\gamma\widehat{e}_2-\beta\widehat{e}_3,~~~\nabla^{B_2}_{\widehat{e}_1}\widehat{e}_3=-\beta\widehat{e}_2-\gamma\widehat{e}_3,\nonumber\\
&\nabla^{B_2}_{\widehat{e}_2}\widehat{e}_1=\nabla^{B_2}_{\widehat{e}_2}\widehat{e}_2=\nabla^{B_2}_{\widehat{e}_2}\widehat{e}_3=0,\nonumber\\
&\nabla^{B_2}_{\widehat{e}_3}\widehat{e}_1=\nabla^{B_2}_{\widehat{e}_3}\widehat{e}_2=\nabla^{B_2}_{\widehat{e}_3}\widehat{e}_3=0.
\end{align}
\end{lem}
\begin{lem} The curvature $R^{B_2}$ of the Bott connection $\nabla^{B_2}$ of $(G_2,g)$ is given by
\begin{align}
&R^{B_2}(\widehat{e}_1,\widehat{e}_2)\widehat{e}_1=R^{B_2}(\widehat{e}_1,\widehat{e}_2)\widehat{e}_2=R^{B_2}(\widehat{e}_1,\widehat{e}_2)\widehat{e}_3=0,\nonumber\\
&R^{B_2}(\widehat{e}_1,\widehat{e}_3)\widehat{e}_1=R^{B_2}(\widehat{e}_1,\widehat{e}_3)\widehat{e}_2=R^{B_2}(\widehat{e}_1,\widehat{e}_3)\widehat{e}_3=0,\nonumber\\
&R^{B_2}(\widehat{e}_2,\widehat{e}_3)\widehat{e}_1=0,~~~R^{B_2}(\widehat{e}_2,\widehat{e}_3)\widehat{e}_2=-\alpha\gamma\widehat{e}_2+\alpha\beta\widehat{e}_3,~~~R^{B_2}(\widehat{e}_2,\widehat{e}_3)\widehat{e}_3=\alpha\beta\widehat{e}_2+\alpha\gamma\widehat{e}_3.
\end{align}
\end{lem}
By (2.5), we have
\begin{align}
&\rho^{B_2}(\widehat{e}_1,\widehat{e}_1)=\rho^{B_2}(\widehat{e}_1,\widehat{e}_2)=\rho^{B_2}(\widehat{e}_1,\widehat{e}_3)=0,\nonumber\\
&\rho^{B_2}(\widehat{e}_2,\widehat{e}_1)=0,~~~\rho^{B_2}(\widehat{e}_2,\widehat{e}_2)=-\alpha\beta,~~~\rho^{B_2}(\widehat{e}_2,\widehat{e}_3)=-\alpha\gamma,\nonumber\\
&\rho^{B_2}(\widehat{e}_3,\widehat{e}_1)=0,~~~\rho^{B_2}(\widehat{e}_3,\widehat{e}_2)=-\alpha\gamma,~~~\rho^{B_2}(\widehat{e}_3,\widehat{e}_3)=\alpha\beta.
\end{align}
Then,
\begin{align}
&\widetilde{\rho}^{B_2}(\widehat{e}_1,\widehat{e}_1)=\widetilde{\rho}^{B_2}(\widehat{e}_1,\widehat{e}_2)=\widetilde{\rho}^{B_2}(\widehat{e}_1,\widehat{e}_3)=0,\nonumber\\
&\widetilde{\rho}^{B_2}(\widehat{e}_2,\widehat{e}_2)=-\alpha\beta,~~~\widetilde{\rho}^{B_2}(\widehat{e}_2,\widehat{e}_3)=-\alpha\gamma,~~~\widetilde{\rho}^{B_2}(\widehat{e}_3,\widehat{e}_3)=\alpha\beta.
\end{align}
By (2.7), we have
\begin{align}
&(L^{B_2}_Vg)(\widehat{e}_1,\widehat{e}_1)=0,~~~(L^{B_2}_Vg)(\widehat{e}_1,\widehat{e}_2)=\mu_2\gamma-\mu_3\beta,~~~(L^{B_2}_Vg)(\widehat{e}_1,\widehat{e}_3)=\mu_2\beta+\mu_3\gamma,\nonumber\\
&(L^{B_2}_Vg)(\widehat{e}_2,\widehat{e}_2)=(L^{B_2}_Vg)(\widehat{e}_2,\widehat{e}_3)=(L^{B_2}_Vg)(\widehat{e}_3,\widehat{e}_3)=0.
\end{align}
Then, if $(G_2,g,V)$  is an affine Ricci soliton associated to the Bott connection $\nabla^{B_2}$, by (2.8), we have the following five equations:
\begin{eqnarray}
       \begin{cases}
        \mu=0\\[2pt]
        \mu_2\gamma-\mu_3\beta=0 \\[2pt]
       \mu_2\beta+\mu_3\gamma=0\\[2pt]
       \alpha\beta-\mu=0\\[2pt]
       \alpha\gamma=0\\[2pt]
       \end{cases}
\end{eqnarray}
By solving (7.13), we get
\begin{thm}
$(G_2, g, V)$ is an affine Ricci soliton associated to the Bott connection $\nabla^{B_2}$ if and only if
$\mu=\alpha=\mu_2=\mu_3=0,~~~\gamma\neq0$.
\end{thm}

\vskip 0.5 true cm
\noindent{\bf 7.3 Affine Ricci solitons of $G_3$}\\
\vskip 0.5 true cm
\begin{lem} The Bott connection $\nabla^{B_2}$ of $G_3$ is given by
\begin{align}
&\nabla^{B_2}_{\widehat{e}_1}\widehat{e}_1=0,~~~\nabla^{B_2}_{\widehat{e}_1}\widehat{e}_2=-\gamma\widehat{e}_3,~~~\nabla^{B_2}_{\widehat{e}_1}\widehat{e}_3=-\beta\widehat{e}_2,\nonumber\\
&\nabla^{B_2}_{\widehat{e}_2}\widehat{e}_1=\nabla^{B_2}_{\widehat{e}_2}\widehat{e}_2=\nabla^{B_2}_{\widehat{e}_2}\widehat{e}_3=0,\nonumber\\
&\nabla^{B_2}_{\widehat{e}_3}\widehat{e}_1=\nabla^{B_2}_{\widehat{e}_3}\widehat{e}_2=\nabla^{B_2}_{\widehat{e}_3}\widehat{e}_3=0.
\end{align}
\end{lem}
\begin{lem} The curvature $R^{B_2}$ of the Bott connection $\nabla^{B_2}$ of $(G_3,g)$ is given by
\begin{align}
&R^{B_2}(\widehat{e}_1,\widehat{e}_2)\widehat{e}_1=R^{B_2}(\widehat{e}_1,\widehat{e}_2)\widehat{e}_2=R^{B_2}(\widehat{e}_1,\widehat{e}_2)\widehat{e}_3=0,\nonumber\\
&R^{B_2}(\widehat{e}_1,\widehat{e}_3)\widehat{e}_1=R^{B_2}(\widehat{e}_1,\widehat{e}_3)\widehat{e}_2=R^{B_2}(\widehat{e}_1,\widehat{e}_3)\widehat{e}_3=0,\nonumber\\
&R^{B_2}(\widehat{e}_2,\widehat{e}_3)\widehat{e}_1=R^{B_2}(\widehat{e}_2,\widehat{e}_3)\widehat{e}_2=0,~~~R^{B_2}(\widehat{e}_2,\widehat{e}_3)\widehat{e}_3=\alpha\beta\widehat{e}_2.
\end{align}
\end{lem}
By (2.5), we have
\begin{align}
&\rho^{B_2}(\widehat{e}_1,\widehat{e}_1)=\rho^{B_2}(\widehat{e}_1,\widehat{e}_2)=\rho^{B_2}(\widehat{e}_1,\widehat{e}_3)=0,\nonumber\\
&\rho^{B_2}(\widehat{e}_2,\widehat{e}_1)=\rho^{B_2}(\widehat{e}_2,\widehat{e}_2)=\rho^{B_2}(\widehat{e}_2,\widehat{e}_3)=0,\nonumber\\
&\rho^{B_2}(\widehat{e}_3,\widehat{e}_1)=\rho^{B_2}(\widehat{e}_3,\widehat{e}_2)=0,~~~\rho^{B_2}(\widehat{e}_3,\widehat{e}_3)=\alpha\beta.
\end{align}
Then,
\begin{align}
&\widetilde{\rho}^{B_2}(\widehat{e}_1,\widehat{e}_1)=\widetilde{\rho}^{B_2}(\widehat{e}_1,\widehat{e}_2)=\widetilde{\rho}^{B_2}(\widehat{e}_1,\widehat{e}_3)=0,\nonumber\\
&\widetilde{\rho}^{B_2}(\widehat{e}_2,\widehat{e}_2)=\widetilde{\rho}^{B_2}(\widehat{e}_2,\widehat{e}_3)=0,~~~\widetilde{\rho}^{B_2}(\widehat{e}_3,\widehat{e}_3)=\alpha\beta.
\end{align}
By (2.7), we have
\begin{align}
&(L^{B_2}_Vg)(\widehat{e}_1,\widehat{e}_1)=0,~~~(L^{B_2}_Vg)(\widehat{e}_1,\widehat{e}_2)=-\mu_3\beta,~~~(L^{B_2}_Vg)(\widehat{e}_1,\widehat{e}_3)=\mu_2\gamma,\nonumber\\
&(L^{B_2}_Vg)(\widehat{e}_2,\widehat{e}_2)=(L^{B_2}_Vg)(\widehat{e}_2,\widehat{e}_3)=(L^{B_2}_Vg)(\widehat{e}_3,\widehat{e}_3)=0.
\end{align}
Then, if $(G_3,g,V)$  is an affine Ricci soliton associated to the Bott connection $\nabla^{B_2}$, by (2.8), we have the following four equations:
\begin{eqnarray}
       \begin{cases}
       \mu=0\\[2pt]
        \mu_3\beta=0 \\[2pt]
       \mu_2\gamma=0\\[2pt]
       \alpha\beta-\mu=0\\[2pt]
       \end{cases}
\end{eqnarray}
By solving (7.20), we get \\
\begin{thm}
$(G_3, g, V)$ is an affine Ricci soliton associated to the Bott connection $\nabla^{B_2}$ if and only if\\
\begin{eqnarray}
&&(1)\mu=\beta=\mu_2\gamma=0;\nonumber\\
&&(2)\mu=\alpha=\mu_3=\mu_2\gamma=0,~~~\beta\neq0.
\end{eqnarray}
\end{thm}
\vskip 0.5 true cm
\noindent{\bf 7.4 Affine Ricci solitons of $G_4$}\\
\vskip 0.5 true cm
\begin{lem} The Bott connection $\nabla^{B_2}$ of $G_4$ is given by
\begin{align}
&\nabla^{B_2}_{\widehat{e}_1}\widehat{e}_1=0,~~~\nabla^{B_2}_{\widehat{e}_1}\widehat{e}_2=-\widehat{e}_2+(2\eta-\beta)\widehat{e}_3,~~~\nabla^{B_2}_{\widehat{e}_1}\widehat{e}_3=-\beta\widehat{e}_2+\widehat{e}_3,\nonumber\\
&\nabla^{B_2}_{\widehat{e}_2}\widehat{e}_1=\nabla^{B_2}_{\widehat{e}_2}\widehat{e}_2=\nabla^{B_2}_{\widehat{e}_2}\widehat{e}_3=0,\nonumber\\
&\nabla^{B_2}_{\widehat{e}_3}\widehat{e}_1=\nabla^{B_2}_{\widehat{e}_3}\widehat{e}_2=\nabla^{B_2}_{\widehat{e}_3}\widehat{e}_3=0.
\end{align}
\end{lem}
\begin{lem} The curvature $R^{B_2}$ of the Bott connection $\nabla^{B_2}$ of $(G_4,g)$ is given by
\begin{align}
&R^{B_2}(\widehat{e}_1,\widehat{e}_2)\widehat{e}_1=R^{B_2}(\widehat{e}_1,\widehat{e}_2)\widehat{e}_2=R^{B_2}(\widehat{e}_1,\widehat{e}_2)\widehat{e}_3=0,\nonumber\\
&R^{B_2}(\widehat{e}_1,\widehat{e}_3)\widehat{e}_1=R^{B_2}(\widehat{e}_1,\widehat{e}_3)\widehat{e}_2=R^{B_2}(\widehat{e}_1,\widehat{e}_3)\widehat{e}_3=0,\nonumber\\
&R^{B_2}(\widehat{e}_2,\widehat{e}_3)\widehat{e}_1=0,~~~R^{B_2}(\widehat{e}_2,\widehat{e}_3)\widehat{e}_2=\alpha\widehat{e}_2+\alpha(\beta-2\eta)\widehat{e}_3,~~~R^{B_2}(\widehat{e}_2,\widehat{e}_3)\widehat{e}_3=\alpha\beta\widehat{e}_2-\alpha\widehat{e}_3.
\end{align}
\end{lem}
By (2.5), we have
\begin{align}
&\rho^{B_2}(\widehat{e}_1,\widehat{e}_1)=\rho^{B_2}(\widehat{e}_1,\widehat{e}_2)=\rho^{B_2}(\widehat{e}_1,\widehat{e}_3)=0,\nonumber\\
&\rho^{B_2}(\widehat{e}_2,\widehat{e}_1)=0,~~~\rho^{B_2}(\widehat{e}_2,\widehat{e}_2)=\alpha(2\eta-\beta),~~~\rho^{B_2}(\widehat{e}_2,\widehat{e}_3)=\alpha,\nonumber\\
&\rho^{B_2}(\widehat{e}_3,\widehat{e}_1)=0,~~~\rho^{B_2}(\widehat{e}_3,\widehat{e}_2)=\alpha,~~~\rho^{B_2}(\widehat{e}_3,\widehat{e}_3)=\alpha\beta.
\end{align}
Then,
\begin{align}
&\widetilde{\rho}^{B_2}(\widehat{e}_1,\widehat{e}_1)=\widetilde{\rho}^{B_2}(\widehat{e}_1,\widehat{e}_2)=\widetilde{\rho}^{B_2}(\widehat{e}_1,\widehat{e}_3)=0,\nonumber\\
&\widetilde{\rho}^{B_2}(\widehat{e}_2,\widehat{e}_2)=\alpha(2\eta-\beta),~~~\widetilde{\rho}^{B_2}(\widehat{e}_2,\widehat{e}_3)=\alpha,~~~\widetilde{\rho}^{B_2}(\widehat{e}_3,\widehat{e}_3)=\alpha\beta.
\end{align}
By (2.7), we have
\begin{align}
&(L^{B_2}_Vg)(\widehat{e}_1,\widehat{e}_1)=0,~~~(L^{B_2}_Vg)(\widehat{e}_1,\widehat{e}_2)=-\mu_2-\mu_3\beta,~~~(L^{B_2}_Vg)(\widehat{e}_1,\widehat{e}_3)=\mu_2(\beta-2\eta)-\mu_3,\nonumber\\
&(L^{B_2}_Vg)(\widehat{e}_2,\widehat{e}_2)=(L^{B_2}_Vg)(\widehat{e}_2,\widehat{e}_3)=(L^{B_2}_Vg)(\widehat{e}_3,\widehat{e}_3)=0.
\end{align}
Then, if $(G_4,g,V)$  is an affine Ricci soliton associated to the Bott connection $\nabla^{B_2}$, by (2.8), we have the following six equations:
\begin{eqnarray}
       \begin{cases}
       \mu=0\\[2pt]
       \mu_2+\mu_3\beta=0\\[2pt]
       \mu_2(\beta-2\eta)-\mu_3=0 \\[2pt]
      \alpha=0\\[2pt]
       \alpha(2\eta-\beta)+\mu=0\\[2pt]
              \alpha\beta-\mu=0\\[2pt]
       \end{cases}
\end{eqnarray}
By solving (7.26), we get \\
\begin{thm}
$(G_4, g, V)$ is an affine Ricci soliton associated to the Bott connection $\nabla^{B_2}$ if and only if
\begin{eqnarray*}
&&(1)\mu=\alpha=\mu_2+\mu_3\eta=0,~~~\beta=\eta;\nonumber\\
&&(2)\mu=\alpha=\mu_3=\mu_2=0,~~~\beta\neq\eta.
\end{eqnarray*}
\end{thm}

\vskip 0.5 true cm
\noindent{\bf 7.5 Affine Ricci solitons of $G_5$}\\
\vskip 0.5 true cm
\begin{lem} The Bott connection $\nabla^{B_2}$ of $G_5$ is given by
\begin{align}
&\nabla^{B_2}_{\widehat{e}_1}\widehat{e}_1=0,~~~\nabla^{B_2}_{\widehat{e}_1}\widehat{e}_2=0,~~~\nabla^{B_2}_{\widehat{e}_1}\widehat{e}_3=\beta\widehat{e}_2,\nonumber\\
&\nabla^{B_2}_{\widehat{e}_2}\widehat{e}_1=0,~~~\nabla^{B_2}_{\widehat{e}_2}\widehat{e}_2=\delta\widehat{e}_3,~~~\nabla^{B_2}_{\widehat{e}_2}\widehat{e}_3=\delta\widehat{e}_2,\nonumber\\
&\nabla^{B_2}_{\widehat{e}_3}\widehat{e}_1=-\alpha\widehat{e}_1,~~~\nabla^{B_2}_{\widehat{e}_3}\widehat{e}_2=0,~~~\nabla^{B_2}_{\widehat{e}_3}\widehat{e}_3=0.
\end{align}
\end{lem}
\begin{lem} The curvature $R^{B_2}$ of the Bott connection $\nabla^{B_2}$ of $(G_5,g)$ is given by
\begin{align}
&R^{B_2}(\widehat{e}_1,\widehat{e}_2)\widehat{e}_1=0,~~~R^{B_2}(\widehat{e}_1,\widehat{e}_2)\widehat{e}_2=\beta\delta\widehat{e}_2,~~~R^{B_2}(\widehat{e}_1,\widehat{e}_2)\widehat{e}_3=-\beta\delta\widehat{e}_3,\nonumber\\
&R^{B_2}(\widehat{e}_1,\widehat{e}_3)\widehat{e}_1=0,~~~R^{B_2}(\widehat{e}_1,\widehat{e}_3)\widehat{e}_2=-\beta\delta\widehat{e}_3,~~~R^{B_2}(\widehat{e}_1,\widehat{e}_3)\widehat{e}_3=-\beta(\alpha+\delta)\widehat{e}_2,\nonumber\\
&R^{B_2}(\widehat{e}_2,\widehat{e}_3)\widehat{e}_1=0,~~~R^{B_2}(\widehat{e}_2,\widehat{e}_3)\widehat{e}_2=-\delta^2\widehat{e}_3,~~~R^{B_2}(\widehat{e}_2,\widehat{e}_3)\widehat{e}_3=-(\beta\gamma+\delta^2)\widehat{e}_2.
\end{align}
\end{lem}
By (2.5), we have
\begin{align}
&\rho^{B_2}(\widehat{e}_1,\widehat{e}_1)=\rho^{B_2}(\widehat{e}_1,\widehat{e}_2)=\rho^{B_2}(\widehat{e}_1,\widehat{e}_3)=0,\nonumber\\
&\rho^{B_2}(\widehat{e}_2,\widehat{e}_1)=0,~~~\rho^{B_2}(\widehat{e}_2,\widehat{e}_2)=\delta^2,~~~\rho^{B_2}(\widehat{e}_2,\widehat{e}_3)=0,\nonumber\\
&\rho^{B_2}(\widehat{e}_3,\widehat{e}_1)=\rho^{B_2}(\widehat{e}_3,\widehat{e}_2)=0,~~~\rho^{B_2}(\widehat{e}_3,\widehat{e}_3)=-(\beta\gamma+\delta^2).
\end{align}
Then,
\begin{align}
&\widetilde{\rho}^{B_2}(\widehat{e}_1,\widehat{e}_1)=\widetilde{\rho}^{B_2}(\widehat{e}_1,\widehat{e}_2)=\widetilde{\rho}^{B_2}(\widehat{e}_1,\widehat{e}_3)=0,\nonumber\\
&\widetilde{\rho}^{B_2}(\widehat{e}_2,\widehat{e}_2)=\delta^2,~~~\widetilde{\rho}^{B_2}(\widehat{e}_2,\widehat{e}_3)=0,~~~\widetilde{\rho}^{B_2}(\widehat{e}_3,\widehat{e}_3)=-(\beta\gamma+\delta^2).
\end{align}
By (2.7), we have
\begin{align}
&(L^{B_2}_Vg)(\widehat{e}_1,\widehat{e}_1)=0,~~~(L^{B_2}_Vg)(\widehat{e}_1,\widehat{e}_2)=\mu_3\beta,~~~(L^{B_2}_Vg)(\widehat{e}_1,\widehat{e}_3)=-\mu_1\alpha,\nonumber\\
&(L^{B_2}_Vg)(\widehat{e}_2,\widehat{e}_2)=2\mu_3\delta,~~~(L^{B_2}_Vg)(\widehat{e}_2,\widehat{e}_3)=-\mu_2\delta,~~~(L^{B_2}_Vg)(\widehat{e}_3,\widehat{e}_3)=0.
\end{align}
Then, if $(G_5,g,V)$  is an affine Ricci soliton associated to the Bott connection $\nabla^{B_2}$, by (2.8), we have the following six equations:
\begin{eqnarray}
       \begin{cases}
       \mu=0\\[2pt]
       \mu_3\beta=0\\[2pt]
       \mu_1\alpha=0\\[2pt]
       \mu_3\delta+\delta^2+\mu=0\\[2pt]
       \mu_2\delta=0\\[2pt]
       \delta^2+\beta\gamma+\mu=0\\[2pt]
       \end{cases}
\end{eqnarray}
By solving (7.32), we get \\
\begin{thm}
$(G_5, g, V)$ is an affine Ricci soliton associated to the Bott connection $\nabla^{B_2}$ if and only if
$\mu=\delta=\gamma=\mu_1=\mu_3\beta=0,~~~ \alpha\neq 0.$
\end{thm}
\vskip 0.5 true cm
\noindent{\bf 7.6 Affine Ricci solitons of $G_6$}\\
\vskip 0.5 true cm
\begin{lem} The Bott connection $\nabla^{B_2}$ of $G_6$ is given by
\begin{align}
&\nabla^{B_2}_{\widehat{e}_1}\widehat{e}_1=0,~~~\nabla^{B_2}_{\widehat{e}_1}\widehat{e}_2=\alpha\widehat{e}_2+\beta\widehat{e}_3,~~~\nabla^{B_2}_{\widehat{e}_1}\widehat{e}_3=\gamma\widehat{e}_2+\delta\widehat{e}_3,\nonumber\\
&\nabla^{B_2}_{\widehat{e}_2}\widehat{e}_1=\nabla^{B_2}_{\widehat{e}_2}\widehat{e}_2=\nabla^{B_2}_{\widehat{e}_2}\widehat{e}_3=0,\nonumber\\
&\nabla^{B_2}_{\widehat{e}_3}\widehat{e}_1=\nabla^{B_2}_{\widehat{e}_3}\widehat{e}_2=\nabla^{B_2}_{\widehat{e}_3}\widehat{e}_3=0.
\end{align}
\end{lem}
\begin{lem} The curvature $R^{B_2}$ of the Bott connection $\nabla^{B_2}$ of $(G_6,g)$ is given by
\begin{align}
&R^{B_2}(\widehat{e}_s,\widehat{e}_t)\widehat{e}_p=0,
\end{align}
for any $(s,t,p)$.
\end{lem}
By (2.5), we have
\begin{align}
&\rho^{B_2}(\widehat{e}_s,\widehat{e}_t)=0,
\end{align}
for any pairs $(s,t)$.
Similarly,
\begin{align}
&\widetilde{\rho}^{B_2}(\widehat{e}_s,\widehat{e}_t)=0,
\end{align}
for any pairs $(s,t)$.
By (2.7), we have
\begin{align}
&(L^{B_2}_Vg)(\widehat{e}_1,\widehat{e}_1)=0,~~~(L^{B_2}_Vg)(\widehat{e}_1,\widehat{e}_2)=\mu_2\alpha+\mu_3\gamma,~~~(L^{B_2}_Vg)(\widehat{e}_1,\widehat{e}_3)=-\mu_2\beta-\mu_3\delta,\nonumber\\
&(L^{B_2}_Vg)(\widehat{e}_2,\widehat{e}_2)=(L^{B_2}_Vg)(\widehat{e}_2,\widehat{e}_3)=(L^{B_2}_Vg)(\widehat{e}_3,\widehat{e}_3)=0.
\end{align}
Then, if $(G_6,g,V)$  is an affine Ricci soliton associated to the Bott connection $\nabla^{B_2}$, by (2.8), we have the following three equations:
\begin{eqnarray}
       \begin{cases}
       \mu=0\\[2pt]
       \mu_2\alpha+\mu_3\gamma=0\\[2pt]
       \mu_2\beta+\mu_3\delta=0\\[2pt]
       \end{cases}
\end{eqnarray}
By solving (7.38), we get \\
\begin{thm}
$(G_6, g, V)$ is an affine Ricci soliton associated to the Bott connection $\nabla^{B_2}$ if and only if\\
\begin{eqnarray*}
&&(1)\mu=\mu_2=\mu_3=0,~~~\alpha\gamma-\beta\delta=0,~~~\alpha+\delta\neq0;\nonumber\\
&&(2)\mu=0,~~~\alpha=\beta,~~~\delta=\gamma,~~~\mu_2\neq0,~~~\mu_3\neq0;\nonumber\\
&&(3)\mu=\alpha+\beta=\delta+\gamma=0,~~~\mu_2\neq0,~~~\mu_3\neq0.
\end{eqnarray*}
\end{thm}
\vskip 0.5 true cm
\noindent{\bf 7.7 Affine Ricci solitons of $G_7$}\\
\vskip 0.5 true cm
\begin{lem} The Bott connection $\nabla^{B_2}$ of $G_7$ is given by
\begin{align}
&\nabla^{B_2}_{\widehat{e}_1}\widehat{e}_1=0,~~~\nabla^{B_2}_{\widehat{e}_1}\widehat{e}_2=-\beta\widehat{e}_2-\beta\widehat{e}_3,~~~\nabla^{B_2}_{\widehat{e}_1}\widehat{e}_3=\beta\widehat{e}_2+\beta\widehat{e}_3,\nonumber\\
&\nabla^{B_2}_{\widehat{e}_2}\widehat{e}_1=\alpha\widehat{e}_1,~~~\nabla^{B_2}_{\widehat{e}_2}\widehat{e}_2=\delta\widehat{e}_3,~~~\nabla^{B_2}_{\widehat{e}_2}\widehat{e}_3=\delta\widehat{e}_2,\nonumber\\
&\nabla^{B_2}_{\widehat{e}_3}\widehat{e}_1=-\alpha\widehat{e}_1,~~~\nabla^{B_2}_{\widehat{e}_3}\widehat{e}_2=-\delta\widehat{e}_3,~~~\nabla^{B_2}_{\widehat{e}_3}\widehat{e}_3=-\delta\widehat{e}_2.
\end{align}
\end{lem}
\begin{lem} The curvature $R^{B_2}$ of the Bott connection $\nabla^{B_2}$ of $(G_7,g)$ is given by
\begin{align}
&R^{B_2}(\widehat{e}_1,\widehat{e}_2)\widehat{e}_1=0,~~~R^{B_2}(\widehat{e}_1,\widehat{e}_2)\widehat{e}_2=\beta(2\delta-\alpha)(\widehat{e}_2+\widehat{e}_3),\nonumber\\
&R^{B_2}(\widehat{e}_1,\widehat{e}_2)\widehat{e}_3=\beta(\alpha-2\delta)(\widehat{e}_2+\widehat{e}_3),~~~R^{B_2}(\widehat{e}_1,\widehat{e}_3)\widehat{e}_1=0,\nonumber\\
&R^{B_2}(\widehat{e}_1,\widehat{e}_3)\widehat{e}_2=\beta(\alpha-2\delta)(\widehat{e}_2+\widehat{e}_3),~~~R^{B_2}(\widehat{e}_1,\widehat{e}_3)\widehat{e}_3=\beta(2\delta-\alpha)(\widehat{e}_2+\widehat{e}_3),\nonumber\\
&R^{B_2}(\widehat{e}_2,\widehat{e}_3)\widehat{e}_1=0,~~~R^{B_2}(\widehat{e}_2,\widehat{e}_3)\widehat{e}_2=\beta\gamma(\widehat{e}_2+\widehat{e}_3),\nonumber\\
&R^{B_2}(\widehat{e}_2,\widehat{e}_3)\widehat{e}_3=-\beta\gamma(\widehat{e}_2+\widehat{e}_3).
\end{align}
\end{lem}
By (2.5), we have
\begin{align}
&\rho^{B_2}(\widehat{e}_1,\widehat{e}_1)=\rho^{B_2}(\widehat{e}_1,\widehat{e}_2)=\rho^{B_2}(\widehat{e}_1,\widehat{e}_3)=0,\nonumber\\
&\rho^{B_2}(\widehat{e}_2,\widehat{e}_1)=0,~~~\rho^{B_2}(\widehat{e}_2,\widehat{e}_2)=-\beta\gamma,~~~\rho^{B_2}(\widehat{e}_2,\widehat{e}_3)=\beta\gamma,\nonumber\\
&\rho^{B_2}(\widehat{e}_3,\widehat{e}_1)=0,~~~\rho^{B_2}(\widehat{e}_3,\widehat{e}_2)=\beta\gamma,~~~\rho^{B_2}(\widehat{e}_3,\widehat{e}_3)=-\beta\gamma.
\end{align}
Then,
\begin{align}
&\widetilde{\rho}^{B_2}(\widehat{e}_1,\widehat{e}_1)=\widetilde{\rho}^{B_2}(\widehat{e}_1,\widehat{e}_2)=\widetilde{\rho}^{B_2}(\widehat{e}_1,\widehat{e}_3)=0,\nonumber\\
&\widetilde{\rho}^{B_2}(\widehat{e}_2,\widehat{e}_2)=-\beta\gamma,~~~\widetilde{\rho}^{B_2}(\widehat{e}_2,\widehat{e}_3)=\beta\gamma,~~~\widetilde{\rho}^{B_2}(\widehat{e}_3,\widehat{e}_3)=-\beta\gamma.
\end{align}
By (2.7), we have
\begin{align}
&(L^{B_2}_Vg)(\widehat{e}_1,\widehat{e}_1)=0,~~~(L^{B_2}_Vg)(\widehat{e}_1,\widehat{e}_2)=\mu_1\alpha-\mu_2\beta+\mu_3\beta,~~~(L^{B_2}_Vg)(\widehat{e}_1,\widehat{e}_3)=\mu_2\beta-\mu_1\alpha-\mu_3\beta,\nonumber\\
&(L^{B_2}_Vg)(\widehat{e}_2,\widehat{e}_2)=2\mu_3\delta,~~~(L^{B_2}_Vg)(\widehat{e}_2,\widehat{e}_3)=-(\mu_2+\mu_3)\delta,~~~(L^{B_2}_Vg)(\widehat{e}_3,\widehat{e}_3)=2\mu_2\delta.
\end{align}
Then, if $(G_7,g,V)$  is an affine Ricci soliton associated to the Bott connection $\nabla^{B_2}$, by (2.8), we have the following six equations:
\begin{eqnarray}
       \begin{cases}
       \mu=0\\[2pt]
       \mu_1\alpha-\mu_2\beta+\mu_3\beta=0\\[2pt]
      \mu_2\beta-\mu_3\beta-\mu_1\alpha=0\\[2pt]
       \mu_3\delta+\mu_2\delta=0\\[2pt]
       \mu_3\delta-\beta\gamma+\mu=0\\[2pt]
       \mu_3\delta+\mu_1\beta+\mu_2\delta-2\delta^2-\beta\gamma-\alpha\delta=0\\[2pt]
       \end{cases}
\end{eqnarray}
By solving (7.44), we get \\
\begin{thm}
$(G_7, g, V)$ is an affine Ricci soliton associated to the Bott connection $\nabla^{B_2}$ if and only if\\
\begin{eqnarray*}
&&(1)\mu=\alpha=\beta=\mu_2=\mu_3=0,~~~ \delta\neq 0;\nonumber\\
&&(2)\mu=\alpha=\gamma=\mu_2=\mu_3=0,~~~\delta\neq 0,~~~~\beta\neq 0;\nonumber\\
&&(3)\mu=\gamma=\delta=0,~~~\alpha\neq0,~~~\mu_2\beta-\mu_3\beta-\mu_1\alpha=0;\nonumber\\
&&(4)\mu=\gamma=\mu_1=\mu_2=\mu_3=0,~~~\alpha\neq0,~~~\delta\neq0.\nonumber\\
\end{eqnarray*}
\end{thm}
\indent Specially, let $V=0$, we get the following corollary:
\begin{cor}
(I)~~$(G_1, g, V)$ is an affine Einstein associated to the Bott connection $\nabla^{B_2}$ if and only if $\mu=\beta=0$;\nonumber\\
(II)~~$(G_2, g, V)$ is an affine Einstein associated to the Bott connection $\nabla^{B_2}$ if and only if $\mu=\alpha=0$;\nonumber\\
(III)~~$(G_3, g, V)$ is an affine Einstein associated to the Bott connection $\nabla^{B_2}$ if and only if
$\mu=\alpha\beta=0$;\nonumber\\
(IV)~~$(G_4, g, V)$ is an affine Einstein associated to the Bott connection $\nabla^{B_2}$ if and only if $\mu=\alpha=0$;\nonumber\\
(V)~~$(G_5, g, V)$ is an affine Einstein associated to the Bott connection $\nabla^{B_2}$ if and only if $\mu=\delta=\gamma=0,~~~\alpha\neq0$;\nonumber\\
(VI)~~$(G_6, g, V)$ is an affine Einstein associated to the Bott connection $\nabla^{B_2}$ if and only if $ \mu=0,~~~\alpha\gamma-\beta\delta=0,~~~\alpha+\delta\neq0$;\nonumber\\
(VII)~~$(G_7, g, V)$ is an affine Einstein associated to the Bott connection $\nabla^{B_2}$ if and only if $\mu=\gamma=\alpha+2\delta=0,~~~\alpha+\delta\neq0$.\nonumber\\
\end{cor}
\section{Affine Ricci Solitons associated to the perturbed Bott connection on three-dimensional Lorentzian Lie groups with the third distribution}
Similarly, by the above calculations, we always obtain $\mu=0$. In order to get the affine Ricci soliton with non zero $\mu$, we introduce the perturbed Bott connection $\widetilde{\nabla}^{B_2}$ in the following. Let $\widehat{e}_1^*$ be the dual base of $e_1$. We define on $G_{i=1,\cdot\cdot\cdot,7}$
\begin{align}
\widetilde{\nabla}^{B_2}_XY=\nabla^{B_2}_XY+a_0\widehat{e}_1^*(X)\widehat{e}_1^*(Y)e_1,
\end{align}
where $a_0$ is a non zero real number. Then
\begin{align}
\widetilde{\nabla}^{B_2}_{\widehat{e}_1}\widehat{e}_1=a_0\widehat{e}_1,~~~\widetilde{\nabla}^{B_2}_{\widehat{e}_s}\widehat{e}_t=\nabla^{B_2}_{\widehat{e}_s}\widehat{e}_t,
\end{align}
where $s$ and $t$ does not equal 1. We define
\begin{align}
(\widetilde{L}_V^{B_2}g)(X,Y):=g(\widetilde{\nabla}^{B_2}_XV,Y)+g(X,\widetilde{\nabla}^{B_2}_YV),
\end{align}
for vector fields $X,Y,V$. Then we have for $G_{i=1,\cdot\cdot\cdot,7}$
\begin{align}
(\widetilde{L}_V^{B_2}g)(\widehat{e}_1,\widehat{e}_1)=2a_0\mu_1,~~~(\widetilde{L}_V^{B_2}g)(\widehat{e}_s,\widehat{e}_t)=(L_V^{B_2}g)(\widehat{e}_s,\widehat{e}_t),
\end{align}
where $s$ and $t$ does not equal 1.
\begin{defn}$(G_i,V,g)$ is called the affine Ricci soliton associated to the connection $\widetilde{\nabla}^{B_2}$ if it satisfies
\begin{align}
(\widetilde{L}_V^{B_2}g)(X,Y)+2\widetilde{\overline{\rho}}^{B_2}(X,Y)+2\mu g(X,Y)=0.
\end{align}
\end{defn}
For $(G_1,\widetilde{\nabla}^{B_2})$, we have
\begin{align}
\widetilde{R}^{B_2}(\widehat{e}_1,\widehat{e}_3)\widehat{e}_1=\alpha(a_0-\beta) \widehat{e}_1,~~~\widetilde{R}^{B_2}(\widehat{e}_2,\widehat{e}_3)\widehat{e}_1=-a_0\beta\widehat{e}_1,~~~\widetilde{R}^{B_2}(\widehat{e}_s,\widehat{e}_t)\widehat{e}_p=R^{B_2}(\widehat{e}_s,\widehat{e}_t)\widehat{e}_p,
\end{align}
for $(s,t,p)\neq(1,3,1),(2,3,1)$.\\
\begin{align}
\widetilde{\overline{\rho}}^{B_2}(\widehat{e}_1,\widehat{e}_3)=\frac{\alpha(a_0-\beta)}{2},~~~\widetilde{\overline{\rho}}^{B_2}(\widehat{e}_s,\widehat{e}_t)=\widetilde{\rho}^{B_2}(\widehat{e}_s,\widehat{e}_t),
\end{align}
for the pair $(s,t)\neq(1,3)$. If $(G_1,g,V)$ is an affine Ricci soliton associated to the connection $\widetilde{\nabla}^{B_2}$, then by (8.5), we have
\begin{eqnarray}
       \begin{cases}
        a_0\mu_1+\mu=0 \\[2pt]
       \alpha\beta-\mu_1\alpha-\mu_3\beta=0\\[2pt]
       \mu_2\beta+\mu_1\alpha+\alpha(a_0-\beta)=0\\[2pt]
      \mu_3\alpha-\beta^2+\mu=0\\[2pt]
       \mu_3\alpha+\mu_2\alpha=0\\[2pt]
      \beta^2+\mu_2\alpha-\mu=0\\[2pt]
       \end{cases}
\end{eqnarray}
Solve (8.8), we get
\begin{thm}
$(G_1,V,g)$ is an affine Ricci soliton associated to the connection $\widetilde{\nabla}^{B_2}$ if and only if
$\mu_1=\beta-\frac{\beta(\beta^2-\mu)}{\alpha},~~~\mu_2=\frac{\mu-\beta^2}{\alpha},~~~\mu_3=\frac{\beta^2-\mu}{\alpha},~~~a_0=\frac{\beta(2\mu-\beta^2)}{\alpha^2},~~~a_0(\beta^3-\beta\gamma-\alpha\beta)-\alpha\mu=0$.
\end{thm}
For $(G_2,\widetilde{\nabla}^{B_2})$, we have
\begin{align}
\widetilde{R}^{B_2}(\widehat{e}_2,\widehat{e}_3)\widehat{e}_1=-a_0\alpha \widehat{e}_1,~~~\widetilde{R}^{B_2}(\widehat{e}_s,\widehat{e}_t)\widehat{e}_p=R^{B_2}(\widehat{e}_s,\widehat{e}_t)\widehat{e}_p,
\end{align}
for $(s,t,p)\neq(2,3,1)$.\\
\begin{align}
\widetilde{\overline{\rho}}^{B_2}(\widehat{e}_s,\widehat{e}_t)=\widetilde{\rho}^{B_2}(\widehat{e}_s,\widehat{e}_t),
\end{align}
for any pairs $(s,t)$. If $(G_2,g,V)$ is an affine Ricci soliton associated to the connection $\widetilde{\nabla}^{B_2}$, then by (8.5), we have
\begin{eqnarray}
       \begin{cases}
       a_0\mu_1+\mu=0  \\[2pt]
       \mu_2\gamma-\mu_3\beta=0 \\[2pt]
      \mu_2\beta+\mu_3\gamma=0\\[2pt]
      \alpha\beta-\mu=0\\[2pt]
       \alpha\gamma=0\\[2pt]
       \end{cases}
\end{eqnarray}
Solve (8.11), we get
\begin{thm}
$(G_2,V,g)$ is an affine Ricci soliton associated to the connection $\widetilde{\nabla}^{B_2}$ if and only if
$\mu=\alpha=\mu_1=\mu_2=\mu_3=0,~~~\gamma\neq0.$
\end{thm}
For $(G_3,\widetilde{\nabla}^{B_2})$, we have
\begin{align}
\widetilde{R}^{B_2}(\widehat{e}_2,\widehat{e}_3)\widehat{e}_1=-a_0\alpha \widehat{e}_1,~~~\widetilde{R}^{B_2}(\widehat{e}_s,\widehat{e}_t)\widehat{e}_p=R^{B_2}(\widehat{e}_s,\widehat{e}_t)\widehat{e}_p,
\end{align}
for $(s,t,p)\neq(2,3,1)$.\\
\begin{align}
\widetilde{\overline{\rho}}^{B_2}(\widehat{e}_s,\widehat{e}_t)=\widetilde{\rho}^{B_2}(\widehat{e}_s,\widehat{e}_t),
\end{align}
for any pairs $(s,t)$. If $(G_3,g,V)$ is an affine Ricci soliton associated to the connection $\widetilde{\nabla}^{B_2}$, then by (8.5), we have
\begin{eqnarray}
       \begin{cases}
       a_0\mu_1+\mu=0\\[2pt]
       \mu_3\beta=0 \\[2pt]
      \mu_2\gamma=0\\[2pt]
      \mu=0\\[2pt]
       \alpha\beta-\mu=0\\[2pt]
       \end{cases}
\end{eqnarray}
Solve (8.14), we get
\begin{thm}
$(G_3,V,g)$ is an affine Ricci soliton associated to the connection $\widetilde{\nabla}^{B_2}$ if and only if
\begin{eqnarray*}
&&(1)\mu=\mu_1=\beta=\mu_2\gamma=0;\nonumber\\
&&(2)\beta\neq 0,~~~\alpha=\mu=\mu_1=\mu_3=\mu_2\gamma=0.
\end{eqnarray*}
\end{thm}
For $(G_4,\widetilde{\nabla}^{B_2})$, we have
\begin{align}
\widetilde{R}^{B_2}(\widehat{e}_2,\widehat{e}_3)\widehat{e}_1=-a_0\alpha \widehat{e}_1,~~~\widetilde{R}^{B_2}(\widehat{e}_s,\widehat{e}_t)\widehat{e}_p=R^{B_2}(\widehat{e}_s,\widehat{e}_t)\widehat{e}_p,
\end{align}
for $(s,t,p)\neq(2,3,1)$.\\
\begin{align}
\widetilde{\overline{\rho}}^{B_2}(\widehat{e}_s,\widehat{e}_t)=\widetilde{\rho}^{B_2}(\widehat{e}_s,\widehat{e}_t),
\end{align}
for any pairs $(s,t)$. If $(G_4,g,V)$ is an affine Ricci soliton associated to the connection $\widetilde{\nabla}^{B_2}$, then by (8.5), we have
\begin{eqnarray}
       \begin{cases}
       a_0\mu_1+\mu=0\\[2pt]
       \mu_3\beta+\mu_2=0\\[2pt]
       \mu_2(\beta-2\eta)-\mu_3=0\\[2pt]
       \alpha(2\eta-\beta)+\mu=0\\[2pt]
      \alpha=0\\[2pt]
      \alpha\beta-\mu=0\\[2pt]
       \end{cases}
\end{eqnarray}
Solve (8.17), we get
\begin{thm}
$(G_4,V,g)$ is an affine Ricci soliton associated to the connection $\widetilde{\nabla}^{B_2}$ if and only if
\begin{eqnarray*}
&&(1)\alpha=\mu=\mu_1=\mu_3=\mu_2=0,~~~\beta\neq\eta;\nonumber\\
&&(2)\alpha=\mu=\mu_1=0,~~~\mu_2+\mu_3\eta=0,~~~\beta=\eta.
\end{eqnarray*}
\end{thm}

For $(G_5,\widetilde{\nabla}^{B_2})$, we have
\begin{align}
\widetilde{R}^{B_2}(\widehat{e}_1,\widehat{e}_3)\widehat{e}_1=-a_0\alpha \widehat{e}_1,~~~\widetilde{R}^{B_2}(\widehat{e}_2,\widehat{e}_3)\widehat{e}_1=-a_0\gamma\widehat{e}_1,~~~
\widetilde{R}^{B_2}(\widehat{e}_s,\widehat{e}_t)\widehat{e}_p=R^{B_2}(\widehat{e}_s,\widehat{e}_t)\widehat{e}_p,
\end{align}
for $(s,t,p)\neq(1,3,1),(2,3,1)$.\\
\begin{align}
\widetilde{\overline{\rho}}^{B_2}(\widehat{e}_2,\widehat{e}_1)=-\frac{a_0\alpha}{2},~~~
\widetilde{\overline{\rho}}^{B_2}(\widehat{e}_s,\widehat{e}_t)=\widetilde{\rho}^{B_2}(\widehat{e}_s,\widehat{e}_t),
\end{align}
for the pair $(s,t)\neq(1,3)$. If $(G_5,g,V)$ is an affine Ricci soliton associated to the connection $\widetilde{\nabla}^{B_2}$, then by (8.5), we have
\begin{eqnarray}
       \begin{cases}
       a_0\mu_1+\mu=0 \\[2pt]
      \mu_3\beta=0 \\[2pt]
      \mu_1\alpha+a_0\alpha=0\\[2pt]
      \mu_3\delta+\delta^2+\mu=0\\[2pt]
       \mu_2\delta=0\\[2pt]
      \delta^2+\beta\gamma+\mu=0\\[2pt]
       \end{cases}
\end{eqnarray}
Solve (8.20), we get
\begin{thm}
$(G_5,V,g)$ is an affine Ricci soliton associated to the connection $\widetilde{\nabla}^{B_2}$ if and only if
$\beta=\alpha=\mu_2=\mu_3=0,~~~\delta\neq0,~~~\mu=-\delta^2,~~~\mu_1=\frac{\delta^2}{a_0}.$
\end{thm}
For $(G_6,\widetilde{\nabla}^{B_2})$, we have
\begin{align}
\widetilde{R}^{B_2}(\widehat{e}_s,\widehat{e}_t)\widehat{e}_p=R^{B_2}(\widehat{e}_s,\widehat{e}_t)\widehat{e}_p,
\end{align}
for any $(s,t,p)$.\\
\begin{align}
\widetilde{\overline{\rho}}^{B_2}(\widehat{e}_s,\widehat{e}_t)=\widetilde{\rho}^{B_2}(\widehat{e}_s,\widehat{e}_t),
\end{align}
for any pairs $(s,t)$. If $(G_6,g,V)$ is an affine Ricci soliton associated to the connection $\widetilde{\nabla}^{B_2}$, then by (8.5), we have
\begin{eqnarray}
       \begin{cases}
       a_0\mu_1+\mu=0 \\[2pt]
      \mu_2\alpha+\mu_3\gamma=0\\[2pt]
       \mu_2\beta+\mu_3\delta=0\\[2pt]
       \mu=0\\[2pt]
       \end{cases}
\end{eqnarray}
Solve (8.23), we get
\begin{thm}
$(G_6,V,g)$ is an affine Ricci soliton associated to the connection $\widetilde{\nabla}^{B_2}$ if and only if
\begin{eqnarray*}
&&(1)\mu=\mu_1=\mu_2=\mu_3=0,~~~\alpha\gamma-\beta\delta=0,~~~\alpha+\delta\neq0;\nonumber\\
&&(2)\mu=\mu_1=0,~~~\alpha=\beta,~~~\delta=\gamma,~~~\mu_2\neq0,~~~\mu_3\neq0;\nonumber\\
&&(3)\mu=\mu_1=\alpha+\beta=\delta+\gamma=0,~~~\mu_2\neq0,~~~\mu_3\neq0.
\end{eqnarray*}
\end{thm}
For $(G_7,\widetilde{\nabla}^{B_2})$, we have
\begin{align}
&\widetilde{R}^{B_2}(\widehat{e}_1,\widehat{e}_2)\widehat{e}_1=a_0\alpha \widehat{e}_1,~~~\widetilde{R}^{B_2}(\widehat{e}_1,\widehat{e}_3)\widehat{e}_1=-a_0\alpha\widehat{e}_1,\nonumber\\
&\widetilde{R}^{B_2}(\widehat{e}_2,\widehat{e}_3)\widehat{e}_1=a_0\gamma\widehat{e}_1,~~~\widetilde{R}^{B_2}(\widehat{e}_s,\widehat{e}_t)\widehat{e}_p=R^{B_2}(\widehat{e}_s,\widehat{e}_t)\widehat{e}_p,
\end{align}
for $(s,t,p)\neq(1,2,1),(1,3,1),(2,3,1)$.\\
\begin{align}
\widetilde{\overline{\rho}}^{B_2}(\widehat{e}_1,\widehat{e}_2)=\frac{a_0\alpha}{2},~~~\widetilde{\overline{\rho}}^{B_2}(\widehat{e}_2,\widehat{e}_3)=-\frac{a_0\alpha}{2},~~~\widetilde{\overline{\rho}}^{B_2}(\widehat{e}_s,\widehat{e}_t)=\widetilde{\rho}^{B_2}(\widehat{e}_s,\widehat{e}_t),
\end{align}
for the pair $(s,t)\neq(1,2),(1,3)$. If $(G_7,g,V)$ is an affine Ricci soliton associated to the connection $\widetilde{\nabla}^{B_2}$, then by (8.5), we have
\begin{eqnarray}
       \begin{cases}
        a_0\mu_1+\mu=0 \\[2pt]
       \mu_2\beta-\mu_1\alpha-\mu_3\beta-a_0\alpha=0\\[2pt]
       \mu_3\delta-\beta\gamma+\mu=0\\[2pt]
       (\mu_2+\mu_3)\delta-\beta\gamma=0\\[2pt]
    \mu_2\delta-\beta\gamma-\mu=0\\[2pt]
       \end{cases}
\end{eqnarray}
Solve (8.26), we get
\begin{thm}
$(G_7,V,g)$ is an affine Ricci soliton associated to the connection $\widetilde{\nabla}^{B_2}$ if and only if
\begin{eqnarray*}
&&(1)\alpha=\beta=\mu_2+\mu_3=0,~~~\delta\neq0,~~~\mu=-\mu_3\delta,~~~\mu_1=\frac{\mu_3\delta}{a_0};\nonumber\\
&&(2)\alpha=\gamma=\mu=\mu_1=\mu_2=\mu_3=0,~~~\beta\neq0;\nonumber\\
&&(3)\alpha\neq0,~~~\delta=\gamma=\mu=\mu_1=0,~~~(\mu_2+\mu_3)\beta-a_0\alpha=0;\nonumber\\
&&(4)\alpha\neq0,~~~\delta\neq0,~~~\gamma=0,~~~\mu_1=-\frac{\mu}{a_0},~~~\mu_2=\frac{\mu}{\delta},~~~\mu_3=-\frac{\mu}{\delta},~~~2a_0\beta\mu+\mu\alpha\delta-a_0^2\alpha\delta=0.
\end{eqnarray*}
\end{thm}

\section*{Acknowledgements}
The second author was supported in part by  NSFC No.11771070. The authors are deeply grateful to the referees for their valuable comments and helpful suggestions.


\section*{References}
\begin{thebibliography}{00}
\bibitem{A} A. N. Siddiqui, B. Y. Chen, O. Bahadir, Statistical solitons and inequalities for statistical warped product submanifolds, Mathematics 7 (2019), 797.
\bibitem{F} F. Baudoin. Sub-Laplacians and hypoelliptic operators on totally geodesic Riemannian foliations.
In Geometry, analysis and dynamics on sub-Riemannian manifolds. Vol. 1, EMS Ser.
Lect. Math. pages 259-321. Eur. Math. Soc. Z$\ddot{u}$rich, 2016.
\bibitem{G} G. Calvaruso, Einstein-like metrics on three-dimensional homogeneous Lorentzian manifolds, Geom.
Dedicata 127 (2007) 99-119.
\bibitem{J} J. A. Alvarez L\'{o}pez and P. Tondeur. Hodge decomposition along the leaves of a Riemannian
foliation. J. Funct. Anal. 99(2):443-458, 1991.
\bibitem{L} L. A. Cordero, P. E. Parker, Left-invariant Lorentzian metrics on 3-dimensional Lie groups, Rend. Mat. Appl. 17 (1997), 129-155.
\bibitem{M} M. Crasmareanu, A new approach to gradient Ricci solitons and generalizations, Filomat 32 (2018), 3337-3346.
\bibitem{N} N. Halammanavar, K. Devasandra, Kenmotsu manifolds admitting Schouten-van Kampen connection, Facta Univ. Ser. Math. Inform. 34 (2019), 23-34.
\bibitem{RK} R. K. Hladky. Connections and curvature in sub-Riemannian geometry. Houston J. Math.,
38(4):1107-1134, 2012.
\bibitem{RS} R. S.Hamilton, The Ricci flowonsurfaces, in: J. A. Isenberg, (Ed.), Mathematics andgeneral relativity, Contemp. Math. 71,Amer. Math. Soc. Providence, RI, Santa Cruz, CA, 1988, pp. 237-262.
\bibitem{S1} S. Hui, R. Prasad, D. Chakraborty, Ricci solitons on Kenmotsu manifolds with respect to quarter symmetric non-metric $\phi$-connection, Ganita 67 (2017), 195-204.
\bibitem{S2} S. Y. Perktas, A. Yildiz, On quasi-Sasakian 3-manifolds with respect to the Schouten-van Kampen connection, Int. Elec. J. Geom. 13 (2020), 62-74.
\bibitem{W} W. Batat, K. Onda, Algebraic Ricci solitons of three-dimensional Lorentzian Lie groups, J. Geom.
Phys. 114 (2017), 138-152.
\bibitem{YH} Y. Han, A. De, P. Zhao, On a semi-quasi-Einstein manifold, J. Geom. Phys. 155 (2020), 103739.
\bibitem{YY1} Y. Wang, Curvature of multiply warped products with an affine connection, Bull. Korean Math. Soc. 50 (2013), 1567-1586.
\bibitem{YY2} Y. Wang, Affine Ricci Solitons of Three-Dimensional Lorentzian Lie Groups, Journal of Nonlinear Mathematical Physics. 2021.
\bibitem{YS} Y. Wang, S. Wei, Gauss-Bonnet theorems in the affine group and the group of rigid motions of the Minkowski plane, Sci. China Math. 2020.


\end{thebibliography}
\end{document}